\documentclass [10pt,a4paper]{article}
\usepackage{amssymb}
\usepackage{enumerate,verbatim}
\usepackage{amsmath}
\usepackage{amsfonts,mathrsfs}
\usepackage{amsthm,wasysym}
\usepackage{epsfig,lpic}
\usepackage[applemac]{inputenc}
\usepackage[english]{babel}

\usepackage{pdfsync}
\usepackage{xcolor,hyperref,refcheck}
\usepackage[T1]{fontenc}
\usepackage{mathtools}

\usepackage[width=1.0\textwidth]{caption}

\definecolor{green}{rgb}{0,0.5,0}

\hypersetup{backref,colorlinks=true,linkcolor=blue,citecolor=blue}

\usepackage{marvosym}

\newcounter{teoremaganso}
\newcounter{appendix}

\newcounter{coryganso}

\renewenvironment{abstract}{\small\quotation\noindent
 {\bfseries \abstractname .}}{\endquotation \par}

\newenvironment{trivlistparm}[2]{\trivlist
\item[{#1 #2}]}{\endtrivlist}

\newenvironment{prooftext}[1]{\trivlistparm{\bfseries}{#1}}{\Qed\endtrivlistparm}
\newenvironment{prova}{\trivlistparm{\bfseries}{Proof.}}{\Qed\endtrivlistparm}

\catcode`\@=11

\def\resetthefootnote{\renewcommand{\thefootnote}{\@arabic\c@footnote} }
\def\@principiremex#1{\trivlist
 \item[\hskip \labelsep{\bfseries #1\ \thetheo.}]\ignorespaces}
\def\opar@principiremex#1[#2]{\trivlist
 \item[\hskip \labelsep{\bfseries #1\ \thetheo\ (#2).}]\ignorespaces}

\newcommand{\newTHEOremrom}[2]{\newenvironment{#1}{\refstepcounter{theo}\@ifnextchar[{\opar@principiremex{#2}}
{\@principiremex{#2}}}{\qedB\endtrivlist}} \catcode`\@=12
\DeclareMathSymbol{\square}{\mathord}{AMSa}{"03}
\newcommand{\qedB}{\nopagebreak\hspace*{\fill}$\square$\par}
\newcommand{\Qed}{\nopagebreak\hspace*{\fill}{\vrule width6pt height6pt depth0pt}\par}

\newtheorem {theo} {Theorem} [section]
\newtheorem {prop} [theo] {Proposition}
\newtheorem {cory} [theo] {Corollary}
\newtheorem {lem} [theo] {Lemma}
\newtheorem {bigtheo} [teoremaganso] {Theorem}
\newtheorem {bigcory} [teoremaganso] {Corollary}

\newTHEOremrom {defi} {Definition}
\newTHEOremrom {obs} {Remark}
\newTHEOremrom {ex} {Example}

\newcommand{\refc}[1]{\mbox{$(\ref{#1})$}}
\newcommand{\secc}[1]{Section~\ref{#1}}
\newcommand{\teoc}[1]{Theorem~\ref{#1}}
\newcommand{\propc}[1]{Proposition~\ref{#1}}
\newcommand{\coryc}[1]{Corollary~\ref{#1}}
\newcommand{\lemc}[1]{Lemma~\ref{#1}}
\newcommand{\defic}[1]{Definition~\ref{#1}}
\newcommand{\obsc}[1]{Remark~\ref{#1}}
\newcommand{\exc}[1]{Example~\ref{#1}}
\newcommand{\figc}[1]{Figure~\ref{#1}}
\newcommand{\tabc}[1]{Table~\ref{#1}}

\newcommand{\N}{\ensuremath{\mathbb{N}}}
\newcommand{\Z}{\ensuremath{\mathbb{Z}}}
\newcommand{\R}{\ensuremath{\mathbb{R}}}

\newcommand{\T}{\boldsymbol{T}}
\newcommand{\Q}{\ensuremath{\mathbb{Q}}}
\newcommand{\C}{\ensuremath{\mathbb{C}}}

\newcommand{\F}{\ensuremath{\mathcal{F}}}

\newcommand{\np}{{\ensuremath{{\hat\mu}}}}
\newcommand{\n}{{\ensuremath{\ell}}}

\newcommand{\cc}{\ensuremath{\mathscr{C}}}
\newcommand{\E}{\ensuremath{\mathcal{E}}}

\def\map#1#2#3{\mbox{${#1}\!:{#2}\longrightarrow{#3}$}}

\newcommand{\op}{\ensuremath{\mbox{\rm o}}}

\newcommand{\mc}[1]{\mathcal{#1}}

\newcommand{\dsp}{\displaystyle}
\newcommand{\gorro}{\hat}

\title{\bf Asymptotic expansion of the Dulac map and time for \\ unfoldings of hyperbolic saddles: Coefficient properties
\footnotetext{This work has been partially funded by the Ministry of Science, Innovation and Universities of Spain through the grants PGC2018-095998-B-I00 and MTM2017-86795-C3-2-P and by the Agency for Management of University and Research Grants of Catalonia through the grants 2017SGR1725 and 2017SGR1617.
}}

\author{D. Mar\'{\i}n and J. Villadelprat
\\*[.1truecm]
{\small \textsl{Departament de Matem{\`a}tiques, Edifici Cc,
Universitat Aut{\`o}noma de Barcelona,}}\\*[-.05truecm]
{\small\textsl{08193 Cerdanyola del Vall\`es (Barcelona), Spain}}
\\*[-.05truecm]
{\small \textsl{Centre de Recerca Matem\`atica, Edifici Cc, Campus de Bellaterra,}}\\*[-.05truecm]
{\small \textsl{08193 Cerdanyola del Vall\`es (Barcelona), Spain}}
\\*[.1truecm]
{\small \textsl{Departament d'Enginyeria Inform{\`a}tica i Matem{\`a}tiques, ETSE,}}
\\*[-.05truecm]
{\small \textsl{Universitat Rovira i Virgili, 43007 Tarragona, Spain}}}
                  
\date{\today}

\begin{document}
\maketitle

\begin{abstract}
We consider a $\cc^\infty$ family of planar vector fields $\{X_{\hat\mu}\}_{\hat\mu\in\hat W}$ having a hyperbolic saddle and we study the Dulac map $D(s;\hat\mu)$ and the Dulac time $T(s;\hat\mu)$ from a transverse section at the stable separatrix to a transverse section at the unstable separatrix, both at arbitrary distance from the saddle. Since the hyperbolicity ratio $\lambda$ of the saddle plays an important role, we consider it as an independent parameter, so that $\hat\mu=(\lambda,\mu)\in \hat W=(0,+\infty)\times W$, where $W$ is an open subset of $\R^N.$ For each $\hat\mu_0\in\hat W$ and $L>0$, the functions $D(s;\hat\mu)$ and $T(s;\hat\mu)$ have an asymptotic expansion at $s=0$ and $\hat\mu\approx\hat\mu_0$ with the remainder being uniformly $L$-flat with respect to the parameters. The principal part of both asymptotic expansions is given in a monomial scale containing a deformation of the logarithm, the so-called Ecalle-Roussarie compensator. In this paper we are interested in the coefficients of these monomials, which are functions depending on $\hat\mu$ that can be shown to be $\cc^\infty$ in their respective domains and ``universally'' defined, meaning that their existence is stablished before fixing the flatness $L$ and the unfolded parameter $\hat\mu_0.$ Each coefficient has its own domain and it is of the form $((0,+\infty)\setminus D)\times W$, where~$D$ a discrete set of rational numbers at which a resonance of the hyperbolicity ratio $\lambda$ occurs. In our main result, \teoc{A}, we give the explicit expression of some of these coefficients and to this end a fundamental tool is the employment of a sort of incomplete Mellin transform. With regard to these coefficients we also prove that they have poles of order at most two at $D\times W$ and we give the corresponding residue, that plays an important role when compensators appear in the principal part. Furthermore we prove a result, \coryc{analitico}, showing that in the analytic setting each coefficient given in \teoc{A} is meromorphic on $(0,+\infty)\times W$ and has only poles, of order at most two, along $D\times W.$
\end{abstract}

\tableofcontents

\section{Introduction and statements of main results}

\begin{figure}[t]
   \centering
  \begin{lpic}[l(0mm),r(0mm),t(0mm),b(5mm)]{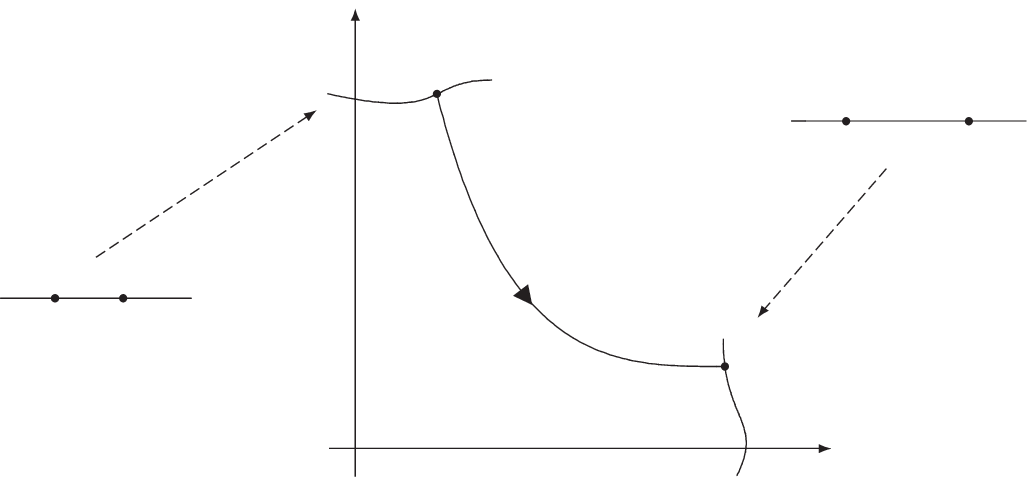}
   \lbl[l]{5,15.5;$0$}   
   \lbl[l]{12,16;$s$}   
   \lbl[l]{18,32;$\sigma_1$}
   \lbl[l]{28,40;$\Sigma_1$}      
   \lbl[l]{31.5,45;$x_2$}   
   \lbl[l]{38,42;$\sigma_1(s)$}   
   \lbl[l]{52,25;$\varphi(\,\cdot\,,\sigma_1(s))$}
   \lbl[l]{75.5,12;$\sigma_2(D(s))=\varphi(T(s),\sigma_1(s))$}  
   \lbl[l]{75,-2;$\Sigma_2$}   
   \lbl[l]{82,1;$x_1$}
   \lbl[l]{80,26;$\sigma_2$}
   \lbl[l]{85,33.5;$0$}
   \lbl[l]{96,33.5;$D(s)$}                   
   \end{lpic}
  \caption{Definition of $T(\,\cdot\,;\np)$ and $D(\,\cdot\,;\np)$, where $\varphi(t,p;\np)$ is the solution of~$X_\np$ passing through the point $p\in\mathscr U$ at time $t=0.$}
  \label{DefTyR}
 \end{figure}

In this paper we consider $\cc^\infty$ unfoldings of planar vector fields with a hyperbolic saddle. The study of the so-called Dulac map of the saddle has attracted the attention of many authors (see for instance \cite{DER,DRR1,DRR2,IlyYak,Mourtada,Roussarie89} and references there in) due, among other reasons, to its close connection with Hilbert's 16th problem (see \cite{Ily,Roussarie} for details). If $\np$ is the parameter unfolding, the \emph{Dulac map} $D(\,\cdot\,;\np)$ of the saddle is the transition map from a transverse section~$\Sigma_1$ at its stable separatrix $W_1$ to a transverse section~$\Sigma_2$ at its unstable separatrix $W_2$, whereas the \emph{Dulac time} $T(\,\cdot\,;\np)$ is the time that spends the flow to do this transition, see \figc{DefTyR}. In a recent paper \cite{MV20} we prove a general result for studying the asymptotic expansions of $D(s;\np)$ and $T(s;\np)$ at $s=0$, where $s$ is the variable parameterizing the transverse section $\Sigma_1$ and $s=0$ corresponds to the intersection point $W_1\cap\Sigma_1.$ In short, this general result gives a remainder that behaves well (i.e., uniformly on the parameters $\np)$ with respect to $\partial_s$ and provides a detailed description of the monomials appearing in the principal part. A key feature of this principal part is that the monomials can be ordered as $s\to 0^+$. This is a very important result for the theoretical point of view because it enables to bound the number of limit cycles or critical periodic orbits bifurcating from a polycycle. However there are specific problems where it is not only interesting to bound this number but also to determine from which parameters~$\np$ these bifurcations occur. Having explicit expressions of the coefficients of the monomials in the principal part is crucial for this purpose, see for instance \cite{ShaZeg94,Swy99} for limit cycles and \cite{MMV03,MMV2} for critical periodic orbits. The present paper is addressed to this issue. There are two features to be noted with regard to the hypothesis on the unfolding under consideration. On the one hand we suppose that the saddle is at the origin and, more significant, that the separatrices lay on the coordinate axis for all $\np$. It is important to point out that there is no loss of generality in assuming this since we prove in \cite[Lemma 4.3]{MV20} that there exists a smooth diffeomorphism, depending on the parameters, that straightens the two segments of the separatrices joining the points $W_1\cap\Sigma_1$ and $W_2\cap\Sigma_2$ with the saddle. That being said, we suppose on the other hand that the vector field has poles along the axis. The reason why we permit this ``polar'' factor is because, when dealing with polynomial vector fields, a special attention must be paid to the study of those polycycles with vertices at infinity in the Poincaré disc. The factor can come from the line at infinity in a saddle at infinity or, more generally, appear in a divisor after desingularizing a non-elementary singular point. We remark that (by means of a reparametrization of time) this factor can be neglected to study the Dulac map but, on the contrary, this cannot be done when dealing with the Dulac time.

The present paper is the continuation of \cite{MV19} and \cite{MV20} and concludes our contribution to the study of the theoretical aspects of the asymptotic expansion of the Dulac map and Dulac time of an unfolding of a hyperbolic saddle. Naturally the results that we shall obtain in this paper are strongly related with our previous ones. For reader's convenience we shall recall the essential results and definitions from \cite{MV19,MV20} in order to ease the legibility. 
Before that let us specify the hypothesis that we shall use throughout the paper. Setting $\np\!:=(\lambda,\mu)\in\hat W\!:=(0,+\infty)\times W$ with $W$ an open set of $\R^N,$ we consider the family of vector fields $\{X_{\np}\}_{\np\in\hat W}$ with
\begin{equation}\label{X}
 X_\np({x_1},{x_2})\!:=\frac{1}{x_1^{n_1}x_2^{n_2}}\Big({x_1}P_1({x_1},{x_2};\np)\partial_{x_1}+{x_2}P_2({x_1},{x_2};\np)\partial_{x_2}\Big),
\end{equation}
where
\begin{itemize}
\item $n\!:=(n_1,n_2)\in\Z^2_{\geq 0},$
\item $P_1$ and $P_2$ belong to $\mathscr C^{\infty}(\mathscr U\!\times\!\hat W)$ for some open set $\mathscr U$ of $\R^2$ containing the origin, 
\item $P_1({x_1},0;\np)>0$ and $P_2(0,{x_2};\np)<0$ for all $({x_1},0),(0,{x_2})\in\mathscr U$ and $\np\in \hat W,$
\item $\lambda=-\frac{P_2(0,0;\np)}{P_1(0,0;\np)}$.
\end{itemize}
Moreover, for $i=1,2,$ let \map{\sigma_i}{(-\varepsilon,\varepsilon)\times \hat W}{\Sigma_i} be a $\mathscr C^{\infty}$~transverse section to~$X_{\np}$ at $x_i=0$ defined by
 \[
  \sigma_i(s;\np)=\bigl(\sigma_{i1}(s;\np),\sigma_{i2}(s;\np)\bigr)
 \]
such that $\sigma_1(0,\np)\in\{(0,x_2);x_2>0\}$ and $\sigma_2(0,\np)\in\{(x_1,0);x_1>0\}$
for all $\np\in \hat W.$ We denote the Dulac map and Dulac time of~$X_\np$ from $\Sigma_1$ to $\Sigma_2$ by $D(\,\cdot\,;\np)$ and $T(\,\cdot\,;\np)$, respectively (see \figc{DefTyR}). 
Of course, in order that these functions are well defined for $s>0$ small enough, the open set $\mathscr U$ must contain the corner
\[
 \left\{(x_1,0);\,x_1\in\big[0,\sigma_{21}(0)\big]\right\}\cup\left\{(0,x_2);x_2\in \big[0,\sigma_{12}(0)\big]\right\}.
\] 
\begin{obs}\label{def_int}
For convenience, taking $\rho>0$ small enough, we define the open intervals
\[
 I_1\!:=\big(-\rho,\sigma_{12}(0)+\rho\big)\text{ and }I_2\!:=\big(-\rho,\sigma_{21}(0)+\rho\big)
\]
and assume in what follows that $\mathscr U$ contains $\big(\{0\}\!\times\! I_1\big)\cup\big(I_2\times\{0\}\big).$ Note then that, for $i=1,2$ and any $k\in\Z_{\geq 0}$, the map $(u,\np)\mapsto\partial_1^kP_i(0,u;\np)$ is $\cc^{\infty}$ on $I_1\!\times\!\hat W$ and the map $(u,\np)\mapsto\partial_2^kP_i(u,0;\np)$ is $\cc^{\infty}$ on $I_2\!\times\!\hat W$. Moreover $0\in I_i$ for $i=1,2.$ This technical observation will be important later on.
\end{obs}

\begin{defi}\label{defi_fun}
Consider $K\in\Z_{\geq 0}\cup\{+\infty\}$ and an open subset $U\subset\hat W\subset \R^{N+1}.$ We say that a function $\psi(s;\np)$ belongs to the class $\mathscr C^K_{s>0}(U)$, respectively $\E^K(U),$ if there 
exist an open neighbourhood $\Omega$ of 
\[ 
 \{(s,\np)\in\R^{N+2};s=0,\np\in U\}=\{0\}\times U
\]  
in $\R^{N+2}$ such that $(s,\np)\mapsto \psi(s;\np)$ is $\mathscr C^K$ on $\Omega\cap\big((0,+\infty)\times U\big),$ respectively $\Omega$. Finally we denote
         \[
           \E_+^K(U)\!:=\{\psi(s;\np)\in\E^K(U);\,\psi(0;\np)>0\text{ for all $\np\in U$}\}.
         \]
Here the letter $\E$ stands for functions in $\mathscr C^K_{s>0}(U)$ having \emph{extension} to $s=0.$            
\end{defi}

More formally, the definition of $\mathscr C^K_{s>0}(U)$ and $\E^K(U)$ must be thought in terms of germs with respect to relative neighborhoods of $\{0\}\times U$ in $(0,+\infty)\times U$. In doing so these sets become rings and we have the inclusions $\mathscr C^K(U)\subset\E^K(U)\subset\mathscr C^K_{s>0}(U)$. 

We can now introduce the notion of (finitely) flatness that we shall use in the sequel. 

\begin{defi}\label{defi2} 
Consider $K\in\Z_{\geq 0}\cup\{+\infty\}$ and an open subset $U\subset\hat W\subset\R^{N+1}.$ Given $L\in\R$ and $\np_0\in U$, we say that $\psi(s;\np)\in\mathscr C^K_{s>0}(U)$ is \emph{$(L,K)$-flat with respect to $s$ at $\np_0$}, and we write $\psi\in\F_L^K(\np_0)$, if for each $\nu=(\nu_0,\ldots,\nu_{N+1})\in\Z_{\geq 0}^{N+2}$  with $|\nu|=\nu_0+\cdots+\nu_{N+1}\leqslant K$ there exist a neighbourhood~$V$ of $\np_0$ and $C,s_0>0$ such that
\begin{equation*}
 \left|\frac{\partial^{|\nu|}\psi(s;\np)}{\partial s^{\nu_0}
 \partial\np_1^{\nu_1}\cdots\partial\np_{N+1}^{\nu_{N+1}}}\right|\leqslant C s^{L-\nu_0}
 \text{ for all $s\in(0,s_0)$ and $\np\in V$.} 
\end{equation*}
If $W$ is a (not necessarily open) subset of $U$ then define $\F_L^K(W)\!:=\bigcap_{\np_0\in W}\F_L^K(\np_0).$
\end{defi}
 
The principal part of the Dulac map and Dulac time will be expressed in terms of the following deformation of the logarithm. 

\begin{defi}%\label{defi_comp}
The function defined for $s>0$ and $\alpha\in\R$ by means of
 \[
  \omega(s;\alpha)\:=
  \left\{
   \begin{array}{ll}
    \frac{s^{-\alpha}-1}{\alpha} & \text{if $\alpha\neq 0,$}\\[2pt]
    -\log s & \text{if $\alpha=0,$}
   \end{array}
  \right.
 \]
is called the \emph{Ecalle-Roussarie compensator}. 
\end{defi}
\begin{figure}[t]
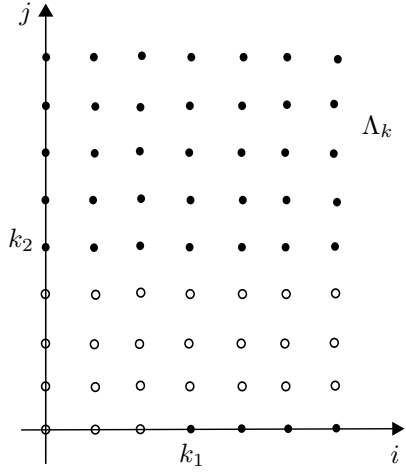

 \centering
 \begin{lpic}[l(0mm),r(0mm),t(0mm),b(5mm)]{lambda(0.75)}
   \lbl[l]{65,2;$i$}   
   \lbl[l]{60,60;$\Lambda_k$}   
   \lbl[l]{28,2;$k_1$}      
   \lbl[l]{0,80;$j$}
   \lbl[l]{-2,40;$k_2$}
 \end{lpic}
 \caption{The filled dots are points $(i,j)\in\Z_{\ge 0}^2$ in the set $\Lambda_k$ for $k=(k_1,k_2)$.}\label{Lambda}  
\end{figure}

\begin{defi}\label{alldefi}
Given any $k=(k_1,k_2)\in\Z^2_{\geq 0}$, throughout the paper we shall use the following notation:
\begin{itemize}
\item $\Lambda_{k}\!:=(\Z_{\ge k_1}\!\times\{0\})\cup(\Z_{\ge0}\!\times \Z_{\ge k_2})$, see \figc{Lambda}.
\item $D_{ij}^k\!:=\big\{\lambda>0:\text{there exits }(i',j')\in\Lambda_{k}\setminus\{(i,j)\}\text{ such that } i+\lambda j=i'+\lambda j'\big\}$.% for each $(i,j)\in\Lambda_k$.  
\item $\mathscr B_{\lambda,L}^k\!:=\big\{(i,j)\in\Lambda_k:i+\lambda j\leqslant L\big\}$ for each $L\in\R$ and $\lambda>0$.
\item $D_L^k\!:=\big\{\lambda>0:\text{there exits }(i,j)\in\mathscr B_{\lambda,L}^k\text{ such that }\lambda\in D_{ij}^k\big\}.$ 
\item For $\lambda=p/q\in\Q_{>0}$ with $\gcd(p,q)=1$ and $(i,j)\in\Lambda_k$,
\begin{equation*}
\mathscr A_{ij\lambda}^k\!:=\left\{
\begin{array}{cl}
 \emptyset & \text{ if $(i+rp,j-rq)\in\Lambda_k$ for some $r\in\N,$} \\[5pt]
 \left\{r\in\Z_{\geq 0}\,:\,(i-rp,j+rq)\in\Lambda_k\right\} & \text{ otherwise.}
\end{array}
\right.
\end{equation*}
\end{itemize}
Observe that if $k_2=0$ then $\Lambda_k=\Z^2_{\ge 0}=\Lambda_0$ regardless of the value of $k_1$. One can prove on the other hand, see \cite[Remark 3.3]{MV20}, that $D_{ij}^k$ and $D_L^k$ are discrete subsets of $\Q_{>0}$.
\end{defi} 
Let us point out that in the previous definition $k$ stands always for a two-dimensional vector with components in $\Z_{\ge 0}$. That being said, if $k=(0,0)$ then we write $\Lambda_{0}$, $D_{ij}^0$, $\mathscr B_{\lambda,L}^0$, $D_L^0$ and $\mathscr A_{ij\lambda}^0$ for shortness.

For the reader's convenience we merge Theorems~A and B of \cite{MV20} in the following result. In its statement we use the notation introduced so far and denote
\begin{equation*}
 T_0(\np)=\left\{\begin{array}{cl}0 & \text{if } n\neq (0,0),\\[3pt] 
 \frac{-1}{P(0,0;\np)} & \text{if } n=(0,0),\end{array}\right.
\end{equation*}
where recall that the components of $n=(n_1,n_2)\in\Z_{\geq 0}^2$ are the orders of the poles of $X_\np$ along the axis.

\begin{theo}\label{oldA} 
Let $D(s;\np)$ and $T(s;\np)$ be, respectively, the Dulac map and the Dulac time of the hyperbolic saddle \refc{X} from $\Sigma_1$ and $\Sigma_2$. 

\begin{enumerate}[$(a)$]

\item For each $(i,j)\in\Lambda_0$ there exists $\Delta_{ij}\in\mathscr C^\infty\big(((0,+\infty)\setminus D_{ij}^0)\times W\big)$ such that, for every $L>0$ and $\lambda_0>0,$ the following hold: 
\begin{enumerate}
\item[$(a1)$] If $\lambda_0\notin D^0_{L-\lambda_0}$ then
\[
D(s;\np)=s^\lambda\sum_{(i,j)\in\mathscr B_{\lambda_0,L-\lambda_0}^0}\Delta_{ij}(\np)s^{i+\lambda j}+\F_L^\infty(\{\lambda_0\}\times W). 
\]

\item[$(a2)$] If $\lambda_0\in D^0_{L-\lambda_0}$ then there exists a neighbourhood $\hat U$ of $\{\lambda_0\}\times W$ such that
\[
D(s;\np)=s^\lambda\sum_{(i,j)\in\mathscr B_{\lambda_0,L-\lambda_0}^0}\boldsymbol{\Delta}_{ij}^{\lambda_0}\big(\omega(s;\alpha);\np\big)s^{i+\lambda j}+\F_L^\infty(\{\lambda_0\}\times W),
\]
where $\lambda_0=p/q$ with $\gcd(p,q)=1$, $\alpha(\np)=p-\lambda q$ and $\boldsymbol{\Delta}_{ij}^{\lambda_0}(w;\np)\in\cc^\infty(\hat U)[w]$ with
\[
\boldsymbol{\Delta}_{ij}^{\lambda_0}(w;\np)=\sum_{r\in\mathscr A_{ij\lambda_0}^0}\Delta_{i-rp,j+rq}(\np)(1+\alpha w)^r\text{ for $\lambda\neq\lambda_0$.}
\]

\end{enumerate}

Moreover $\Delta_{00}(\np)>0$ for all $\np\in\hat W.$ 

\item For each $(i,j)\in\Lambda_n$ there exists $T_{ij}\in\mathscr C^\infty\big(((0,+\infty)\setminus D_{ij}^n)\times W\big)$ such that, for every $L>0$ and $\lambda_0>0,$ the following hold:
\begin{enumerate}
\item[$(b1)$] If $\lambda_0\notin D^n_L$ then
\[
T(s;\np)=T_0(\np)\log s+\sum_{(i,j)\in\mathscr B_{\lambda_0,L}^n}T_{ij}(\np)s^{i+\lambda j}+\F_L^\infty(\{\lambda_0\}\times W).
\]

\item[$(b2)$] If $\lambda_0\in D^n_L$ then there exists a neighbourhood $\hat U$ of $\{\lambda_0\}\times W$ such that
\[
T(s;\np)=T_0(\np)\log s+\sum_{(i,j)\in\mathscr B_{\lambda_0,L}^n}\T_{ij}^{\lambda_0}\big(\omega(s;\alpha);\np\big)s^{i+\lambda j}+\F_L^\infty(\{\lambda_0\}\times W),
\]
where $\lambda_0=p/q$ with $\gcd(p,q)=1$, $\alpha(\np)=p-\lambda q$ and $\T_{ij}^{\lambda_0}(w;\np)\in\cc^\infty(\hat U)[w]$ with
 \[
\T_{ij}^{\lambda_0}(w;\np)=\sum_{r\in\mathscr A_{ij\lambda_0}^n}T_{i-rp,j+rq}(\np)(1+\alpha w)^r\text{ for $\lambda\neq\lambda_0$.}
\]

\end{enumerate}

\end{enumerate}

\end{theo}

For every $(i,j)\in\Lambda_n$, \teoc{oldA} shows that $T_{ij}(\lambda,\mu)$ is $\cc^\infty$ on $((0,+\infty)\setminus D_{ij}^n)\times W$. We will prove, see \lemc{polos}, that for each $\lambda_0\in D_{ij}^n$ there exists $\ell\in\Z_{\geq 0}$ such that $\np\mapsto (\lambda-\lambda_0)^\ell T_{ij}(\np)$ extends $\cc^\infty$ to $\{\lambda_0\}\times W.$ Moreover the number~$\ell$, which depends on $(i,j)$, $\lambda_0$ and $n=(n_1,n_2)$, is bounded by $i+j.$ Hence, roughly speaking, the coefficient $T_{ij}(\lambda,\mu)$ has poles of order at most $i+j$ along $D_{ij}^n\times W.$ Likewise, by \lemc{polos} as well, it follows that $\Delta_{ij}(\lambda,\mu)$ has poles of order at most $i+j$ along $D_{ij}^0\times W.$ 

One of the main goals in this paper is to obtain explicit formulas for some of the coefficients $\Delta_{ij}$ 
of the Dulac map and some of the coefficients $T_{ij}$ of the Dulac time. More concretely, we will give the expressions of $\Delta_{ij}$ for $(i,j)\in\{(0,0),(0,1),(1,0),(1,1)\}$ and $T_{ij}$ for $(i,j)\in\{(n_1,0),(n_1+1,0),(0,n_2),(0,n_2+1)\}$.
This information is relevant because the corresponding monomials $s^{i+\lambda j}$ are the first (as $s\to 0^+$) that appear in the asymptotic expansion of the Dulac map (see \teoc{3punts}) and the Dulac time (see \teoc{9punts}).
With this aim in view we next define some functions that depend uniquely on $P_i(x_1,x_2;\np)$, for $i=1,2$, and $n=(n_1,n_2)$, see \refc{X}. The latter is fixed, whereas the dependence on $\np$ will be omitted for shortness.
\begin{equation}\label{def_fun}
\begin{array}{ll}
\dsp L_1(u)\!:=\exp\int_0^u\left(\frac{P_1(0,z)}{P_2(0,z)}+\frac{1}{\lambda}\right)\frac{dz}{z} & 
\dsp L_2(u)\!:=\exp\int_0^u\left(\frac{P_2(z,0)}{P_1(z,0)}+{\lambda}\right)\frac{dz}{z} \\[15pt]
\dsp M_1(u)\!:=L_1(u)\partial_1\!\left(\frac{P_1}{P_2}\right)(0,u)&  
\dsp M_2(u)\!:=L_2(u)\partial_2\!\left(\frac{P_2}{P_1}\right)(u,0)\\[15pt]
\dsp A_1(u)\!:=\frac{L_1^{n_1}(u)}{P_2(0,u)}& 
\dsp A_2(u)\!:=\frac{L_2^{n_2}(u)}{P_1(u,0)}\\[15pt]
\dsp B_1(u)\!:=n_1%\frac{L_1^{n_1}(u)}{P_2(0,u)}
A_1(u)\gorro{M}_1(1/\lambda,u)
& 
\dsp B_2(u)\!:=n_2%\frac{L_2^{n_2}(u)}{P_1(u,0)}
A_2(u)\gorro{M}_2(\lambda,u)
\\[3pt]
\dsp \hspace{2.8truecm} +L_1^{n_1+1}(u)\partial_1P_2^{-1}(0,u) & 
\dsp \hspace{2.8truecm}+L_2^{n_2+1}(u)\partial_2P_1^{-1}(u,0)
\\[15pt]
\dsp C_1(u)\!:=L_1^2(u)\partial_1^2P_2^{-1}(0,u)
&
\dsp C_2(u)\!:=L_2^2(u)\partial_2^2P_1^{-1}(u,0)
\\[3pt]
\dsp \hspace{2.4truecm}+2L_1(u)\gorro{M}_1(1/\lambda,u)\partial_1P_2^{-1}(0,u)&
\dsp \hspace{2.4truecm}+2L_2(u)\gorro{M}_2(\lambda,u)\partial_2P_1^{-1}(u,0)
\end{array}
\end{equation}
Here, given $\alpha\in\R\setminus\Z_{\ge 0}$ and a real valued function $f(x)$ that is $\cc^\infty$ in an open interval containing $x=0$, $\hat f(\alpha,x)$ is a sort of incomplete Mellin transform that we will introduce in Appendix~\ref{apMellin}. In this regard we point out, see \lemc{fun_ok}, that the functions $L_i(u),$ $M_i(u)$ and $A_i(u)$ are $\cc^{\infty}$ on an interval $I_i$ that contains $u=0$ for $i=1,2.$ On the other hand, for shortness as well, in the statement of our main result we use the compact notation $\sigma_{ijk}$ for the $k$th derivative at $s=0$ of the $j$th component of $\sigma_i(s;\np)$, i.e., 
\[
 \sigma_{ijk}(\np)\!:=\partial^k_s\sigma_{ij}(0;\np).
\]
In particular we consider the following real values (where once again we omit the dependence on $\np$):
\begin{equation}\label{def_S}
\begin{array}{l}
\dsp S_1\!:
=\frac{\sigma_{112}}{2\sigma_{111}}-\frac{\sigma_{121}}{\sigma_{120}}\left(\frac{P_1}{P_2}\right)\!(0,\sigma_{120})-\frac{\sigma_{111}}{L_1(\sigma_{120})}\gorro{M}_1(1/\lambda,\sigma_{120})\\[20pt]
\dsp S_2\!:
=\frac{\sigma_{222}}{2\sigma_{221}}-\frac{\sigma_{211}}{\sigma_{210}}\left(\frac{P_2}{P_1}\right)\!(\sigma_{210},0)-{\frac{\sigma_{221}}{L_2(\sigma_{210})}}\gorro{M}_2(\lambda,\sigma_{210}).
\end{array}
\end{equation}

We are now in position to state the main result of the present paper, \teoc{A}, which provides the explicit expression of the above-mentioned coefficients, see points $(b)$ and $(c)$. In addition to that we also establish in point $(a)$ a factorization property among the coefficients $\Delta_{ij}$ and $T_{ij}$ that holds for arbitrary $(i,j)$. This factorization is along the lines of the one given by Roussarie (see \cite[Theorem F]{Roussarie86} or \cite[\S 5.1.3]{Roussarie}) for the coefficients of the local Dulac map.

\begin{bigtheo}\label{A}
Assume $n\neq (0,0)$ and let $D(s;\np)$ and $T(s;\np)$ be, respectively, the Dulac map and the Dulac time of the hyperbolic saddle \refc{X} from $\Sigma_1$ and $\Sigma_2$. Consider moreover the coefficients $\Delta_{ij}$ and $T_{ij}$ given by \teoc{oldA}. Then the following assertions hold:
\begin{enumerate}[$(a)$]

\item\label{a} There exists a sequence $\{\Omega_{ij}\}_{(i,j)\in\Lambda_0}$ with $\Omega_{ij}\in\cc^{\infty}\big(((0,+\infty)\setminus D_{i0}^0)\times W\big)$ such that if $(i,j)\in\Lambda_0$ then
\begin{align*}
 \Delta_{ij}(\np)&=\Omega_{ij}(\np)\Delta_{0j}(\np)\text{ for all $\np\in\hat W$ with $\lambda\notin D_{ij}^0$,}
 \intertext{and if $(i,j)\in\Lambda_n$ with $j>0$ then}
 T_{ij}(\np)&=\Omega_{i,j-1}(\np)T_{0j}(\np)\text{ for all $\np\in\hat W$ with $\lambda\notin D_{ij}^n\cup D_{i0}^0\subset D_{ij}^0$.}
\end{align*}
%Moreover $D_{ij}^n\cup D_{i0}^0\subset D_{ij}^0$. 

\item\label{b} The coefficients $\Delta_{ij}$ for $(i,j)\in\{(0,0),(0,1),(1,0),(1,1)\}$ of the Dulac map are given by
\[
\Delta_{00}(\np)=\frac{\sigma_{111}^\lambda\sigma_{120}}{L_1^\lambda(\sigma_{120})}\frac{L_2(\sigma_{210})}{\sigma_{221}\sigma_{210}^\lambda},\quad
\Delta_{01}(\np)=-\Delta_{00}^2S_2,\quad 
\Delta_{10}(\np)=\Delta_{00}\lambda S_1\text{ and }
\Delta_{11}(\np)=-2\Delta_{00}^2\lambda S_1 S_2,
\]
where each equality is valid for all $\np\in\hat W$ with $\lambda\notin D_{ij}^0$. 
In particular, $\Omega_{10}(\np)=\lambda S_1$ and $\Omega_{11}(\np)=2\lambda S_1$.

\item\label{c} The coefficients $T_{ij}$ for  $(i,j)\in\{(n_1,0),(n_1+1,0),(0,n_2),(0,n_2+1)\}$ of the Dulac time are given by
\begin{align*}
&T_{n_1,0}(\np)
=-\frac{\sigma_{111}^{n_1}\sigma_{120}^{n_2}}{L_1^{n_1}(\sigma_{120})}\gorro{A}_1(n_1/\lambda-n_2,\sigma_{120}),\\[8pt]
&T_{0,n_2}(\np)
=\Delta_{00}^{n_2}\frac{\sigma_{210}^{n_1}\sigma_{221}^{n_2}}{L_2^{n_2}(\sigma_{210})}\gorro A_2(n_2\lambda-n_1,\sigma_{210}),\\[8pt]
&T_{n_1+1,0}(\np)
=-\sigma_{111}^{n_1}\sigma_{120}^{n_2}\left(\frac{\sigma_{121}}{\sigma_{120}P_2(0,\sigma_{120})}+\frac{n_1S_1}{L_1^{n_1}(\sigma_{120})}\gorro{A}_1(n_1/\lambda-n_2,\sigma_{120})\right.\\
&\hspace{5.75truecm}+\left.\frac{\sigma_{111}}{L_1^{n_1+1}(\sigma_{120})}\gorro{B}_1\big((n_1+1)/\lambda-n_2,\sigma_{120}\big)\right),\\[8pt]
&T_{0,n_2+1}(\np)
=\Delta_{00}^{n_2+1}\sigma_{210}^{n_1}\sigma_{221}^{n_2}\left(\frac{\sigma_{211}}{\sigma_{210}P_1(\sigma_{210},0)}+\frac{\sigma_{221}}{L_2^{n_2+1}(\sigma_{210})}\gorro B_2\big(\lambda(n_2+1)-n_1,\sigma_{210}\big)\right),
\end{align*}
where each equality is valid for all $\np\in\hat W$ with $\lambda\notin D_{ij}^n,$ except for the third one in which the values $\lambda=\frac{1}{k}$, $k=1,2,\ldots,\lceil\frac{n_2}{n_1+1}\rceil-1,$ must be excluded as well.
Moreover, if $n_1=0$ then
\begin{align*}
&T_{20}(\np)
=-\sigma_{120}^{n_2}\left(\frac{\sigma_{122}\sigma_{120}+(n_2-1)\sigma_{121}^2}{2\sigma_{120}^2P_2(0,\sigma_{120})}+\frac{\sigma_{121}^2}{2\sigma_{120}}\partial_2P_2^{-1}(0,\sigma_{120})
+\frac{\sigma_{121}\sigma_{111}}{\sigma_{120}}
\partial_1P_2^{-1}(0,\sigma_{120})
\right.\\
&\left.\hspace{4cm}+
\frac{\sigma_{111}^2}{2L_1^2(\sigma_{120})}\gorro{C}_1(2/\lambda-n_2,\sigma_{120})+\frac{\sigma_{111}S_1}{L_1(\sigma_{120})}\gorro{B}_1(1/\lambda-n_2,\sigma_{120})
\right),
\end{align*}
for all $\np\in\hat W$ with $\lambda\notin D_{20}^n\cup\left\{\frac{1}{k};\,k=1,2,\ldots,\lceil \frac{n_2}{2}\rceil-1\right\}.$
Finally if $n_2=0$ then
\begin{align*}
&T_{02}(\np)
=\Delta_{00}^2\sigma_{210}^{n_1}\left(
\frac{\sigma_{212}\sigma_{210}+(n_1-1)\sigma_{211}^2}{2\sigma_{210}^2P_1(\sigma_{210},0)}+\frac{\sigma_{211}^2}{2\sigma_{210}}\partial_1P_1^{-1}(\sigma_{210},0)+\frac{\sigma_{211}\sigma_{221}}{\sigma_{210}}\partial_2P_1^{-1}(\sigma_{210},0)\right.\\
&\hspace{5.6cm}\left.+\frac{\sigma_{221}^2}{2L_2^2(\sigma_{210})}\gorro{C}_2(2\lambda-n_1,\sigma_{210})-\frac{\sigma_{211}S_2}{2\sigma_{210}P_1(\sigma_{210},0)}\right)
\end{align*}
for all $\np\in\hat W$ with $\lambda\notin D_{02}^n$. 
\end{enumerate}
\end{bigtheo}

We point out that the coefficients $T_{ij}(\np)$ depend on $\np$ but also on $n=(n_1,n_2).$ We do not specify this dependence in the notation for the sake of shortness. This is the reason why, for instance, the expression for $T_{n_1+1,0}(\np)$ does not follow by replacing $n_1$ by $n_1+1$ in the expression for $T_{n_1,0}(\np)$.

The employment of the incomplete Mellin transform introduced in Appendix~\ref{apMellin} allows us to generalise and unify several formulas that we obtained previously in \cite{MMV03,MV06} under more restrictive hypothesis. With regard to the hypothesis, in those papers we restrict ourselves to the analytic setting (see \obsc{analytic_setting} below) and, more restraining, we assume that the family of vector fields $\{X_{\np}\}_{\np\in\hat W}$ in \refc{X} verifies the \emph{family linearization property} (FLP, for short), which means that $\{X_{\np}\}_{\np\in\hat W}$ is locally analytically equivalent to its linear part. In the present paper we do not require the FLP assumption and we consider the smooth setting instead of the analytic one. Furthermore the expressions for the coefficients that we obtain in those papers are only valid for hyperbolicity ratios varying in a specific range.  By using the properties of the incomplete Mellin transform proved in  \teoc{L8} we can get through this constrain as well. Let us exemplify this by noting that if $n_1=0$ and $n_2>0$ then
 \begin{align*}
   T_{0n_2}(\np)
&=\left(\frac{\sigma_{221}\Delta_{00}}{L_2(\sigma_{210})}\right)^{n_2}\gorro A_2(n_2\lambda,\sigma_{210})\\[5pt]
&=\left(\frac{\sigma_{221}\Delta_{00}}{L_2(\sigma_{210})}\right)^{n_2}\!\left(
 -\frac{A_2(0)}{n_2\lambda}+\sigma_{210}^{n_2\lambda}\int_0^{\sigma_{210}}
 \left(A_2(u)-A_2(0)\right)u^{-n_2\lambda}\frac{du}{u}
\right).
 \end{align*}
Here the first equality follows by $(c)$ in \teoc{A} (and it is valid for all $\lambda\notin D_{0n_2}^n=\frac{\N}{n_2} $, see \obsc{domains} below), whereas the second one follows by applying $(b)$ in \teoc{L8} with $k=1$ and assuming $n_2\lambda<1$ additionally. In \cite{MMV03} we study the case $\{n_1=0,n_2>0\}$ and the integral expression for $T_{0n_2}$ obtained after the second equality is precisely the one that we give in that paper, which only holds for $\lambda\in (0,\frac{1}{n_2})$ because the integrand has a pole of order $n_2\lambda+1$ at $u=0.$ Similarly, if $n_1=0$ and $n_2>0$ then
\begin{align*}
T_{10}(\np)
&=-\sigma_{120}^{n_2}\left(\frac{\sigma_{121}}{\sigma_{120}P_2(0,\sigma_{120})}
   +\frac{\sigma_{111}}{L_1(\sigma_{120})}\gorro{B}_1(1/\lambda-n_2,\sigma_{120})\right)\\[5pt]
&=-\sigma_{120}^{n_2}\left(\frac{\sigma_{121}}{\sigma_{120}P_2(0,\sigma_{120})}
   +\frac{\sigma_{111}}{L_1(\sigma_{120})}\sigma_{120}^{1/\lambda-n_2}\int_0^{\sigma_{120}}B_1(u)\,u^{n_2-1/\lambda}\frac{du}{u}\right). 
\end{align*}
In this case the first equality follows by $(c)$ in \teoc{A} (and it is valid as long as $\lambda\notin D_{10}^n=\frac{1}{\N_{\geq n_2}}$, see \obsc{domains} below) and the second one follows by applying $(b)$ in \teoc{L8} with $k=0$ provided that $1/\lambda-n_2<0$. The integral expression for $T_{10}$ obtained after the second equality is precisely the one that we give in \cite{MMV03}, which only converges for $\lambda\in (\frac{1}{n_2},+\infty).$ In \cite{MV06} we extend the results in \cite{MMV03} to arbitrary $n=(n_1,n_2)$ but still in the analytic setting and under the FLP assumption. The coefficient formulas given in that paper are also particular cases of the ones in \teoc{A}.

\begin{obs}\label{domains}
For the reader's convenience we specify the sets $D_{ij}^0$ and $D_{ij}^n$ corresponding to the coefficients in points $(b)$ and~$(c)$ in \teoc{A}. Taking \defic{alldefi} into account one can readily get that 
\[
D_{00}^0=\emptyset,\; D_{01}^0=\N,\; D_{10}^0=\frac{1}{\N}\text{ and }D_{11}^0=\N\cup\frac{1}{\N}
\]
for the coefficients of the Dulac map. Similarly, for the coefficients of the Dulac time, we have $D_{00}^n=\emptyset$,
\[
D_{n_1,0}^n=
\bigcup_{i=1}^{n_1}\frac{i}{\N_{\geq n_2}},
\; 
D_{0,n_2}^n=\left\{\begin{array}{cl}
\frac{\N_{\geq n_1}}{n_2} & \text{ if $n_2\geqslant 1$,}\\[5pt]
\emptyset & \text{ if $n_2=0$,}
\end{array}
\right.
\;
D_{n_1+1,0}^n=\bigcup_{i=1}^{n_1+1}\frac{i}{\N_{\ge n_2}}
\text{ and }
D_{0,n_2+1}^n=\frac{\N_{\ge n_1}}{n_2+1}\cup\N,
\]
together with $D_{20}^n=\frac{2}{\N_{\geq n_2}}$ for $n_1=0$ and 
$D_{02}^n=\frac{\N}{2}$ for $n_2=0$.
\end{obs}

As we already mentioned, by \lemc{polos} we know that the coefficients $\Delta_{ij}(\lambda,\mu)$ and $T_{ij}(\lambda,\mu)$ have poles of order at most $i+j$ along $\{\lambda_0\}\times W$ with $\lambda_0\in D_{ij}^0$ and along $\{\lambda_0\}\times W$ with $\lambda_0\in D_{ij}^n$, respectively. This general result will be proved in \secc{sec:poles}. In that section we sharpen this upper bound for the coefficients given in points $(b)$ and $(c)$ of \teoc{A} and we also compute the corresponding residues. This information is of relevance because these residues are the values at $\lambda_0$ of the leading coefficients of the polynomials $\boldsymbol{\Delta}_{ij}^{\lambda_0}(w;\np)$ and $\T_{ij}^{\lambda_0}(w;\np)$ in \teoc{3punts} and \teoc{9punts}, respectively. We illustrate this in \exc{ex1} for the Dulac map.

\begin{obs}\label{analytic_setting}
In this paper, foreseeing future applications, we will sometimes consider the analytic setting. By \emph{analytic setting} we mean that, for $i=1,2$, the function $P_i(x_1,x_2;\np)$ in \refc{X} is analytic on $V\times\hat W$ and that the parametrization $\sigma_i(s;\np)$ of the transverse section $\Sigma_i$ is analytic on $(-\varepsilon,\varepsilon)\times\hat W.$ Note in particular, see \obsc{def_int}, that $\partial_1^kP_i(0,u;\np)\in\cc^{\omega}(I_1\!\times\!\hat W)$ and $\partial_2^kP_i(u,0;\np)\in\cc^{\omega}(I_2\!\times\!\hat W)$ for $i=1,2$ and $k\in\Z_{\geq 0}$.
\end{obs}

In view of the above discussion about the poles of the coefficients, it is reasonable to expect that in the analytic setting the coefficients are meromorphic. In the present paper we are able to prove that this is the case for the coefficients considered in \teoc{A}. The following constitutes our second main result:

\begin{bigcory}\label{analitico}
In the analytic setting the following assertions hold:
\begin{enumerate}[$(a)$]
\item For each $(i,j)\in\{(0,0),(1,0),(0,1),(1,1)\}$, the coefficient $\Delta_{ij}$ of the Dulac map is meromorphic on 
$\hat W=((0,+\infty)\times W$ and has only poles, of order at most two, along $D_{ij}^0\times W.$

\item For each $(i,j)\in\{(0,0),(n_1,0),(0,n_2),(n_1+1,0),(0,n_2+1)\}$, the coefficient $T_{ij}$ of the Dulac time is meromorphic on $\hat W=((0,+\infty)\times W$ and has only poles, of order at most two, along $D_{ij}^n\times W.$ This is also the case for $(i,j)=(2,0)$ and $(i,j)=(0,2)$ assuming $n_1=0$ and $n_2=0$, respectively. 

\end{enumerate}
\end{bigcory}

Taking this partial result into account in the analytic setting we conjecture that for arbitrary $(i,j)$ the coefficient $\Delta_{ij}(\lambda,\mu)$ of the Dulac map is meromorphic on $(0,+\infty)\times W$ with poles along $\lambda\in D_{ij}^0$ and that the coefficient $T_{ij}(\lambda,\mu)$ of the Dulac time is meromorphic on $(0,+\infty)\times W$ with poles along $\lambda\in D_{ij}^n$.

The paper is organized in the following way. \secc{primera} is mainly devoted to prove \teoc{A}. Once this is done, and as an intermediate step towards the proof of \coryc{analitico}, at the end of \secc{primera} we show that, in the analytic setting, the coefficients $\Delta_{ij}$ and $T_{ij}$ listed in $(a)$ and $(b)$ of \teoc{A}, respectively, are analytic in their domains (see \propc{pre-analitico}). In \secc{sec:poles} we study the poles and residues of the coefficients. We begin by proving the above-mentioned \lemc{polos}, which constitutes a general result about the order of the poles. Next we prove a bunch of propositions that give the order of the pole and the respective residue for each coefficient listed in points $(a)$ and $(b)$ of \teoc{A}. Finally we conclude the section with the proof of \coryc{analitico}. \secc{quarta} aims at future applications of the tools developed so far. The main result of this paper, \teoc{A}, is intended to be applied in combination with \teoc{oldA}, that gathers our main results in \cite{MV20}. For this reason, and in order to ease the applicability, in \secc{quarta} we particularise \teoc{oldA} to specify the first monomials appearing in the asymptotic expansion of the Dulac map $D(s;\np)$, see \teoc{3punts}, and the Dulac time $T(s;\np)$, see \teoc{9punts}, for arbitrary hyperbolicity ratio~$\lambda_0$. By ``first monomials'' we mean as $s\to 0^+$, more concretely with respect to the strict partial order $\prec_{\lambda_0}$ introduced in \cite[Definition 1.7]{MV20}. It is here, dealing with a resonant hyperbolicity ratio $\lambda_0=p/q$, where the compensator $\omega(s;p-\lambda q)$ comes into play and the residues of the poles are needed, see \exc{ex1}.

\section{Proof of \teoc{A}}\label{primera}

For the reader's convenience we state first a result that we proved in a previous paper, see \cite[Corollary~2.2]{MV20}. In its statement we follow the notation introduced in Definitions~\ref{defi_fun} and \ref{defi2}.

\begin{lem}\label{new-factor1} 
Consider $f(s;\np)\in\E^K(U)$ with $K\in\N$
and any $m\in\N$ with $m\leqslant K$. Then the following hold:
\begin{enumerate}[$(a)$]
\item There exist $f_i(\np)\in\cc^{K-i}(U)$, $i=0,1,\ldots,m-1$, and $g(s;\np)\in\E^{K-m}(U)$ such that $$f(s;\np)=\sum_{i=0}^{m-1}f_i(\np)s^i+s^mg(s;\np).$$ 
\item For any $L\geqslant 0,$ $\E^{K}(U)\subset\cc^{K'}(U)[s]+\F_L^{K'}(U)$ provided that $K\geqslant K'+L$. 
\end{enumerate}
\end{lem}

The previous statement is aimed to study the flatness of the remainder in the asymptotic expansions that we shall deal with. The proof of $(a)$ shows in fact, see \cite{MV20}, that if $f\in\cc^K(I\times U)$ with~$I$ an open interval of $\R$ containing~$0$ then $g\in\cc^{K-m}(I\times U)$. We prove next that this result has its obvious analytic and smooth analogous. From now on, for simplicity in the exposition, we shall use $\varpi\in\{\infty,\omega\}$ as a wild card in $\cc^\varpi$ for the smooth class~$\cc^\infty$ and the analytic class $\cc^\omega$. 

\begin{lem}\label{new-factor2}
Let us consider an open interval $I$ of $\R$ containing $0$, an open subset $U$ of $\R^N$ and $m\in\N$.  
If  $f(s;\nu)\in\cc^\varpi(I\times U)$ with $\varpi\in\{\infty,\omega\}$ then there exists $g(s;\nu)\in\cc^\varpi(I\times U)$ such that 
\[
 f(s;\nu)=\sum_{i=0}^{m-1}\frac{\partial_s^if(0;\nu)}{i!}s^i+s^mg(s;\nu).
\]
\end{lem}

\begin{prova}
Given $\varpi\in\{\infty,\omega\}$, we claim that if $f(s;\nu)\in\cc^\varpi(I\times U)$ verifies $f(0;\nu)=0$ for all $\nu\in U$ then there exists $q(s;\nu)\in\cc^\varpi(I\times U)$ such that $f(s;\nu)=sq(s;\nu).$ In order to prove the claim note first that the existence of $q$ in a neighbourhood of any $(s_0,\nu_0)\in I\times U$ with $s_0\neq 0$ is clear. Moreover this function is uniquely defined on $(I\setminus\{0\})\times U.$ If $s_0=0$ then there exist~$\cc^\varpi$ functions $q(s;\nu)$ and $r(\nu)$ in a neighbourhood $V$ of $(0,\nu_0)$ in $\R^{N+1}$ such that $f(s;\nu)=sq(s;\nu)+r(\nu)$. Indeed, the case $\varpi=\omega$ follows by the Weierstrass Division Theorem (see \cite[Theorem 1.8]{Greuel} or \cite[Theorem 6.1.3]{Krantz}), whereas the case $\varpi=\infty$ is a consequence of the Malgrange Division Theorem (see \cite[Theorem 2]{Nirenberg} for instance). Furthermore, due to $r(\nu)=f(0;\nu)=0$, we get that $f(s;\nu)=sq(s;\nu)$. Hence for each $\nu_0\in U$ there exist a neighbourhood $V_{\nu_0}$ of $(0,\nu_0)$ in $\R^{N+1}$ and a function $q_{\nu_0}\in\cc^\varpi(V_{\nu_0})$ such that $f(s;\nu)=sq_{\nu_0}(s;\nu).$ Since $q_{\nu_0}(s;\nu)=\frac{f(s;\nu)}{s}$ for all $(s,\nu)\in V_{\nu_0}$ with $s\neq 0,$ we conclude that $q_{\nu_1}=q_{\nu_2}$ whenever $V_{\nu_1}\cap V_{\nu_2}\neq\emptyset$. This proves the claim.

The desired result follows from the claim by using induction on $m.$ More precisely, for the base case $m=1$ we apply the claim to $f(s;\nu)-f(0;\nu)$. For the inductive step we apply the claim to $g(s;\nu)-g(0;\nu)$, where $g$ is the remainder for the inductive hypothesis. In this way one can prove the existence of functions $f_i\in\cc^\varpi(U)$ and $g\in\cc^\varpi(I\times U)$ verifying that $f(s;\nu)=\sum_{i=0}^{m-1}f_i(\nu)s^i+s^mg(s;\nu).$
From here one can readily see that $f_i(\nu)=\frac{\partial_s^if(0;\nu)}{i!}$ and this completes the proof.
\end{prova}

In the next lemma we show that the regularity assumptions on the vector field \refc{X}, see Remarks~\ref{def_int} and~\ref{analytic_setting}, are transferred to the functions defined in \refc{def_fun}.

\begin{lem}\label{fun_ok}
Fix $\varpi\in\{\infty,\omega\}$ and let us assume the following:
\begin{enumerate}[$(a)$]
\item $P_1(u,0;\np)$ and $P_2(0,u;\np)$ are non-vanishing functions on 
         $I_2\!\times\!\hat W$ and $I_1\!\times\!\hat W$, respectively. 
%\item $P_i(0,u;\np)\in\cc^{\varpi}(I_1\!\times\!\hat W)$ and 
%         $P_i(u,0;\np)\in\cc^{\varpi}(I_2\!\times\!\hat W)$ for $i=1,2$.  
\item $\partial_1^kP_i(0,u;\np)\in\cc^{\varpi}(I_1\!\times\!\hat W)$ and 
         $\partial_2^kP_i(u,0;\np)\in\cc^{\varpi}(I_2\!\times\!\hat W)$ for $i=1,2$ and $k=0,1,2.$
\end{enumerate}
Then, for $i=1,2$, the functions $L_i(u;\np)$, $M_i(u;\np)$ and $A_i(u;\np)$ given in \refc{def_fun} are $\cc^\varpi$ on $I_i\times\hat W$. Moreover, 
\begin{enumerate}[1.]
\item the functions $B_1(u;\np)$ and $C_1(u;\np)$ are $\cc^\varpi$ on 
         $I_1\times((0,+\infty)\setminus\frac{1}{\N})\times W$, and
\item the functions $B_2(u;\np)$ and $C_2(u;\np)$ are $\cc^\varpi$ on $I_2\times((0,+\infty)\setminus\N)\times W$.
\end{enumerate}
\end{lem}

\begin{prova}
Since $\frac{P_2(0,0;\np)}{P_1(0,0;\np)}=-\lambda$ by definition, the application of \lemc{new-factor2} with $m=1$ implies that $L_i(u;\np)$ is $\cc^\varpi(I_i\times\hat W)$ for $i=1,2$. In its turn this shows that $A_i(u;\np)$ and $M_i(u;\np)$ are $\cc^\varpi(I_i\times\hat W)$ for $i=1,2$. Then, by \teoc{L8}, we can assert that $\hat M_i(\alpha,u;\np)$ is $\cc^\varpi$ on $(\R\setminus\Z_{\geq 0})\times I_i\times\hat W.$ More precisely, we use assertion~$(a)$ for the case $\varpi=\infty$ and assertion~$(d)$ for the $\varpi=\omega.$ This easily implies, see \refc{def_fun}, that the assertions 1 and 2 in the statement are true and completes the proof of the result.
\end{prova}

All the assertions except the last one in the next result are proved in \cite[Lemma A.2]{MV20}. The last one follows as a particular case of assertion $(c)$ in \cite[Lemma A.3]{MV20}.

\begin{lem}\label{FLK} 
Let $U$ and $U'$ be open sets of $\R^N$ and $\R^{N'}$ respectively and consider $W\subset U$ and $W'\subset U'.$ 
Then the following holds:
\begin{enumerate}[$(a)$]
\item $\F_L^K(W)\subset\F_L^K(\hat W)$ for any $\hat W\subset W$ and $\bigcap_n\F_L^K(W_n)=\F_L^K\left(\bigcup_n W_n\right)$.
\item $\F_L^K(W)\subset\F_L^K(W\times W')$.
\item $\mathscr C^K(U)\subset\E^K(U)\subset\F_0^K(W)$.
\item If $K\geqslant K'$ and $L\geqslant L'$ then $\F_L^K(W)\subset\F_{L'}^{K'}(W)$.
\item $\F_L^K(W)$ is closed under addition.
\item If $f\in\F_L^K(W)$ and $\nu\in\Z_{\ge0}^{N+1}$ with $|\nu|\leqslant K$ then 
         $\partial^\nu f\in\F_{L-\nu_0}^{K-|\nu|}(W)$.
\item $\F_L^K(W)\cdot\F_{L'}^K(W)\subset\F_{L+L'}^K(W)$.
\item Assume that \map{\phi}{U'}{U} is a $\mathscr C^K$ function with $\phi(W')\subset W$ and let us take 
        $g\in\F_{L'}^K(W')$ with $L'>0$ and verifying $g(s;\eta)>0$ for all $\eta\in W'$ and $s>0$ small enough. 
        Consider also any $f\in\F_L^K(W)$. Then $h(s;\eta)\!:=f(g(s;\eta);\phi(\eta))$ is a well-defined function 
        that belongs to $\F_{LL'}^K(W')$. 
\item If $\alpha\in\cc^K(U)$ then $s^\alpha\in\F_L^K(\{\nu\in U:\alpha(\nu)>L\})$.        
\end{enumerate}
\end{lem}

By applying the previous lemmas we can now prove the following:

\begin{lem}\label{composition}
Let $V$ an open set of $\R^N$ and consider a polynomial $Q(\,\cdot\,;\nu)$ with coefficients in $\cc^K(V)$ such that $Q(0;\nu)>0$ for all $\nu\in V.$ Let us also take $L>0$ and $L'\geqslant 1$ together with $\alpha\in\cc^K(V)$ such that $\alpha(\nu)>0$ for all $\nu\in V$. Then the following holds:
\begin{enumerate}[$(a)$]
\item $\left(sQ(s)+\F_{L+1}^K(V)\right)^\alpha\subset s^\alpha Q^\alpha(s)+\F_{L}^K(V)$, and 
\item $\F_{L'}^K(V)\circ \left(s^\alpha Q(s)+\F_L^K(V)\right)\subset \F_L^K\left(\left\{\nu\in V:\alpha(\nu)>L/L'\right\}\right).$
\end{enumerate}
\end{lem}

\begin{prova} 
In order to prove $(a)$ note first that 
\begin{equation}\label{comp_eq1}
 (sQ(s)+\F_{L+1}^K(V))^\alpha\subset s^\alpha(Q(s)+\F_{L}^K(V))^\alpha\subset s^\alpha Q^\alpha(s)(1+\F_{L}^K(V))^\alpha.
\end{equation}
Indeed, this follows by using twice $(g)$ in \lemc{FLK}. More concretely, in the first equality together with the fact that $1/s\in\F_{-1}^K(V)$, whereas in the second one noting also that $1/Q(s)\in\mc E^K(V)\subset\F_0^K(V)$. On the other hand, by using Lemmas~\ref{new-factor1} and \ref{FLK}, 
\[
g(x)\!:=(1+x)^\alpha-1\in s\E^\infty(V)\subset \F_1^\infty(V)\F_0^\infty(V)\subset \F_1^\infty(V),
\]
Thus $g\circ\F_{L}^K(V)\in\F_{L}^K(V)$ by $(h)$ in \lemc{FLK} and, therefore, $(1+\F_{L}^K(V))^\alpha\subset 1+\F_{L}^K(V)$. Taking this into account, the assertion in $(a)$ follows from~\refc{comp_eq1} noting that $s^\alpha Q^\alpha(s)\F_{L}^K(V)\subset\F_0^K(V)\F_{L}^K(V)\subset\F_{L}^K(V)$ due to $s^\alpha\in\F_0^K(V)$ by $(i)$ in \lemc{FLK}.

Let us turn next to the assertion in $(b)$. To this end note that $s^\alpha Q(s)\in\F_{L/L'}^K(V\cap\{\alpha>L/L'\})$ by~$(i)$ in \lemc{FLK}. On the other hand, due to $L'>1,$ $\F_{L}^K(V)\subset \F_{L/L'}^K(V)\subset\F_{L/L'}^K(V\cap\{\alpha>L/L'\})$ by~$(d)$ and $(a)$ in \lemc{FLK}. Thus, by $(e)$ in \lemc{FLK}, 
\[
s^\alpha Q(s)+\F_L^K(V)\subset\F_{L/L'}^K(V\cap\{\alpha>L/L'\}).
\]
On account of this and that, by $(a)$ in \lemc{FLK} again, $\F_{L'}^K(V)\subset\F_{L'}^K(V\cap\{\alpha>L/L'\})$, the application of $(h)$ in \lemc{FLK} shows that
\begin{align*}
\F_{L'}^K(V)\circ \left(s^\alpha Q(s)+\F_L^K(V)\right)&\subset
\F_{L'}^K(V\cap\{\alpha>L/L'\})\circ \F_{L/L'}^K(V\cap\{\alpha>L/L'\})\\
&\subset \F_L^K(V\cap\{\alpha>L/L'\}).
\end{align*}
This completes the proof of the result.
\end{prova}

We only need one more technical result in order to tackle the proof of \teoc{A}. It will be a consequence of the following easy observation.

\begin{obs}\label{rm1}
If $\sum_{i=1}^ma_ix^{\lambda_i}+\psi(x)=0$ for all $x\in (0,\varepsilon),$ where $\lambda_i\in\R$ with $\lambda_1<\lambda_2<\cdots<\lambda_m$, $a_1,a_2,\ldots,a_m\in\R$ and $\psi(x)=\mathrm{o}(x^{\lambda_m})$ then $a_1=a_2=\cdots=a_m=0$.
\end{obs}

\begin{lem}\label{combinacio}
Consider $\alpha,\beta\in\R\setminus\Z$ with $\alpha-\beta\notin\Z$ and two functions $f$ and $g$ that are $\cc^K$ on the interval $(-\delta,\delta)$ with $K>-\min(\alpha,\beta).$ If there exists $c\in\R$ satisfying that $x^\alpha f(x)+x^\beta g(x)=c$ for all $x\in (0,\delta)$ then $c=0.$
\end{lem}

\begin{prova}
Suppose that $\alpha<\beta$ and $n\!:=\min\{i\in\Z_{\geq 0}:\alpha+i>0\}.$ Hence $K\geqslant n$ and by applying Taylor's theorem we can write 
\[
 f(x)=a_0+a_1x+\ldots+a_nx^n+x^{n}R_1(x)\text{ and }
g(x)=b_0+b_1x+\ldots+b_nx^n+x^{n}R_2(x),
\]
with $\lim_{x\to 0}R_i(x)=0.$ Let us also set $\kappa\!:=\min\{i\in\Z_{\geq 0}:\beta+i>\alpha+n\}.$ Note then that $\kappa\in \{0,1\ldots,n\}.$ If we define $\psi(x)\!:=(b_{\kappa}x^\kappa+b_{\kappa+1}x^{\kappa+1}+\ldots+b_nx^n)x^\beta+x^{n}(x^\alpha R_1(x)+x^\beta R_2(x))$ then, on account of the assumption $x^\alpha f(x)+x^\beta g(x)=c,$ we get that
\[
 -cx^0+a_0x^\alpha+a_1x^{\alpha+1}+\ldots+a_nx^{\alpha+n}+b_0x^\beta+b_1x^{\beta+1}+\ldots+b_{\kappa-1}x^{\beta+\kappa-1}+\psi(x)=0
\]
for all $x\in (0,\delta).$ Taking the definition of $n$ and $\kappa$ into account, note that $\psi(s)=\op(x^{0})$, $\psi(s)=\op(x^{\alpha+n})$ and $\psi(s)=\op(x^{\beta+\kappa-1})$. Moreover all the exponents in $x^0,x^\alpha,x^{\alpha+1},\ldots,x^{\alpha+n},x^\beta,x^{\beta+1},\ldots,x^{\beta+\kappa-1}$ are different by the hypothesis on~$\alpha$ and~$\beta$, so that they can be ordered. Thus, on account of \obsc{rm1}, we can assert that all their coefficients are equal to zero, in particular $c=0.$ 
\end{prova}

\begin{prooftext}{Proof of \teoc{A}.}
Note first that by \teoc{oldA} we have two well defined sequences $\{\Delta_{ij}\}_{(i,j)\in\Lambda_0}$ and $\{T_{ij}\}_{(i,j)\in\Lambda_n}$ with $\Delta_{ij}\in\mathscr C^\infty\big(((0,+\infty)\setminus D_{ij}^0)\times W\big)$
and $T_{ij}\in\mathscr C^\infty\big(((0,+\infty)\setminus D_{ij}^n)\times W\big)$ where, by applying \cite[Lemma~3.2]{MV20}, $D_{ij}^0$ and $D_{ij}^n$ are discrete sets of rational numbers in $(0,+\infty).$ In order to prove the assertions in \refc{a}, for each $(i,j)\in\Lambda_0$ and $\np\in ((0,+\infty)\setminus D_{i0}^0)\times W$ we define $\Omega_{ij}(\np)$ by means of 
\begin{equation}\label{Aeq33}
 \left(1+\sum_{i=1}^\infty\frac{\Delta_{i0}(\np)}{\Delta_{00}(\np)}s^i\right)^{j+1}=\sum_{i=0}^\infty\Omega_{ij}(\np)s^i,
\end{equation}
where the equality must be thought in the ring of formal power series in $s$. 
Hence $\Omega_{ij}\in\Q\left[\frac{\Delta_{10}}{\Delta_{00}},\frac{\Delta_{20}}{\Delta_{00}},\ldots,\frac{\Delta_{i0}}{\Delta_{00}}\right]$ for each fixed $(i,j)\in\Lambda_0$. One can verify, see \defic{alldefi}, that $D_{i0}^0=\bigcup_{\ell=1}^i\frac{\ell}{\N}$ and thus $\cup_{k=1}^iD_{k0}^0=D_{i0}^0.$ Consequently, since $\Delta_{00}>0$ on $\hat W$ by $(a)$ in \teoc{oldA}, we can assert that $$\Omega_{ij}\in\mathscr C^\infty\big(((0,+\infty)\setminus D_{i0}^0)\times W\big).$$
That being said, our first goal is to prove that if $(i,j)\in\Lambda_0$ then
\begin{align}\label{Aeq30}
 &\Delta_{ij}(\np)-\Omega_{ij}(\np)\Delta_{0j}(\np)=0\text{ for all $\np\in\hat W$ with $\lambda\notin D_{ij}^0$,}\\
 \intertext{and that if $(i,j)\in\Lambda_n$ with $j>0$ then}\label{Aeq31}
 &T_{ij}(\np)-\Omega_{i,j-1}(\np)T_{0j}(\np)=0\text{ for all $\np\in\hat W$ with $\lambda\notin D_{ij}^n\cup D_{i0}^0$.}
\end{align}
To this aim let us note that the function on the left hand side of the equality in \refc{Aeq30}, respectively \refc{Aeq31}, is $\mathscr C^\infty$ in a neighbourhood of any $\np_\star=(\lambda_\star,\mu_\star)\in (0,+\infty)\times W$ with $\lambda_\star$ outside the discrete set $D_{ij}^0\cup D_{i0}^0\cup D_{0j}^0$, respectively $D_{ij}^n\cup D_{i0}^0\cup D_{0j}^n$. In this regard observe that $D_{ij}^n\subset D_{ij}^0,$ see \defic{alldefi}. It is also easy to show that, for any given any $k\in\Z_{\geq 0}^2,$ we have $D_{i0}^k\subset D_{ij}^k$ and $D_{0j}^k\subset D_{ij}^k.$ Consequently
\[
 D_{ij}^0\cup D_{i0}^0\cup D_{0j}^0=D_{ij}^0\text{ and }D_{ij}^n\cup D_{i0}^0\cup D_{0j}^n=D_{ij}^n\cup D_{i0}^0\subset D_{ij}^0,
\]
so that the function in \refc{Aeq30} is continuous on $((0,+\infty)\setminus D_{ij}^0)\times W$ whereas the function in \refc{Aeq31} is continuous on $((0,+\infty)\setminus (D_{ij}^n\cup D_{i0}^0))\times W$. Since $D_{ij}^0$ and $D_{ij}^n\cup D_{i0}^0$ are discrete sets of rational number in $(0,+\infty)$, it is clear that both identities will follow by continuity once we prove it for any $\np=(\lambda,\mu)\in\hat W$ with $\lambda\notin\Q.$ 

%More concretely, for any given $k=(k_1,k_2)\in\Z_{\geq 0}^2,$ it is not difficult to show that 
%\[
% D_{ij}^k=\frac{\N+\max(k_1-i-1,0)}{j}\cup
% \left(\;\bigcup_{\ell=1}^{j-k_2}\frac{\N}{\ell}\;\right)\cup
% \left(\;\bigcup_{\ell=1}^i\frac{\ell}{\N}\;\right)
%\]
%for all $(i,j)\in\Lambda_k$ with $j>0$, whereas $D_{i0}^k=\bigcup_{\ell=1}^i\frac{\ell}{\N_{\geq k_2}}$ for all $i\geqslant k_1.$ Here we follow the usual convention that a finite union is the empty-set when its index set is empty, so that in particular
% \[
%  %D_{i0}^k=\bigcup_{\ell=1}^i\frac{\ell}{\N_{\geq k_2}}\text{ and }
%  D_{0j}^k=\frac{\N_{\geq k_1}}{j}\cup\left(\;\bigcup_{\ell=1}^{j-k_2}\frac{\N}{\ell}\;\right).
% \]
%Consequently $D_{i0}^k,D_{0j}^k\subset D_{ij}^k.$\footnote{Debería suceder que $D_{i0}^k\cup D_{0j}^k=D_{ij}^k$.} 

The strategy to prove the identities in \refc{b} and \refc{c} will be the same. Indeed, let us write them as 
\begin{equation*}%\label{Aeq1}
\Delta_{ij}(\lambda,\mu)=\tilde\Delta_{ij}(\lambda,\mu)\text{ and }T_{ij}(\lambda,\mu)=\tilde T_{ij}(\lambda,\mu),%\text{ where $\np=(\lambda,\mu)$,}
\end{equation*}
i.e., $\tilde\Delta_{ij}$ and $\tilde T_{ij}$ are the functions on the right hand side of the equalities in the statement we want to prove. As we already mentioned, we know that 
\[
\Delta_{ij}\in\mathscr C^\infty\big(((0,+\infty)\setminus D_{ij}^0)\times W\big)\text{ and }T_{ij}\in\mathscr C^\infty\big(((0,+\infty)\setminus D_{ij}^n)\times W\big) 
\]
by \teoc{oldA}. On the other hand it turns out that there exist $\tilde D_{ij}^0,\tilde D_{ij}^n\subset\Q_{> 0}$ such that 
\[
 \tilde \Delta_{ij}\in\mathscr C^\infty\big(((0,+\infty)\setminus\tilde D_{ij}^0)\times W\big)\text{ and }
 \tilde T_{ij}\in\mathscr C^\infty\big(((0,+\infty)\setminus \tilde D_{ij}^n)\times W\big).
\]
The sets $\tilde D_{ij}^0$ and $\tilde D_{ij}^n$ will be given explicitly later on but at this moment the relevant property is that they are discrete in $(0,+\infty)$ as well. That said, for simplicity in the exposition, let us explain how the proof goes for the identity $T_{0,n_2}(\lambda,\mu)=\tilde T_{0,n_2}(\lambda,\mu)$. Thus, since 
$D_{0,n_2}^n\cup\tilde D_{0,n_2}^n$ is a discrete set of rational numbers in $(0,+\infty)$, for any given $\lambda_\star\notin D_{0,n_2}^n\cup\tilde D_{0,n_2}^n$ there exists a sequence of irrational numbers $(\lambda_k)_{k\in\N}$ such that $\lim_{k\to\infty}\lambda_k=\lambda_\star$. Hence, if we take any $\mu\in W$ then, by continuity, $\lim_{k\to\infty}T_{0,n_2}(\lambda_k,\mu)=T_{0,n_2}(\lambda_\star,\mu)$ and $\lim_{k\to\infty}\tilde T_{0,n_2}(\lambda_k,\mu)=\tilde T_{0,n_2}(\lambda_\star,\mu).$ So it is clear that the validity of the equality $T_{0,n_2}(\lambda,\mu)=\tilde T_{0,n_2}(\lambda,\mu)$ at any $\lambda=\lambda_\star$ which is not inside $D_{0,n_2}^n\cup\tilde D_{0,n_2}^n$ will follow once we prove it for any $\np=(\lambda,\mu)\in\hat W$ with $\lambda\notin\Q.$  This will be precisely our goal to prove each one of the equalities in the statement. As a matter of fact we will show that each equality is true in a neighbourhood of any $\np_0=(\lambda_0,\mu_0)\in\hat W$ with $\lambda_0\notin\Q.$

In addition to the identities in \refc{b} and \refc{c} we shall prove the equality in \refc{Aeq30} for $(i,j)=(i_1,j_1)$ and the equality in \refc{Aeq31} for $(i,j)=(i_2,j_2)$, where $(i_1,j_1)\in\Lambda_0$ and $(i_2,j_2)\in\Lambda_n$ are arbitrary but fixed. To this end, in view of the previous considerations, we fix any $\np_0=(\lambda_0,\mu_0)\in\hat W$ with $\lambda_0\notin\Q.$ Then by \cite[Theorem~A]{MMV08} we know that for each ${K}\in\N$ there exists a $\cc^{K}$ diffeomorphism
$$\Phi(u_1,u_2,\np)=\big(u_1\psi_1(u_1,u_2;\np),u_2\psi_2(u_1,u_2;\np),\np\big),$$ 
defined in an open set $U\times V$ with $(0,0)\in U\subset\R^2$ and $\np_0\in V\subset\hat W$, verifying
\begin{equation}\label{Aeq3}
 \Phi^*X_\np=\frac{P_1(0,0;\np)}{u_1^{n_1}u_2^{n_2}}(u_1\partial_{u_1}-\lambda u_2\partial_{u_2})
\end{equation}
and such that $\psi_i(0,0;\np)=1$, $i=1,2.$ Let us point out that in the forthcoming analysis it will be crucial that $K$ is larger than some fixed quantity $\mathcal N=\mathcal N(\lambda_0,n_1,n_2,i_1,i_2,j_1,j_2).$ We will specify at each step of the proof which is the necessary lower bound for $K$ and, at the end, $\mathcal N$ will be the maximum of them. This provides us with a specific value for $\mathcal N$ (that is not relevant at all) and in what follows we simply suppose that we take a $\cc^{K}$ normalising diffeomorphism $\Phi$ with $K\geqslant\mathcal N.$

For convenience we assume, without lost of generality, that 
\[
 U=\{(u_1,u_2)\in\R^2: |u_1|<\delta\text{ and }|u_2|<\delta\}=(-\delta,\delta)^2
\]
for some $\delta>0$ small enough such that, see \obsc{def_int}, $\Phi\big((-\delta,\delta)^2\times V)\subset(-\rho,\rho)^2\times V.$ Taking $\varepsilon_1,\varepsilon_2\in (0,\delta)$ we consider auxiliary $\cc^{K}$ transverse sections $\Sigma_1^\ell$ and $\Sigma_2^\ell$ to $x_1=0$ and $x_2=0$, see \figc{descomposicio}, parametrized~by
\begin{equation}\label{Aeq2}
 \tau_1(s;\varepsilon_1,\np)\!:=\Phi(s,\varepsilon_1;\np)\text{ and }\tau_2(s;\varepsilon_2,\np)\!:=\Phi(\varepsilon_2,s;\np),
\end{equation}
respectively. From now on, in addition to $\np$, we will also consider $\varepsilon\!:=(\varepsilon_1,\varepsilon_2)$ as parameter. In this respect we remark that $\tau_i(s;\varepsilon_i,\np)$ is a $\cc^{K}$ function on $U\times V$ for $i=1,2.$ Similarly as we did with $\sigma_i$, we denote
 \[
 \tau_{ijk}(\varepsilon_i,\np)\!:=\partial_s^k\tau_{ij}(0;\varepsilon_i,\np)
 \]
and we will write $\tau_{ijk}$ for the sake of shortness. %Notice that then $\tau_{12k}=O(\varepsilon_1)$ and $\tau_{21k}=O(\varepsilon_2)$. 
\begin{figure}[t]
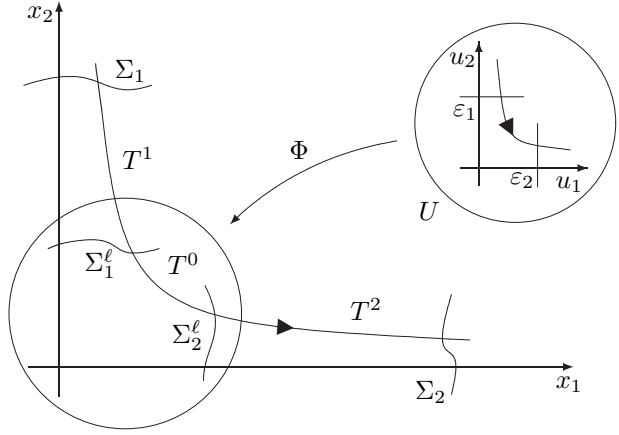

 \centering
 \begin{lpic}[l(0mm),r(0mm),t(0mm),b(5mm)]{descomposicio2} 
   \lbl[l]{2.5,55;$x_2$}
   \lbl[l]{14,47.5;$\Sigma_1$}
   \lbl[l]{10,22;$\Sigma_1^\ell$}
   \lbl[l]{21.5,12.5;$\Sigma_2^\ell$}
   \lbl[l]{37,37;$\Phi$}
   \lbl[l]{54,29;$U$}
   \lbl[l]{53.5,5;$\Sigma_2$}
   \lbl[l]{58,49;$u_2$}
   \lbl[l]{58.5,42;$\varepsilon_1$}
   \lbl[l]{66,33;$\varepsilon_2$}
   \lbl[l]{72,6;$x_1$}
   \lbl[l]{72,32.5;$u_1$}
   \lbl[l]{21,22;$T^0$}
   \lbl[l]{15,36;$T^1$} 
   \lbl[l]{45,16;$T^2$}    
 \end{lpic}
 \caption{Auxiliary transverse sections in the decomposition of $T.$}\label{descomposicio}
\end{figure}

The idea now is to decompose the Dulac map $D(s;\np)$ and the Dulac time $T(s;\np)$ as 
\begin{equation}\label{Aeq7}
D(s)=R_2(D_0(R_1(s)))\text{ and } T(s)=T^1(s)+T^0(R_1(s))+T^2(D_0(R_1(s))).
\end{equation}
Here $R_1(\,\cdot\,;\varepsilon_1,\np)$, $D_0(\,\cdot\,;\varepsilon,\np)$ and $R_2(\,\cdot\,;\varepsilon_2,\np)$ are, respectively, the transitions maps from $\Sigma_1$ to $\Sigma_1^\ell$, from 
$\Sigma_1^\ell$ to $\Sigma_2^\ell$, and from $\Sigma_2^\ell$ to $\Sigma_2$, whereas $T^1(\,\cdot\,;\varepsilon_1,\np)$, $T^0(\,\cdot\,;\varepsilon,\np)$ and $T^2(\,\cdot\,;\varepsilon_2,\np)$ are, respectively, the time that spends the flow to do this transition. It is well known that $D_0$ and $T^0$ are singular at $s=0$, whereas the other ones are regular. We study the latter by 
applying the results obtained in Appendix~\ref{ap_regular} and to this end, see \refc{ap1}, we rewrite the given vector field as
\[
X_\np=\frac{1}{x_1^{n_1}x_2^{n_2}}\left(x_1P_1(x_1,x_2)\partial_{x_1}+x_2P_2(x_1,x_2)\partial_{x_2}\right)
=\frac{1}{x_{i_2}^{n_{i_2}}f_{i_2}(x_{i_1},x_{i_2})}\left(\partial_{x_{i_1}}+h_{i_2}(x_{i_1},x_{i_2})x_{i_2}\partial_{x_{i_2}}\right)
\]
where $(i_1,i_2)\in\{(2,1),(1,2)\}$ and
\begin{equation}\label{Aeq0}
 \begin{array}{lll}
  \dsp f_1(u,v)=\frac{u^{n_2-1}}{P_2(v,u)} & \quad & \dsp h_1(u,v)=\frac{P_1(v,u)}{uP_2(v,u)}
  \\[10pt]
  \dsp f_2(u,v)=\frac{u^{n_1-1}}{P_1(u,v)}  & & \dsp  h_2(u,v)=\frac{P_2(u,v)}{uP_1(u,v)}
 \end{array}
\end{equation}
(At this point, and in what follows, we omit the dependence on the parameters for the sake of shortness when there is no risk of ambiguity. Moreover all though the proof the scripts $1$ and $2$ refer, respectively, to the first and second regular passage.) Setting $I\!:=(0,\delta),$ we apply (twice) \lemc{L3} with $\nu=(\varepsilon_i,\np)\in I\times V$ for $i=1,2.$ In doing so, and taking \lemc{new-factor1} also into account, we can assert that
\begin{equation}\label{Aeq9}
 R_i(s;\varepsilon_i,\np)=\sum_{k=1}^{L_i}R_{ik}(\varepsilon_i,\np)s^k+\F_{{L_i}+1}^0(I\!\times\!V)\text{ and }T^i(s;\varepsilon_i,\np)=\sum_{k=n_i}^{L_i}T^i_{k}(\varepsilon_i,\np)s^k+\F_{{L_i}+1}^0(I\!\times\! V),
\end{equation}
with $R_{ik},T^i_{k}\in\cc^0(I\times V)$ provided that ${K}\geqslant {L_i}+1$ for $i=1,2$. We know furthermore that $R_{i1}>0$. Turning to the assumption $K\geqslant\mathcal N$, let us advance that we will also require that $L_i\geqslant\mathcal N$ for $i=1,2$, which is neither a problem because, as we explained before, $\mathcal N=\mathcal N(\lambda_0,n_1,n_2,i_1,i_2,j_1,j_2)$ and we can take~$K$ large enough from the very beginning.

With regard to the passage from $\Sigma_1^\ell$ to $\Sigma_2^\ell,$ taking \refc{Aeq3} and \refc{Aeq2} into account  (see also \figc{descomposicio}), an easy computation shows that
\begin{align}\label{Aeq13}
  D_0(s)&=ds^\lambda\text{ with $d\!:=\varepsilon_1\varepsilon_2^{-\lambda}$}
\intertext{and}\label{Aeq8}
 T^0(s)&=\int_s^{\varepsilon_2}\left.\frac{u_1^{n_1}u_2^{n_2}}{P_1(0,0)}\right|_{u_2=\varepsilon_1\left(\frac{s}{u_1}\right)^{\lambda}}\frac{du_1}{u_1}=T_1^0s^{n_1}+T_2^0(ds^{\lambda})^{n_2},
\intertext{where}\notag
T_1^0&\!:=\frac{-\varepsilon_1^{n_2}}{(n_1-\lambda n_2)P_1(0,0)}
\text{ and } 
T_2^0\!:=\frac{\varepsilon_2^{n_1}}{(n_1-\lambda n_2)P_1(0,0)}.
\end{align} 
(Here, on account of $\lambda_0\notin\Q$, we reduce $V$ so that $n_1-\lambda n_2\neq 0$ for all $\np\in V.$) 
Hence $D(s)=R_2(dR_1^\lambda(s))$. If we take any strictly positive $\beta(\np)\in\cc^0(V)$ then, due to $R_{11}>0$, 
\begin{equation}\label{Aeq5}
R_1^\beta(s)=s^\beta R_{11}^\beta\left(1+\sum_{k=2}^{L_1}\frac{R_{1k}}{R_{11}}s^{k-1}\right)^\beta+\F_{L_1}^0(I\times V)=s^\beta R_{11}^\beta\sum_{\ell=0}^{L_1-1}\Upsilon^{[\beta]}_\ell s^\ell+\F_{L_1}^0(I\times V),
\end{equation}
where in the first equality we apply by $(a)$ in \lemc{composition} and in the second one we define $\Upsilon^{[\beta]}_\ell=\Upsilon^{[\beta]}_\ell(\varepsilon_1,\np)$ for $\ell=0,1,\ldots,L_1-1$ as the $\cc^0(I\times V)$ functions verifying 
\begin{equation}\label{Aeq18}
\left(1+\sum_{k=2}^{L_1}\frac{R_{1k}}{R_{11}}s^{k-1}\right)^\beta=\sum_{\ell=0}^{L_1-1}\Upsilon^{[\beta]}_\ell s^\ell
+\F_{L_1}^0(I\times V).
\end{equation}
(Here we apply Taylor's theorem at order $L_1$ to the function $x\mapsto (1+x)^\beta$ taking a uniform estimate of the remainder by means of its integral form.) Note in particular that $\Upsilon_0^{[\beta]}=1.$ Taking \refc{Aeq5} with $\beta(\np)=\lambda$ and applying $(b)$ in \lemc{composition} we obtain
\[
D(s)=\textstyle R_2(dR_1^\lambda(s))=\sum\limits_{k=1}^{L_2}R_{2k}d^kR_1^{\lambda k}(s)
+\F_{L_1}^0\!\left(\left\{(\varepsilon,\np)\in I^2\!\times\! V:\lambda>\frac{L_1}{L_2+1}\right\}\right)
\]
Now we choose $L_1$ and $L_2$ such that $\lambda_0>\frac{L_1}{L_2+1}$ and we shrink $V$ if necessary in order that
$\lambda>\frac{L_1}{L_2+1}$ for all $\np\in V.$ In doing so we get that
\[
D(s)=\textstyle R_2(dR_1^\lambda(s))=\sum\limits_{k=1}^{L_2}R_{2k}d^kR_1^{\lambda k}(s)
+\F_{L_1}^0(I^2\!\times\! V).
\]
Next, by taking \refc{Aeq5} with $\beta(\np)=\lambda k$, $k=1,2,\ldots,L_2$,
\begin{align*}
 D(s)&=\sum\limits_{k=1}^{L_2}\sum\limits_{\ell=0}^{L_1-1}R_{2k}R_{11}^{\lambda k}d^k\Upsilon^{[\lambda k]}_\ell s^{\ell+\lambda k}+\F_{L_1}^0(I^2\!\times\! V)\\
&=s^\lambda\sum\limits_{\ell=0}^{L_1-1}\sum\limits_{k=0}^{L_2-1}
R_{2,k+1}R_{11}^{\lambda(k+1)}d^{k+1}\Upsilon_\ell^{[\lambda(k+1)]}s^{\ell+\lambda k}
+\F_{L_1}^0(I^2\!\times\! V).
\end{align*}
Since $\lambda_0\notin\Q$, assertion $(a1)$ in \teoc{oldA} shows that
 \begin{equation}\label{Aeq6}
 \Delta_{\ell k}=R_{2,k+1}R_{11}^{\lambda(k+1)}d^{k+1}\Upsilon_\ell^{[\lambda(k+1)]}
 \text{ for all $(\varepsilon,\np)\in I^2\!\times\! V.$}
 \end{equation}
Here we also take \obsc{rm1} into account, shrinking (if necessary) the neighbourhood $V$ of $\np_0=(\lambda_0,\mu_0)$ in order that all the exponents $\ell+\lambda k$ are different for every $\np\in V.$ At this point it is worth to make the following remarks with regard to the previous equality:        
\begin{itemize}
\item It gives the expression of $\Delta_{ij}$ provided that $0\leqslant i\leqslant L_1-1$, $0\leqslant j\leqslant L_2-1$ 
         and $i+\lambda_0 j<L_1.$ Since we are just interested 
         in $(i,j)\in\{(0,0),(0,1),(1,0),(1,1),(i_1,j_1)\}$, these conditions 
         reduce to specific lower bounds for $L_1$ and $L_2$ that depend only on $\lambda_0,$ $i_1$ and $j_1.$
         For instance, in order to prove that the factorization in \refc{Aeq30} holds for $(i,j)=(i_1,j_1)$ we need that
         \[
          L_1>\max(i_1+\lambda_0j,i_1+1)\text{ and }L_2> j_2+1.
         \]
         This does not constitute a problem because we can take ${K}$, and therefore $L_1$ and $L_2$, arbitrarily large. 
\item The coefficient $\Delta_{\ell k}$ is a function that depends only on $\np$, whereas each function on the right
         hand side of \refc{Aeq6} depends on $\np$ but also on $\varepsilon$. This constitutes a key point that we 
         will exploit in the forthcoming arguments. Particularized to $\ell=0$, from \refc{Aeq13} and \refc{Aeq6} we get that
         \begin{equation}\label{Aeq15}
         \Delta_{0k}=\big(R_{2,k+1}\varepsilon_2^{-\lambda(k+1)}\big)\big(R_{11}^\lambda\varepsilon_1\big)^{k+1}
         \end{equation}
         does not depend on $\varepsilon=(\varepsilon_1,\varepsilon_2).$ Since the first factor does not depend on 
         $\varepsilon_1$ and the second one does not depend on $\varepsilon_2$, taking $k=0$ and using that 
          $\Delta_{00}(\np)\neq 0$ for all $\np\in\hat W,$ we conclude that
         \begin{equation*}
         R_{2,1}(\varepsilon_2,\np)\varepsilon_2^{-\lambda}
         \text{ and }
         R_{11}^\lambda(\varepsilon_1,\np)\varepsilon_1
         \text{ do not depend on $\varepsilon$,}
         \end{equation*}
         which in its turn, again from \refc{Aeq15}, implies that
         \begin{equation}\label{Aeq14}
         R_{2,k+1}(\varepsilon_2,\np)\varepsilon_2^{-\lambda(k+1)}
         \text{ does not depend on $\varepsilon$ for all $k\geqslant 1$.}
         \end{equation}
\end{itemize}

         Since $\Upsilon_0^{[\beta]}=1$ for any function $\beta,$ the factorization in \refc{Aeq6} also shows that 
         \begin{equation}\label{Aeq35}
          \Delta_{\ell k}=\Upsilon_\ell^{[\lambda(k+1)]}\Delta_{0k}.
         \end{equation}
         Consequently
         \begin{align*}
         \sum_{\ell=0}^{L_1-1}\Upsilon^{[\lambda(k+1)]}_\ell s^\ell&+\F_{L_1}^0(I\times V)
         =\left(1+\sum_{\ell=2}^{L_1}\frac{R_{1\ell}}{R_{11}}s^{\ell-1}\right)^{\lambda (k+1)}
         \\[5pt]
         &=\left(\;\sum_{\ell=0}^{L_1-1}\Upsilon^{[\lambda]}_\ell s^{\ell}+\F_{L_1}^0(I\times V)\right)^{k+1}
         =\left(\;\sum_{\ell=0}^{L_1-1}\frac{\Delta_{\ell 0}}{\Delta_{00}} s^{\ell}+\F_{L_1}^0(I\times V)\right)^{k+1}
         \\[5pt]
         &=\left(\;\sum_{\ell=0}^{L_1-1}\frac{\Delta_{\ell 0}}{\Delta_{00}} s^{\ell}\right)^{k+1}+\F_{L_1}^0(I\times V)
         =\sum_{\ell=0}^{L_1-1}\Omega_{\ell k}s^\ell+\F_{L_1}^0(I\times V),
         \end{align*}
         where in the first and second equalities we use the definition of $\Upsilon_\ell^{[\beta]}$ in \refc{Aeq18} with 
         $\beta(\np)=\lambda(k+1)$ and $\beta(\np)=\lambda$, respectively, in the third one we use 
         \refc{Aeq35} with $k=0$, in the fourth one we apply the binomial formula and \lemc{FLK} 
         and, finally, the last one follows from the definition in \refc{Aeq33}. Clearly this implies that
         \begin{equation}\label{Aeq36} 
          \Upsilon_\ell^{[\lambda(k+1)]}=\Omega_{\ell k}\text{ for $\ell=0,1,\ldots,L_1-1$.}
          \end{equation}
         Particularized to $(\ell,k)=(i_1,j_1)$, from \refc{Aeq35} once again we obtain that
         \[
         \Delta_{i_1j_1}=\Upsilon_{i_1}^{[\lambda(j_1+1)]}\Delta_{0j_1}=\Omega_{i_1j_1}\Delta_{0j_1}.
         \]
         This identity holds for all $\np\in V.$ On account of the considerations explained in the beginning of the proof this 
         shows that the assertion in \refc{Aeq30} is true for $(i,j)=(i_1,j_1)$ as desired.                  

We turn now to the study of the coefficients of the Dulac time. For convenience we write it as 
\[
T(s)=T^-(s)+T^+(s),
\]
where we define, recall \refc{Aeq7} and \refc{Aeq8},
\[
T^-(s)\!:=T^1(s)+T^0_1R_1^{n_1}(s)\text{ and }
T^+(s)\!:=\big(T^2(u)+T^0_2 u^{n_2}\big)\big|_{u=D_0(R_1(s))}.
\]
With respect to the first summand we observe that, from \refc{Aeq9} and taking \refc{Aeq5} with $\beta(\np)=n_1$, 
\begin{equation}\label{Aeq10}
 T^-(s)=\sum_{k=n_1}^{L_1-1}T_{k0}^-s^k+\F_{L_1}^0(I\times V)\text{ where }T_{k0}^-\!:=T_k^1+T_1^0R_{11}^{n_1}\Upsilon_{k-n_1}^{[n_1]}.
\end{equation}
On the other hand, from \refc{Aeq9}, we can write $T^2(u)+T^0_2 u^{n_2}=\sum\limits_{k=n_2}^{L_2}\bar T_k^2u^k+\F_{L_2+1}^0(I\times V)$ where
\begin{equation}\label{Aeq20}
\bar T^2_k\!:=\left\{\begin{array}{cc} T^2_k+T^0_2 & \text{if $k=n_2$,}\\[3pt]
T^2_k & \text{if $k>n_2$.}\end{array}\right.
\end{equation}
Consequently, taking \refc{Aeq5} with $\beta(\np)=\lambda$ and applying $(b)$ in \lemc{composition} we obtain
\begin{align}\notag
T^+(s)=\big(T^2(u)+T^0_2 u^{n_2}\big)\big|_{u=dR_1^\lambda(s)}
&=%\textstyle
\sum\limits_{k=n_2}^{L_2}\bar T^2_k d^kR_1^{\lambda k}
%+\F_{L_1}^K\left(\left\{(\varepsilon,\np)\in I^2\!\times\! V:\lambda>\frac{L_1}{L_2+1}\right\}\right)
+\F_{L_1}^0(I^2\!\times\! V)
\\\notag
&=%\textstyle
\sum\limits_{k=n_2}^{L_2}\bar T^2_k d^k \left(
s^{\lambda k}R_{11}^{\lambda k}\sum\limits_{\ell=0}^{L_1-1}\Upsilon_{\ell}^{[\lambda k]}s^\ell+\F_{L_1}^0(I\times V)
\right)
+\F_{L_1}^0(I^2\!\times\! V)\\\label{Aeq11}
&=%\textstyle
\sum\limits_{k=n_2}^{L_2}\sum\limits_{\ell=0}^{L_1-1}T^+_{\ell k}s^{\ell+\lambda k}+\F_{L_1}^0(I^2\!\times\! V).
\end{align}
Here we also use $\lambda>\frac{L_1}{L_2+1}$ for all $\np\in V$ in the first equality, in the second one we take \refc{Aeq5} with $\beta(\np)=\lambda k$, whereas in the last one we use that $d=\varepsilon_1\varepsilon_2^{-\lambda}$ and define
\begin{equation}\label{Aeq21}
T_{\ell k}^+\!:=(\bar T^2_k\varepsilon_2^{-\lambda k})(\varepsilon_1R_{11}^\lambda)^k\Upsilon_{\ell}^{[\lambda k]}.
\end{equation}
%\begin{equation}\label{Aeq21}
%T_{\ell k}^+\!:=T_{0k}^+\Upsilon_{\ell}^{[\lambda k]}\text{ with $T_{0k}^+\!:=\bar T^2_k(d\,R_{11}^\lambda)^k=(\bar T^2_k\varepsilon_2^{-\lambda k})(\varepsilon_1R_{11}^\lambda)^k$.}
%\end{equation}
Note that $T_{\ell 0}^+=0$ for all $\ell\geqslant 1$ due to $\Upsilon_{\ell}^{[0]}=0$ for all $\ell\geqslant 1.$ Consequently, since $T_{\ell k}$ is by definition the coefficient of $s^{\ell+\lambda k}$ in $T(s)=T^-(s)+T^+(s)$, from \refc{Aeq10} and \refc{Aeq11} we get that
\begin{equation}\label{Aeq22}
 T_{\ell k}=\left\{
 \begin{array}{ll}
  T_{\ell k}^+ & \text{ if $k>0.$} \\
  T_{\ell 0}^- & \text{ if $k=0$ and $\ell\geqslant 1.$}
 \end{array}
 \right.
\end{equation}
(To be more precise, the above equality follows from \obsc{rm1} and by applying $(b1)$ in \teoc{oldA} thanks to $\lambda_0\notin\Q$ and shrinking, if necessary, the neighbourhood $V$ of $\np_0=(\lambda_0,\mu_0)$ in order that all the exponents $\ell+\lambda k$ are different for every $\np\in V.)$ Finally, since the coefficient $T_{00}$ only exists in case the that $n_1n_2=0$ and $n\neq (0,0)$ by hypothesis, we have that
\[
 T_{00}=
   \left\{\begin{array}{ll} 
   T_{00}^-& \text{if $n_1=0$,}\\[3pt]
   T_{00}^+& \text{if $n_2=0$.}\end{array}\right.
\]
Similarly as we noted previously for $\Delta_{ij}$, let us remark that since we are only interested in the coefficients~
\[
T_{ij}\text{ with }(i,j)\in\{(n_1,0),(n_1+1,0),(0,n_2),(0,n_2+1),(i_2,j_2)\},
\]
from \refc{Aeq10} and \refc{Aeq11} we get specific lower bounds for $L_1$ and $L_2$ to be satisfied. Once again, this is not a problem because these lower bounds are given in terms of $\lambda_0$, $n_1$, $n_2$, $i_2$ and $j_2$ and, on the other hand, we can take $K,$ and so $L_1$ and $L_2,$ arbitrarily large. For instance, in order to show that the factorization in~\refc{Aeq31} holds for $(i,j)=(i_2,j_2)$ with $j_2>0$ we argue as follows. Precisely due to $j_2>0$, we get that
\[
 T_{i_2j_2}=T_{i_2j_2}^+=T_{0j_2}^+\Upsilon_{i_2}^{[\lambda j_2]}=T_{0j_2}\Omega_{i_2,j_2-1},
\]
where in the first equality we take \refc{Aeq22} into account, the second one follows readily from \refc{Aeq21} thanks to
$\Upsilon_0^{[\lambda j_2]}=1$, and in the last one we apply the identity in \refc{Aeq36}. For this to happen, see also \refc{Aeq11}, we need that
\[
 L_1>\max(i_2+1,i_2+\lambda_0 j_2)\text{ and }L_2>j_2.
\]
This shows the validity of the factorization for all $\np\in V.$ As we explained at the beginning of the proof, this factorization extends to all $\np=(\lambda,\mu)\in\hat W$ with $\lambda\notin D_{i_2j_2}^n\cup D_{i_20}^0$ by continuity and the fact that $D_{i_2j_2}^n\cup D_{i_20}^0$ is a discrete subset of rational numbers in $(0,+\infty).$ 

So far we have proved \refc{Aeq30} and \refc{Aeq31}, which constitute assertion \refc{a} in the statement. In doing so we have also identified all the elements needed to compute $\Delta_{ij}$ and $T_{ij}$ but recall that we must only analyze the cases $(i,j)\in\{(0,0),(1,0),(0,1),(1,1)\}$ and $(i,j)\in\{(n_1,0),(n_1+1,0),(0,n_2),(0,n_2+1)\},$ respectively. 
With this aim in view we shall apply \lemc{L3} to obtain the explicit expressions of the coefficients $R_{i1},$ $R_{i2},$ $T^i_{n_i}$ and $T^i_{n_i+1}$ in \refc{Aeq9} for $i=1,2$. Let us advance that the formulae for $i=1$ and $i=2$ are related by switching $\lambda$ and $1/\lambda$,  $\sigma$ and $\tau$, the subscripts $1$ and $2$ (with the exception of the third subscript~$k$ in $\sigma_{ijk}$ and $\tau_{ijk}$) and by exchanging the order of the variables in the functions $f_i$ and $h_i$.  

For the reader's convenience we sum up in \tabc{atlas} the fundamental information for applying the results in 
Appendix~\ref{ap_regular} to study the regular passages, see \figc{descomposicio}, together with the functions $L_i$ defined in~\refc{def_fun} and the functions $f_i$ and $h_i$ given in \refc{Aeq0}. 
\begin{table}[t]
  \begin{center}
    \begin{tabular}{c|c|c} 
       & \textsc{First regular} & \textsc{Second regular} \\[-1pt]
       & \textsc{passage} & \textsc{passage} \\[3pt]
      \hline
       & &  \\[-7pt]
       $\ell$ & $n_1$ & $n_2$ \\[5pt]
       $\nu$ & $(\varepsilon_1,\np)$ & $(\varepsilon_2,\np)$ \\[7pt]
       $h(x,y)$ & $\frac{P_1(y,x)}{xP_2(y,x)}$ & $\frac{P_2(x,y)}{xP_1(x,y)}$  \\[7pt]
       $H(x,y)$ & $\left(\frac{y}{x}\right)^{\frac{1}{\lambda}}\!\frac{L_1(x)}{L_1(y)}$ 
                          & $\left(\frac{y}{x}\right)^{\lambda}\!\frac{L_2(x)}{L_2(y)}$  \\[7pt]
       $f(x,y)$ & $\frac{x^{n_2-1}}{P_2(y,x)}$ & $\frac{x^{n_1-1}}{P_1(x,y)}$ \\[7pt]
       $\xi(s;\nu)$ & $\big(\sigma_{12}(s;\np),\sigma_{11}(s;\np)\big)$ 
                         & $\big(\tau_{21}(s;\varepsilon_2,\np),\tau_{22}(s;\varepsilon_2,\np)\big)$ \\[7pt]
       $\zeta(s;\nu)$ & $\big(\tau_{12}(s;\varepsilon_1,\np),\tau_{11}(s;\varepsilon_1,\np)\big)$ 
                        & $\big(\sigma_{21}(s;\np),\sigma_{22}(s;\np)\big)$ \\
    \end{tabular}
  \end{center}
  \caption{Information related with the application of the results in Appendix~\ref{ap_regular}.  
               The auxiliary sections $\Sigma_1^\ell$ and $\Sigma_2^\ell$ are given by 
               $\tau_1(s;\varepsilon_1,\np)=\Phi(s,\varepsilon_1;\np)$ 
               and $\tau_2(s;\varepsilon_2,\np)=\Phi(\varepsilon_2,s;\np)$, respectively, see \refc{Aeq2}.}\label{atlas}
\end{table}
On account of this the application of \lemc{L2} yields
\begin{align}\label{Aeq25}
\rho_{11}(x)&=\alpha_{11}x^{\frac{-1}{\lambda}}L_1(x)\text{ with } \alpha_{11}\!:=\frac{\sigma_{111}\sigma_{120}^{\frac{1}{\lambda}}}{L_1(\sigma_{120})}\\
\intertext{for the first regular passage and}\notag
\rho_{21}(x)&=\alpha_{21}x^{-{\lambda}}L_2(x)\text{ with }\alpha_{21}\!:=\frac{\tau_{221}\tau_{210}^{{\lambda}}}{L_2(\tau_{210})}
\end{align}
for the second one. (Here, to be consistent with the previous notation, the subscript $i$ in $\rho_{ij}$ refers to the 
first or second regular passage, whereas $j$ refers to the derivation's order.) Next, by applying \lemc{L3},
\begin{equation}\label{Aeq27}
R_{11}%=\frac{\rho_{11}(\tau_{120})}{\tau_{111}}
=\alpha_{11}\frac{\tau_{120}^{\frac{-1}{\lambda}}L_1(\tau_{120})}{\tau_{111}}
\text{ and }
R_{21}%=\frac{\rho_{21}(\sigma_{210})}{\sigma_{221}}
=\alpha_{21}\frac{L_2(\sigma_{210})}{\sigma_{221}\sigma_{210}^\lambda}.
\end{equation}
Observe at this point that $\alpha_{11}$ does not depend on $\varepsilon$ and that, see \refc{Aeq14}, this is also the case of $R_{11}^\lambda\varepsilon_1.$ From the first equality in \refc{Aeq27}, this implies that $\frac{L_1^\lambda(\tau_{120})}{\tau_{111}^\lambda\tau_{120}}\varepsilon_1$ does not depend on $\varepsilon$. On the other hand, $\tau_{120}=\varepsilon_1\psi_2(0,\varepsilon_1)$ and $\tau_{111}=\psi_1(0,\varepsilon_1),$ see~\refc{Aeq2}, together with $\psi_i(0,0)=L_1(0)=1,$ imply that $\lim_{\varepsilon_1\to 0}\frac{L_1^\lambda(\tau_{120})}{\tau_{111}^\lambda\tau_{120}}\varepsilon_1=1$. Thus $\frac{L_1^\lambda(\tau_{120})}{\tau_{111}^\lambda\tau_{120}}\varepsilon_1=1$ and, consequently, $R_{11}^\lambda\varepsilon_1=\alpha_{11}^\lambda$. In short,
\begin{equation}\label{Aeq16}
\frac{\tau_{111}^\lambda\tau_{120}}{L_1^\lambda(\tau_{120})}=\varepsilon_1
\text{ and }
R_{11}=\alpha_{11}\varepsilon_1^{-1/\lambda}=\frac{\sigma_{111}\sigma_{120}^{1/\lambda}}{L_1(\sigma_{120})}\varepsilon_1^{-1/\lambda}.
\end{equation}
Furthermore, from \refc{Aeq14} again, $R_{21}\varepsilon_2^{-\lambda}$ does not depend on $\varepsilon.$ This implies, on account of the second equality in \refc{Aeq27}, that $\alpha_{21}\varepsilon_2^{-\lambda}$ does not depend on $\varepsilon$ neither. Then, taking $\varepsilon_2\to 0$ exactly as before, we conclude that 
\begin{equation}\label{Aeq17}
 \alpha_{21}=\varepsilon_2^\lambda.
\end{equation}
Therefore $R_{21}=\varepsilon_2^\lambda\frac{L_2(\sigma_{210})}{\sigma_{221}\sigma_{210}^\lambda}$ and consequently, from \refc{Aeq15},
\[
\Delta_{00}(\np)=(R_{11}^\lambda\varepsilon_1)(R_{21}\varepsilon_2^{-\lambda})=\frac{\sigma_{111}^\lambda\sigma_{120}}{L_1^\lambda(\sigma_{120})}\frac{L_2(\sigma_{210})}{\sigma_{221}\sigma_{210}^\lambda}
\text{ for all $\np\in V.$}
\]
On account of the considerations explained in the first paragraph of the proof, this shows the validity of the first equality in \refc{b} for all $\np=(\lambda,\mu)\in (0,+\infty)\times W$. Indeed, following the notation introduced there,~$\tilde\Delta_{00}$ is the function on the right hand side of the above equality, which belongs to $\cc^{\infty}\big((0+\infty)\times W\big)$ by \lemc{fun_ok}, i.e., $\tilde D_{00}^0=\emptyset,$ and we have on the other hand, see \obsc{domains}, $D_{00}^0=\emptyset$ as well.

Next we proceed with the computation of the second order derivatives in \lemc{L2}. Using the first column in \tabc{atlas}, some long but easy computations show that
\begin{align}\notag
\rho_{12}(x)&=\frac{\alpha_{11}}{\sigma_{111}}x^{\frac{-1}{\lambda}}L_1(x)
\bigg(\sigma_{112}-\frac{2\sigma_{121}\sigma_{111}}{\sigma_{120}}\left(\frac{P_1}{P_2}\right)\!(0,\sigma_{120})+2\sigma_{111}\alpha_{11}\int_{\sigma_{120}}^x\underbrace{L_1(u)\,\partial_1\!\left(\frac{P_1}{P_2}\right)\!(0,u)}_{M_1(u)}\, u^{\frac{-1}{\lambda}} \frac{du}{u}\bigg)\\\label{Aeq29}
&
=\alpha_{12}x^{\frac{-1}{\lambda}}L_1(x)+2\alpha_{11}^2x^{\frac{-2}{\lambda}}L_1(x)\hat M_1(1/\lambda,x),
\intertext{ for all $x\in I_1\cap (0,+\infty)$ with}
\alpha_{12}&\!:=\frac{\alpha_{11}}{\sigma_{111}}\left(\sigma_{112}-\frac{2\sigma_{121}\sigma_{111}}{\sigma_{120}}\left(\frac{P_1}{P_2}\right)\!(0,\sigma_{120})\right)-2\alpha_{11}^2\sigma_{120}^{\frac{-1}{\lambda}}\hat M_1(1/\lambda,\sigma_{120}).\label{Aeq19}
\end{align}
Here we use for the first time the properties of the incomplete Mellin transform introduced in Appendix~\ref{Mellin}. 
More concretely, by \lemc{fun_ok}, $M_1(u;\np)\in\cc^\infty(I_1\times\hat W)$ with $0\in I_1$. Hence, by applying \teoc{L8} there exists a unique $\hat M_1(\alpha,u;\np)\in\cc^\infty((\R\setminus\Z_{\geq 0})\times I_1\times\hat W)$ such that 
$\partial_u\big(\hat M_1(\alpha,u)u^{-\alpha}\big)=M_1(u)u^{-\alpha-1}$ for all $u\in I_1\cap(0,+\infty)$. Analogously, taking the second column in \tabc{atlas}, one can also verify that
\begin{align}\label{Aeq28}
\rho_{22}(x)&=\alpha_{22}x^{-\lambda}L_2(x)+2\alpha_{21}^2x^{-2\lambda}L_2(x)\hat M_2(\lambda,x)
\text{ for all $x\in I_2\cap (0,+\infty)$,}
\intertext{with}\notag
 \alpha_{22}&\!:=\frac{\alpha_{21}}{\tau_{221}}\left(\tau_{222}-\frac{2\tau_{211}\tau_{221}}{\tau_{210}}\left(\frac{P_2}{P_1}\right)\!(\tau_{210},0)\right)-2\alpha_{21}^2\tau_{210}^{-\lambda}\gorro{M}_2(\lambda,\tau_{210}).
\end{align}
We claim that $\alpha_{22}=\varepsilon_2^\lambda\varphi_1(\varepsilon_2,\np)$ with $\varphi_1\in\cc^{K}\big((-\delta,\delta)\times V\big).$ Indeed, this is so due to the following facts:
\begin{enumerate}
\item $P_1(x_1,x_2;\np)$ and $P_2(x_1,x_2;\np)$ are $\cc^\infty$ and do not vanish on 
        $x_2=0$ and $x_1=0$, respectively.
\item $L_2(u;\np)$ and $M_2(u;\np)$ are $\cc^\infty(I_2\times\hat W)$ by \lemc{fun_ok} 
         and the first one does not vanish.
\item The parametrization $\tau_2(s;\varepsilon_2,\np)$ of the section 
         $\Sigma_2^\ell$ is defined by means of~$\Phi\in\cc^{K}(U\times V)$, see \refc{Aeq2}, where recall that 
         $U=(-\delta,\delta)\times (-\delta,\delta),$
\item and therefore, the map $(\varepsilon_2,\np)\mapsto\gorro{M}_2(\lambda,\tau_{210};\np)$ belongs to 
         $\cc^{K}\big((-\delta,\delta)\times V\big)$ by $(a)$ 
        in \teoc{L8} since $\lambda\notin\Z_{\geq 0}$ due to $\lambda_0\notin\Q$ and shrinking $V$ if necessary.           
\item $\tau_{221}=\psi_2(\varepsilon_2,0)$ and $\tau_{210}=\varepsilon_2\psi_1(\varepsilon_2,0)$ with $\psi_i(0,0)=1.$
         Moreover, see \refc{Aeq17}, $\alpha_{21}=\varepsilon_2^\lambda.$        
\end{enumerate}
The key point for our purposes will be that, for each fixed $\np$, the function $\varphi_1$ is $\cc^{K}$ in a neighbourhood of $\varepsilon_2=0.$ On account of this, for simplicity in the exposition we will say that $\alpha_{22}=\varepsilon_2^\lambda\varphi_1(\varepsilon_2)$ with $\varphi_1\in\cc^{K}$. In what follows we will deal several times with this type of situation and for shortness we will omit the previous details. More generally, for the same reason, when we write $\varphi_k(\varepsilon_i)$ with $i=1,2$ and any subscript~$k$ we shall mean that $\varphi_k$ is some function depending only on $\varepsilon_i$ and $\np$ that belongs to $\cc^{K}((-\delta,\delta)\times V\big)$.

We are now in position to compute the second order derivatives by means of \lemc{L3}. In this case, for the sake of convenience in the exposition, we begin with the second regular passage. In doing so, and using \tabc{atlas} together with the expressions for $R_{21}$ and $\rho_{22}$ given in \refc{Aeq27} and \refc{Aeq28}, respectively, we get
\begin{align*}
R_{22}=\left(\frac{\sigma_{211}}{\sigma_{210}}\!\left(\frac{P_2}{P_1}\right)\!(\sigma_{210},0)-\frac{\sigma_{222}}{2\sigma_{221}}\right)\alpha_{21}^2\frac{\sigma_{210}^{-2{\lambda}}}{\sigma_{221}^2}L_2^2(\sigma_{210})
&+\frac{\alpha_{22}}{2}\frac{\sigma_{210}^{-{\lambda}}}{\sigma_{221}}L_2(\sigma_{210})\\
&+\alpha_{21}^2\frac{\sigma_{210}^{-2{\lambda}}}{\sigma_{221}}L_2(\sigma_{210})\gorro{M}_2(\lambda,\sigma_{210}).
\end{align*}
This implies that $\alpha_{22}\varepsilon_2^{-2\lambda}$ does not depend on $\varepsilon$ because this is the case for $\sigma_2$ and $R_{22}\varepsilon_2^{-2\lambda}$, see~\refc{Aeq14}, and moreover $\alpha_{21}=\varepsilon_2^\lambda$  from \refc{Aeq17}. Hence the previous claim shows that $\alpha_{22}\varepsilon_2^{-2\lambda}=\varepsilon_2^{-\lambda}\varphi_1(\varepsilon_2)=c$ where $c$ is a constant depending only on $\np.$ Therefore $\varphi_1(\varepsilon_2)=c\varepsilon_2^{\lambda}$. Since $\lambda_0\notin\Q$, we have that $\lambda\notin\Z_{\geq 0}$ for all $\np\in V$ (shrinking $V$ if necessary) and, consequently, $c=0$ because $\varphi_1$ is $\cc^{K}$ in a neighbourhood of $\varepsilon_2=0$ with ${K}$ arbitrarily large. (More precisely it suffices to take ${K}>\lambda_0$ and make smaller $V$ so that $K>\lambda$ for all $\np\in V.$) Accordingly 
\begin{equation}\label{Aeq24}
\alpha_{22}=0
\end{equation}
and, since $\alpha_{21}=\varepsilon_2^\lambda$ on account of \refc{Aeq17},
\begin{equation}\label{Aeq26}
R_{22}=-\varepsilon_2^{2\lambda}
\bigg(
\underbrace{
\frac{\sigma_{222}}{2\sigma_{221}}
-\frac{\sigma_{211}}{\sigma_{210}}\!\left(\frac{P_2}{P_1}\right)\!(\sigma_{210},0)
-\frac{\sigma_{221}}{L_2(\sigma_{210})}\gorro{M}_2(\lambda,\sigma_{210})}_{S_2}
\bigg)
\left(\frac{L_2(\sigma_{210}}{\sigma_{221}\sigma_{210}^{\lambda}}\right)^2.
\end{equation}
Then, using \refc{Aeq15} with $k=1$ and the expression of $R_{11}$ in \refc{Aeq16}, 
\[
 \Delta_{01}=-S_2
 \left(\frac{L_2(\sigma_{210}}{\sigma_{221}\sigma_{210}^{\lambda}}\right)^2
 \left(\frac{\sigma_{111}^\lambda\sigma_{120}}{L_1^\lambda(\sigma_{120})}\right)^2
 =-S_2\Delta_{00}^2\text{ for all $\np\in V.$}
\]
By applying \lemc{fun_ok} and \teoc{L8}, the function $\gorro{M}_2(\lambda,\sigma_{210})$ in $S_2$ is $\cc^\infty$ in a neighbourhood of any $(\lambda_\star,\mu_\star)\in (0,+\infty)\times W$ such that $\lambda_\star\notin\Z_{\geq 0}.$ Thus the function on the right hand side of the above equality, that we denote by~$\tilde\Delta_{01}$ in the second paragraph of the proof, is $\cc^\infty$ on $((0,+\infty)\setminus\tilde D_{01}^n)\times W$ with $\tilde D_{01}^n\!:=\N$. Since we know on the other hand by \teoc{oldA} that $\Delta_{01}\in\cc^\infty (((0,+\infty)\setminus  D_{01}^n)\times W)$ with $D_{01}^n=\N$, see \obsc{domains}, this implies by continuity that the second equality in \refc{b} is true for $(\lambda,\mu)\in \big((0,+\infty)\setminus D_{01}^0\big)\times W.$ Certainly we also use here, and it is essential, that the parameter $\np_0=(\lambda_0,\mu_0)\in\hat W$ with $\lambda_0\notin\Q$ that we fix at the very beginning is arbitrary.

Let us begin now with the computation of $R_{21}$, i.e., the second coefficient of the transition map for the first passage, by means of \lemc{L3}. In this case, using \tabc{atlas} together with \refc{Aeq16} and \refc{Aeq29}, 
we get
\[
R_{12}=\left(\frac{\tau_{121}}{\tau_{120}}\!\left(\frac{P_1}{P_2}\right)\!(0,\tau_{120})-\frac{\tau_{112}}{2\tau_{111}}\right)\!\alpha_{11}^2\underbrace{\frac{\tau_{120}^{{\frac{-2}{\lambda}}}}{\tau_{111}^2}L_1^2(\tau_{120})}_{\varepsilon_1^{{-2/\lambda}}}+\frac{\alpha_{12}}{2}\underbrace{\frac{\tau_{120}^{{\frac{-1}{\lambda}}}}{\tau_{111}}L_1(\tau_{120})}_{\varepsilon_1^{{-1/\lambda}}}+\alpha_{11}^2\underbrace{\frac{\tau_{120}^{{\frac{-2}{\lambda}}}}{\tau_{111}}L_1(\tau_{120})\gorro{M}_1(1/\lambda,\tau_{120})}_{\varepsilon_1^{{-2/\lambda}}\varphi_2(\varepsilon_1)}.
\]
Since $R_{11}=\varepsilon_1^{-1/\lambda}\alpha_{11}$ from \refc{Aeq16} once again and, on the other hand, $\tau_{120}=\varepsilon_1\psi_2(0,\varepsilon_1)$ with $\psi_2(0,0)=1$, it follows that we can write 
\[
 \frac{R_{12}}{R_{11}}=\varphi_3(\varepsilon_1)\varepsilon_1^{-1/\lambda-1}+\frac{\alpha_{12}}{2\alpha_{11}}.
\]
Observe that the quotient $\frac{R_{12}}{R_{11}}$ does not depend on $\varepsilon$ because, from \refc{Aeq18} and \refc{Aeq6}
\[
 \frac{\Delta_{1k}}{\Delta_{0k}}=\Upsilon_1^{[\lambda(k+1)]}=\lambda (k+1)\frac{R_{12}}{R_{11}}.
\]
Since this is also the case for the quotient $\frac{\alpha_{12}}{\alpha_{11}},$ see \refc{Aeq25} and \refc{Aeq19}, it turns out that $\varphi_3(\varepsilon_1)\varepsilon_1^{-1/\lambda-1}=c$ for some constant depending only on $\np$. 
Thus $\varphi_3(\varepsilon_1)=c\varepsilon_1^{1/\lambda+1}$ and, due to $\lambda\approx\lambda_0\notin\Q,$ this implies $c=0$. Therefore, 
\begin{equation}\label{Aeq23}
 \frac{R_{12}}{R_{11}}=\frac{\alpha_{12}}{2\alpha_{11}}=\frac{\sigma_{112}}{2\sigma_{111}}-\frac{\sigma_{121}}{\sigma_{120}}\!\left(\frac{P_1}{P_2}\right)\!(0,\sigma_{120})-\frac{\sigma_{111}}{L_1(\sigma_{120})}\gorro{M}_1(1/\lambda,\sigma_{120})=S_1,
 \end{equation}
where the second equality follows from~\refc{Aeq25} and~\refc{Aeq19} again and the last one from the definition in \refc{def_S}. Hence
\[
 \Delta_{10}=\Delta_{00}\lambda S_1\text{ and }\Delta_{11}=\Delta_{01}2\lambda S_1=-\Delta_{00}^22\lambda S_1S_2\text{ for all $\np\in V.$}
\]
On account of the expression of $S_2$ and $S_1$ given in \refc{Aeq26} and \refc{Aeq23}, respectively, the application of \teoc{L8} shows (following the notation introduced in the first paragraph of the proof) that $\tilde D_{10}^0=\frac{1}{\N}$ and $\tilde D_{11}^0=\N\cup\frac{1}{\N}.$ 
Since these sets coincide with $D_{10}^0$ and $D_{11}^0,$ respectively, this concludes the proof of assertion \refc{b}.

Let us show next the validity of the identities in assertion \refc{c}, that deal with the coefficients of the Dulac time. As before we begin with the study of the regular passages and the computation of the first coefficients of their time functions. With regard to $T^1(s;\varepsilon_1,\np)$ it turns out that
\[
 T^1_{n_1}=\alpha_{11}^{n_1}\int_{\sigma_{120}}^{\tau_{120}}\underbrace{\frac{L_1^{n_1}(x)}{P_2(0,x)}}_{A_1(x)}x^{n_2-\frac{n_1}{\lambda}}\frac{dx}{x}
 =\alpha_{11}^{n_1}\bigg(\underbrace{\tau_{120}^{n_2-\frac{n_1}{\lambda}}\gorro{A}_1(n_1/\lambda-n_2,\tau_{120})}_{\varepsilon_1^{n_2-\frac{n_1}{\lambda}}\varphi_4(\varepsilon_1)}-\sigma_{120}^{n_2-\frac{n_1}{\lambda}}\gorro{A}_1(n_1/\lambda-n_2,\sigma_{120})\bigg).
\]
The first equality above follows by \lemc{L3} taking into account the expression of $\rho_{11}$ in~\refc{Aeq25} and \tabc{atlas}. The second equality follows by applying \teoc{L8} with $A_1(x;\np)$, that belongs to $\cc^\infty(I_1\times\hat W)$ by \lemc{fun_ok}, and the fact that $\tau_{120}=\varepsilon_1\psi_2(0,\varepsilon_1)$ with $\psi_2(0,0)=1$. Then  
 \begin{align*}
 T_{n_10}&=T_{n_10}^-=T^1_{n_1}+T_1^0R_{11}^{n_1}\\[7pt]
 &=
 \alpha_{11}^{n_1}\left(
 \varepsilon_1^{n_2-\frac{n_1}{\lambda}}
 \left(\varphi_4(\varepsilon_1)-\frac{1}{(n_1-\lambda n_2)P_1(0,0)}\right)
 -\sigma_{120}^{n_2-\frac{n_1}{\lambda}}\gorro{A}_1(n_1/\lambda-n_2,\sigma_{120})
 \right)\\[7pt]
% =-\left(\frac{\sigma_{111}\sigma_{120}^{\frac{1}{\lambda}}}{L_1(\sigma_{120})}\right)^{n_1}\sigma_{120}^{n_2-\frac{n_1}{\lambda}}\gorro{A}_1(n_1/\lambda-n_2,\sigma_{120})
 &=-\frac{\sigma_{111}^{n_1}\sigma_{120}^{n_2}}{L_1^{n_1}(\sigma_{120})}\gorro{A}_1(n_1/\lambda-n_2,\sigma_{120}).
 \end{align*}
The first and second equalities above follow from \refc{Aeq22} and \refc{Aeq10}, respectively, and the third one by using~\refc{Aeq8} together with~\refc{Aeq16}. In the last equality we use that $T_{n_10}$, $\alpha_{11}=\frac{\sigma_{111}\sigma_{120}^{1/\lambda}}{L_1(\sigma_{120})}$ and $\sigma_{1}$ do not depend on $\varepsilon$ and this, on account of $\lambda\approx\lambda_0\notin\Q$, implies that $\varphi_4(\varepsilon_1)=\frac{1}{(n_1-\lambda n_2)P_1(0,0)}$. For the reader's convenience let us be more precise in this last implication because we use the same argument repeatedly. The point is that there exists~$c$, not depending on $\varepsilon_1$, such that 
\[
 \varepsilon_1^{n_2-\frac{n_1}{\lambda}}\left(\varphi_4(\varepsilon_1)-\frac{1}{(n_1-\lambda n_2)P_1(0,0)}\right)=c
 \text{ for all $\varepsilon_1$}
\] 
and we know on the other hand that $\varphi_4$ is $\cc^{K}((-\delta,\delta))$ with~${K}$ arbitrarily large. In this case for our purpose we need ${K}>\frac{n_1}{\lambda_0}-n_2$, so that (by shrinking~$V$) we have ${K}>\frac{n_1}{\lambda}-n_2$ for all $\np\in V$. Since $\lambda_0\notin\Q$ we can also assume that $\frac{n_1}{\lambda}-n_2\notin\Z_{\geq 0}$ for all $\np\in V.$ That being said, note then that from the above equality it turns out that $\varphi_4$ is a $\cc^{K}$ function that is written as $\varphi_4(\varepsilon_1)=c\varepsilon_1^{\frac{n_1}{\lambda}-n_2}+\hat c$ with the exponent $\frac{n_1}{\lambda}-n_2$ smaller than ${K}$ and not being in $\Z_{\geq 0}.$ It is evident that this is only possible if $c=0$, as we claimed. Hence
\[
T_{n_10}=\frac{\sigma_{111}^{n_1}\sigma_{120}^{n_2}}{L_1^{n_1}(\sigma_{120})}\gorro{A}_1(n_1/\lambda-n_2,\sigma_{120})\text{ for all $\np\in V.$}
\]
By \teoc{L8}, the function on the right hand side is $\cc^\infty$ in a neighbourhood of any $(\lambda_\star,\mu_\star)\in\hat W$ with $\frac{n_1}{\lambda_\star}-n_2\notin\Z_{\geq 0},$ i.e., $\lambda_\star\notin\tilde D_{n_1,0}^n\!:=\frac{n_1}{\N_{\geq n_2}}.$ Thus $\tilde D_{n_1,0}^n\subset D_{n_1,0}^n=\bigcup_{i=1}^{n_1}\frac{i}{\N_{\geq n_2}},$ see \obsc{domains}, and therefore by continuity the above equality is valid provided that $\lambda\notin D_{n_1,0}^n.$ This proves the first identity in~\refc{c}.  

Regarding the time function $T^2(s;\varepsilon_2,\np)$ of the second regular passage one can check that
\begin{align*}
T^2_{n_2}&=\varepsilon_2^{n_2\lambda}\bigg(\sigma_{210}^{n_1- n_2\lambda }\gorro{A}_2(n_2\lambda-n_1,\sigma_{210})-\underbrace{\tau_{210}^{n_1-\lambda n_2}\gorro{A}_2(n_2\lambda-n_1,\tau_{210})}_{\varepsilon_2^{n_1-n_2\lambda}\varphi_5(\varepsilon_2)}\bigg)\\
&=\varepsilon_2^{n_2\lambda}\sigma_{210}^{n_1-n_2\lambda}\gorro{A}_2(n_2\lambda-n_1,\sigma_{210})
+\varepsilon_2^{n_1}\varphi_5(\varepsilon_2),
\end{align*}
where the first equality follows by \lemc{L3} and on account of $\rho_{21}(x)=\varepsilon_2^\lambda x^{-\lambda}L_2(x)$, and the second equality by applying \teoc{L8} with $A_2(x;\np)$, that belongs to $\cc^\infty(I_2\times\hat W)$ by \lemc{fun_ok}. Hence, taking~\refc{Aeq8} and \refc{Aeq20} into account, 
\[
 \bar T^2_{n_2}=T^2_{n_2}+T^0_2=\varepsilon_2^{n_2\lambda}\sigma_{210}^{n_1-n_2\lambda}\gorro{A}_2(n_2\lambda-n_1,\sigma_{210})+\varepsilon_2^{n_1}\left(
 \varphi_5(\varepsilon_2)+\frac{1}{(n_1-\lambda n_2)P_1(0,0)}
 \right)
\]
and, accordingly, 
 \[
 T_{0,n_2}=T^+_{0,n_2}=(\bar T^2_{n_2}\varepsilon_2^{-n_2\lambda})(R_{11}^\lambda\varepsilon_1)^{n_2}=\sigma_{210}^{n_1-n_2\lambda}\gorro{A}_2(n_2\lambda-n_1,\sigma_{210})\left(\frac{\sigma_{111}^\lambda\sigma_{120}}{L_1^\lambda(\sigma_{120})}\right)^{n_2}, 
 \]
where the first and second equalities follow from \refc{Aeq22} and \refc{Aeq21}, respectively. Finally, in the last equality we use that $\sigma_1$ and $\sigma_2$ do not depend on $\varepsilon$ and that this is also the case for $T_{0,n_2}$ and, see \refc{Aeq16}, $R_{11}\varepsilon_1^{1/\lambda}=\frac{\sigma_{111}\sigma_{120}^{1/\lambda}}{L_1(\sigma_{120})}$. Since $\lambda\approx\lambda_0\notin\Q$, this implies $\varphi_5(\varepsilon_2)=\frac{-1}{(n_1-\lambda n_2)P_1(0,0)}$ and finishes the proof of the second identity in \refc{c}. 

We proceed next with the computation of the coefficient $T_{n_1+1}^1$. To this end we apply \lemc{L3} taking account of \tabc{atlas} and the expressions of $\rho_{11}$, $R_{11}$ and $\rho_{12}$ given in \refc{Aeq25}, \refc{Aeq16} and \refc{Aeq29}, respectively. In doing so we obtain
\begin{align*}
T^1_{n_1+1}&=\underbrace{\left(\varepsilon_1^{\frac{-1}{\lambda}}\alpha_{11}\right)^{n_1+1}\frac{\tau_{121}\tau_{111}^{n_1}\tau_{120}^{n_2-1}}{P_2(0,\tau_{120})}}_{\varepsilon_1^{n_2-\frac{n_1+1}{\lambda}}\varphi_6(\varepsilon_1)}-\frac{\sigma_{121}\sigma_{111}^{n_1}\sigma_{120}^{n_2-1}}{P_2(0,\sigma_{120})}+\alpha_{11}^{n_1+1}\int_{\sigma_{120}}^{\tau_{120}}L_1^{n_1+1}(x)x^{n_2-\frac{n_1+1}{\lambda}}\partial_1P_2^{-1}(0,x)\frac{dx}{x}\\[7pt] 
&\quad+\frac{n_1}{2}\alpha_{11}^{n_1-1}\int_{\sigma_{120}}^{\tau_{120}}L_1^{n_1-1}(x)x^{\frac{-(n_1-1)}{\lambda}}\left(\alpha_{12}x^{\frac{-1}{\lambda}}L_1(x)+2\alpha_{11}^2x^{\frac{-2}{\lambda}}L_1(x)\gorro{M}_1(1/\lambda,x)\right)\frac{x^{n_2-1}}{P_2(0,x)}dx.
\end{align*} 
Here we also use that $\tau_1$ does not depend on $\varepsilon_2$ and that $\tau_{120}$ and $\tau_{121}$ vanish at $\varepsilon_1=0.$ Then some easy manipulations first, on account of the definitions of $A_1$ and $B_1$ given in \refc{def_fun}, and next the application of \teoc{L8} yields to
\begin{align*}
T_{n_1+1}^1&=\varepsilon_1^{n_2-\frac{n_1+1}{\lambda}}\varphi_6(\varepsilon_1)-\frac{\sigma_{121}\sigma_{111}^{n_1}\sigma_{120}^{n_2-1}}{P_2(0,\sigma_{120})}\\
&\quad+\alpha_{11}^{n_1+1}\int_{\sigma_{120}}^{\tau_{120}}B_1(x)x^{n_2-\frac{n_1+1}{\lambda}}\frac{dx}{x}+\frac{n_1\alpha_{12}\alpha_{11}^{n_1-1}}{2}\int_{\sigma_{120}}^{\tau_{120}}A_1(x)x^{n_2-\frac{n_1}{\lambda}}\frac{dx}{x}\\
&=-\frac{\sigma_{121}\sigma_{111}^{n_1}\sigma_{120}^{n_2-1}}{P_2(0,\sigma_{120})}+\varepsilon_1^{n_2-\frac{n_1+1}{\lambda}}\varphi_7(\varepsilon_1)+\varepsilon_1^{n_2-\frac{n_1}{\lambda}}\varphi_8(\varepsilon_1)\\
&\quad-\alpha_{11}^{n_1+1}\sigma_{120}^{n_2-\frac{n_1+1}{\lambda}}\gorro{B}_1\!\left(\frac{n_1+1}{\lambda}-n_2,\sigma_{120}\right)
-\frac{n_1\alpha_{12}\alpha_{11}^{n_1-1}}{2}\sigma_{120}^{n_2-\frac{n_1}{\lambda}}\gorro{A}_1\!\left(\frac{n_1}{\lambda}-n_2,\sigma_{120}\right),
\end{align*}
where in the second equality we also use that $\alpha_{11}$ and $\alpha_{12}$ do not depend on $\varepsilon,$ see \refc{Aeq25} and \refc{Aeq19}, respectively. Notice that
\[
T_{n_1+1,0}=T_{n_1+1,0}^-=T^1_{n_1+1}+T_1^0R_{11}^{n_1}\Upsilon_1^{[n_1]}=
T^1_{n_1+1}+n_1T_1^0R_{11}^{n_1}\frac{R_{12}}{R_{11}}=
T^1_{n_1+1}+n_1T_1^0R_{11}^{n_1}S_1,
\]
where in the first equality we use \refc{Aeq22}, in the second one \refc{Aeq10} with $k=n_1+1$, in the third one the fact that $\Upsilon_1^{[n_1]}=n_1\frac{R_{12}}{R_{11}}$ from \refc{Aeq18}, and in the last one that $S_1=\frac{R_{12}}{R_{11}}=\frac{\alpha_{12}}{2\alpha_{11}}$ from \refc{Aeq23}. On account of this and using also that, from \refc{Aeq8} and \refc{Aeq16}, $T_1^0R_{11}^{n_1}=-\varepsilon_1^{n_2-\frac{n_1}{\lambda}}\frac{\alpha_{11}^{n_1}}{(n_1-\lambda n_2)P_1(0,0)}$  we get
\begin{align*}
 T_{n_1+1,0}=&-\frac{\sigma_{121}\sigma_{111}^{n_1}\sigma_{120}^{n_2-1}}{P_2(0,\sigma_{120})}+\varepsilon_1^{n_2-\frac{n_1+1}{\lambda}}\varphi_7(\varepsilon_1)+\varepsilon_1^{n_2-\frac{n_1}{\lambda}}\varphi_9(\varepsilon_1)\\
&-\alpha_{11}^{n_1+1}\sigma_{120}^{n_2-\frac{n_1+1}{\lambda}}\gorro{B}_1\!\left(\frac{n_1+1}{\lambda}-n_2,\sigma_{120}\right)
-n_1S_{1}\alpha_{11}^{n_1}\sigma_{120}^{n_2-\frac{n_1}{\lambda}}\gorro{A}_1\!\left(\frac{n_1}{\lambda}-n_2,\sigma_{120}\right)\\
=&-\frac{\sigma_{121}\sigma_{111}^{n_1}\sigma_{120}^{n_2-1}}{P_2(0,\sigma_{120})}
-\alpha_{11}^{n_1}\sigma_{120}^{n_2-\frac{n_1}{\lambda}}\left(\alpha_{11}\sigma_{120}^{\frac{-1}{\lambda}}\gorro{B}_1\!\left(\frac{n_1+1}{\lambda}-n_2,\sigma_{120}\right)
+n_1S_{1}\gorro{A}_1\!\left(\frac{n_1}{\lambda}-n_2,\sigma_{120}\right)\right).
\end{align*}
Here we also use that $\sigma_1$, $\alpha_{11}$, $T_{n_1+1,0}$ and $S_1$ do not depend on $\varepsilon$ and  apply \lemc{combinacio} to conclude that 
\[
 \varepsilon_1^{n_2-\frac{n_1+1}{\lambda}}\varphi_7(\varepsilon_1)+\varepsilon_1^{n_2-\frac{n_1}{\lambda}}\varphi_9(\varepsilon_1)=0.
\]  
Then by using the expression of $\alpha_{11}$ in \refc{Aeq25} and an easy manipulation we get that
\begin{align*}
&T_{n_1+1,0}(\np)
=-\sigma_{111}^{n_1}\sigma_{120}^{n_2}\left(\frac{\sigma_{121}}{\sigma_{120}P_2(0,\sigma_{120})}+\frac{n_1S_1}{L_1^{n_1}(\sigma_{120})}\gorro{A}_1(n_1/\lambda-n_2,\sigma_{120})\right.\\
&\hspace{5.75truecm}+\left.\frac{\sigma_{111}}{L_1^{n_1+1}(\sigma_{120})}\gorro{B}_1\big((n_1+1)/\lambda-n_2,\sigma_{120}\big)\right)
\end{align*}
for all $\np\in V.$ The application of \lemc{fun_ok} and \teoc{L8} shows that the function on the right hand side is $\cc^\infty$ in a neighbourhood of any $(\lambda_\star,\mu_\star)\in\hat W$ such that $\left\{\frac{1}{\lambda_\star},\frac{n_1}{\lambda_\star}-n_2,\frac{n_1+1}{\lambda_\star}-n_2\right\}\!\cap\Z_{\geq 0}=\emptyset,$ i.e., 
\[
 \lambda_\star\notin\tilde D_{n_1+1,0}^n\!:=\frac{1}{\N}\cup\frac{n_1}{\N_{\geq n_2}}\cup\frac{n_1+1}{\N_{\geq n_2}}.
\]
Since $D_{n_1+1,0}^n=\bigcup_{i=1}^{n_1+1}\frac{i}{\N_{\ge n_2}}$, see \obsc{domains}, by continuity we can assert that the third identity in \refc{c} is true at any $\np=(\lambda,\mu)\in\hat W$ with $\lambda\notin D_{n_1+1,0}^n\cup\tilde D_{n_1+1,0}^n=D_{n_1+1,0}^n\cup\left\{\frac{1}{k};\,k=1,2,\ldots,\lceil\frac{n_2}{n_1+1}\rceil-1\right\}.$

We begin at this point the computation of the coefficient $T_{n_2+1}.$ To this aim we apply \lemc{L3} using in this case the second column in \tabc{atlas} and the expressions of $R_{21}$, $\rho_{21}$ and $\rho_{22}$. We thus obtain
\begin{align*}
T_{n_2+1}^2=&\,\underbrace{\alpha_{21}^{n_2+1}}_{\varepsilon_2^{\lambda(n_2+1)}}\sigma_{211}\sigma_{221}^{n_2}\left(\frac{L_2(\sigma_{210})}{\sigma_{221}\sigma_{210}^\lambda}\right)^{n_2+1}\frac{\sigma_{210}^{n_1-1}}{P_1(\sigma_{210},0)}
-\underbrace{\frac{\tau_{211}\tau_{221}^{n_2}\tau_{210}^{n_1-1}}{P_1(\tau_{210},0)}}_{\varepsilon_2^{n_1}\varphi_{10}(\varepsilon_2)}\\
&+\frac{1}{2}\underbrace{\alpha_{21}^{n_2+1}}_{\varepsilon_2^{\lambda(n_2+1)}}\int_{\tau_{210}}^{\sigma_{210}}x^{-\lambda(n_2-1)}L_2^{n_2-1}(x)
\Bigg(
n_2\bigg(\alpha_{21}^{-2}\underbrace{\alpha_{22}}_{0}x^{-\lambda}L_2(x)+2x^{-2\lambda}L_2(x)\gorro{M_2}(\lambda,x)\bigg)\frac{x^{n_1-1}}{P_1(x,0)}\\[-8pt]
&\hspace{6truecm}+2x^{-2\lambda}L_2^2(x)x^{n_1-1}\partial_2P_1^{-1}(x,0)
\Bigg)dx,
\end{align*}
where we use that $\alpha_{21}=\varepsilon_2^\lambda$ from \refc{Aeq17}, $\alpha_{22}=0$ from \refc{Aeq24} and the fact that $\tau_{210}$ and $\tau_{211}$ vanish at $\varepsilon_2=0$. Notice on the other hand that, by using \refc{Aeq20}, \refc{Aeq21} and \refc{Aeq22}, 
\[
 T_{0,n_2+1}=T_{0,n_2+1}^+=\left(T_{n_2+1}^2\varepsilon_2^{-\lambda(n_2+1)}\right)(\varepsilon_1R_{11}^\lambda)^{n_2+1},
\]
which in particular shows that $T_{n_2+1}^2\varepsilon_2^{-\lambda(n_2+1)}$ does not depend on $\varepsilon.$ Having said this, note that
\begin{align*}
T_{n_2+1}^2\varepsilon_2^{-\lambda(n_2+1)}=&\,
\frac{\sigma_{211}\sigma_{210}^{n_1-1-\lambda(n_2+1)}}{\sigma_{221}}\frac{L_2^{n_2+1}(\sigma_{210})}{P_1(\sigma_{210},0)}
+\varepsilon_2^{n_1-\lambda(n_2+1)}\varphi_{10}(\varepsilon_2)
\\
&+\int_{\tau_{210}}^{\sigma_{210}}\bigg(
\underbrace{n_2\frac{L_2^{n_2}(x)}{P_1(x,0)}\gorro{M}_2(\lambda,x)+L_2^{n_2+1}(x)\partial_2P_1^{-1}(x,0)}_{B_2(x)}\bigg)x^{n_1-\lambda(n_2+1)}
\frac{dx}{x}\\
=&\,
\frac{\sigma_{211}\sigma_{210}^{n_1-1-\lambda(n_2+1)}}{\sigma_{221}}\frac{L_2^{n_2+1}(\sigma_{210})}{P_1(\sigma_{210},0)}+\varepsilon_2^{n_1-\lambda(n_2+1)}\varphi_{10}(\varepsilon_2)
\\
&+\sigma_{210}^{n_1-\lambda(n_2+1)}\gorro{B}_2(\lambda(n_2+1)-n_1,\sigma_{210})
-\underbrace{\tau_{210}^{n_1-\lambda(n_2+1)}\gorro{B}_2(\lambda(n_2+1)-n_1,\tau_{210})}_{\varepsilon_2^{n_1-\lambda(n_2+1)}\varphi_{11}(\varepsilon_2)}\\
=&\,
\frac{\sigma_{211}\sigma_{210}^{n_1-1-\lambda(n_2+1)}}{\sigma_{221}}\frac{L_2^{n_2+1}(\sigma_{210})}{P_1(\sigma_{210},0)}+\sigma_{210}^{n_1-\lambda(n_2+1)}\gorro{B}_2(\lambda(n_2+1)-n_1,\sigma_{210}),
\end{align*}
where in the second equality we apply \teoc{L8} and in the third one we take advantage of the fact that $T_{n_2+1}^2\varepsilon_2^{-\lambda(n_2+1)}$ and $\sigma_2$ do not depend on $\varepsilon$ to conclude, 
thanks to $\lambda\approx\lambda_0\notin\Q,$ that $\varphi_{10}=\varphi_{11}.$ Hence, due to
$\varepsilon_1R_{11}^\lambda=\frac{\sigma_{111}^\lambda\sigma_{120}}{L_1^\lambda(\sigma_{120})}$ by the second equality in \refc{Aeq16}, we get that
\begin{align*}
T_{0,n_2+1}=&\left(T_{n_2+1}^2\varepsilon_2^{-\lambda(n_2+1)}\right)(\varepsilon_1R_{11}^\lambda)^{n_2+1}\\[3pt]
=&\left(\frac{\sigma_{111}^\lambda\sigma_{120}}{L_1^\lambda(\sigma_{120})}\right)^{n_2+1}
\left(
\frac{\sigma_{211}\sigma_{210}^{n_1-1-\lambda(n_2+1)}}{\sigma_{221}}\frac{L_2^{n_2+1}(\sigma_{210})}{P_1(\sigma_{210},0)}+\sigma_{210}^{n_1-\lambda(n_2+1)}\gorro{B}_2(\lambda(n_2+1)-n_1,\sigma_{210})
\right).
\end{align*}
From here, taking the expression of $\Delta_{00}$ into account, we can assert that
\[
T_{0,n_2+1}(\np)
=\Delta_{00}^{n_2+1}\sigma_{210}^{n_1}\sigma_{221}^{n_2}\left(\frac{\sigma_{211}}{\sigma_{210}P_1(\sigma_{210},0)}+\frac{\sigma_{221}}{L_2^{n_2+1}(\sigma_{210})}\gorro B_2\big(\lambda(n_2+1)-n_1,\sigma_{210}\big)\right)
\]
for all $\np\in V.$ Exactly as in the previous cases, by applying \lemc{fun_ok} and \teoc{L8} it turns out that the function on the right hand side is $\cc^\infty$ on $\big((0,+\infty)\setminus\tilde D_{0,n_2+1}^n\big)\times W$ with $\tilde D_{0,n_2+1}^n\!:=\frac{\N_{\geq n_1}}{n_2+1}$. Furthermore, by \teoc{oldA} we know that the function on the left hand side is $\cc^\infty$ on $\big((0,+\infty)\setminus D_{0,n_2+1}^n\big)\times W$ where, see \obsc{domains}, $D_{0,n_2+1}^n=\frac{\N_{\geq n_1}}{n_2+1}\cup\N.$ Accordingly, due to $\tilde D_{0,n_2+1}^n\subset D_{0,n_2+1}^n,$ by continuity we can conclude that the fourth equality in \refc{c} is true on the given domain.

It only remains to compute $T_{20}$ and $T_{02}$ in the case that $n_1=0$ and $n_2=0$, respectively. Let us consider first the case $n_1=0.$ To this end we begin by computing the coefficient of $s^2$ in the time function $T^1$ of the first regular passage. By applying $(b)$ in \lemc{L3} for the case $\ell=0$ and taking $f(x_1,x_2)=\frac{x_1^{n_2-1}}{P_2(x_2,x_1)}$, see \tabc{atlas}, we know that it is written as $T_2^1=\frac{1}{2}(U_1-V_1+W_1)$ with
\begin{align*}
U_1&=(\tau_{122}R_{11}^2+\tau_{121}R_{12})f(\tau_{120},0)
    +\tau_{121}^2R_{11}^2\partial_1f(\tau_{120},0)+2\tau_{121}\tau_{111}R_{11}^2\partial_2f(\tau_{120},0)\\
&=\varepsilon_1^{n_2-1/\lambda}\varphi_{12}(\varepsilon_1)
    +\varepsilon_{1}^{n_2-2/\lambda}\varphi_{13}(\varepsilon_1),
 \\[7pt]  
V_1&=\sigma_{122}f(\sigma_{120},0)+\sigma_{121}^2\partial_1f(\sigma_{120},0) +2\sigma_{121}\sigma_{111}\partial_2f(\sigma_{120},0) \\
&=\frac{\sigma_{122}\sigma_{120}^{n_2-1}}{2P_2(0,\sigma_{120})}
+\frac{\sigma_{121}^2\sigma_{120}^{n_2-2}}{2}\left(\frac{n_2-1}{P_2(0,\sigma_{120})}+\sigma_{120}\partial_2P_2^{-1}(0,\sigma_{120})\right)
+\sigma_{121}\sigma_{111}\sigma_{120}^{n_2-1}\partial_1P_2^{-1}(0,\sigma_{120})
\intertext{and}
W_1&=\int_{\sigma_{120}}^{\tau_{120}}\left((\alpha_{11}x^{\frac{-1}{\lambda}}L_1(x))^ 2\partial_{2}^2f(x,0)+\big(\alpha_{12}x^{\frac{-1}{\lambda}}L_1(x)+2\alpha_{11}^2x^{\frac{-2}{\lambda}}L_1(x)\gorro{M}_1(1/\lambda,x)\big)\partial_2f(x,0)\right)dx
\\
&=\alpha_{11}^2\int_{\sigma_{120}}^{\tau_{120}}C_1(x)x^{n_2-\frac{2}{\lambda}}\frac{dx}{x}+\alpha_{12}\int_{\sigma_{120}}^{\tau_{120}}B_1(x)x^{n_2-\frac{1}{\lambda}}\frac{dx}{x}
\\
&=
\alpha_{11}^2\bigg(\underbrace{\tau_{120}^{n_2-\frac{2}{\lambda}}\gorro{C}_1(2/\lambda-n_2,\tau_{120})}_{\varepsilon_1^{n_2-2/\lambda}\varphi_{14}(\varepsilon_1)}-\sigma_{120}^{n_2-\frac{2}{\lambda}}\gorro{C}_1(2/\lambda-n_2,\sigma_{120})\bigg)
\\[-8pt]
&\hspace{5truecm}+\alpha_{12}\bigg(\underbrace{\tau_{120}^{n_2-\frac{1}{\lambda}}\gorro{B}_1(1/\lambda-n_2,\tau_{120})}_{\varepsilon_1^{n_2-1/\lambda}\varphi_{15}(\varepsilon_1)}-\sigma_{120}^{n_2-\frac{1}{\lambda}}\gorro{B}_1(1/\lambda-n_2,\sigma_{120})\bigg).
\end{align*}
Let us note that to rearrange $U_1$ we use that $R_{11}=\alpha_{11}\varepsilon_1^{-1/\lambda}$ and $R_{12}=\frac{1}{2}\alpha_{12}\varepsilon_1^{-1/\lambda}$ from \refc{Aeq16} and \refc{Aeq23}, respectively, and moreover that $\tau_{122},$ $\tau_{120}$ and $\tau_{121}$ vanish at $\varepsilon_1=0.$ On the other hand, to simplify~$W_1$ we apply \teoc{L8} and use that, in this case, $B_1(x)=L_1(x)
\partial_1P_2^{-1}(0,x)$ due to $n_1=0.$ By the same reason, using also \refc{Aeq10} and \refc{Aeq22}, we get that
\[
 T_{20}=T_{20}^-=T^1_2+T^0_1\Upsilon^{[0]}_2=T^1_2=\frac{1}{2}(U_1-V_1+W_1)
\] 
since $\Upsilon^{[0]}_2=0.$ This shows in particular that $U_1-V_1+W_1$ does not depend on $\varepsilon$ and, since this is also the case for $\alpha_{11}$ and $\alpha_{12},$ we can assert that
\[
\varepsilon_1^{n_2-1/\lambda}(\varphi_{12}(\varepsilon_1)+\alpha_{12}\varphi_{15}(\varepsilon_1))+\varepsilon_{1}^{n_2-2/\lambda}(\varphi_{13}(\varepsilon_1)+\alpha_{11}^2\varphi_{14}(\varepsilon_1))=0
\]
by applying \lemc{combinacio} and using that $\lambda\approx\lambda_0\notin\Q.$ Finally, since $\alpha_{11}=\frac{\sigma_{111}\sigma_{120}^{1/\lambda}}{L_1(\sigma_{120})}$ and $\alpha_{12}=2\alpha_{11}S_1$ by~\refc{Aeq25} and \refc{Aeq23}, respectively, we obtain that  
\begin{align*}
T_{20}(\np)=&-\frac{\sigma_{122}\sigma_{120}^{n_2-1}}{2P_2(0,\sigma_{120})}
-\frac{\sigma_{121}^2\sigma_{120}^{n_2-2}}{2}\left(\frac{n_2-1}{P_2(0,\sigma_{120})}+\sigma_{120}\partial_2P_2^{-1}(0,\sigma_{120})\right)
-\sigma_{121}\sigma_{111}\sigma_{120}^{n_2-1}\partial_1P_2^{-1}(0,\sigma_{120})
\\[5pt]
&-\frac{\sigma_{111}^2\sigma_{120}^{n_2}}{2L_1^2(\sigma_{120})}\gorro{C}_1(2/\lambda-n_2,\sigma_{120})
-S_1\frac{\sigma_{111}\sigma_{120}^{n_2}}{L_1(\sigma_{120})}\gorro{B}_1(1/\lambda-n_2,\sigma_{120})
\end{align*}
for all $\np\in V.$ By applying \lemc{fun_ok} and \teoc{L8} we have that $\gorro{C}_1(2/\lambda-n_2,\sigma_{120})$ is $\cc^\infty$ in a neighbourhood of any $(\lambda_\star,\mu_\star)\in\hat W$ such that $\left\{1/{\lambda_\star},2/{\lambda_\star}-n_2\right\}\cap\Z_{\geq 0}=\emptyset.$ The condition for the function~$S_1$, see \refc{def_S}, and $\gorro{B}_1(1/\lambda-n_2,\sigma_{120})$ is $1/\lambda_\star\notin\Z_{\geq 0}$ and $1/\lambda_\star-n_2\notin\Z_{\geq 0}$, respectively. Therefore the function on the right hand side in the above equality is $\cc^\infty$ on $\big((0,+\infty)\setminus\tilde D_{20}^n\big)\times W$ with 
$\tilde D_{20}^n\!:=\frac{1}{\N}\cup \frac{2}{\N_{\geq n_2}}$. Due to $D_{20}^n=\frac{2}{\N_{\geq n_2}}$ from \obsc{domains}, we get that $D_{20}^n\cup\tilde D_{20}^n=D_{20}^n\cup\left\{\frac{1}{k};\,k=1,2,\ldots,\lceil\frac{n_2}{2}\rceil-1\right\}$ and, on account of the considerations in the second paragraph of the proof, this shows that the above equality is true in the domain given in the statement.

Let us turn finally to the computation of $T_{02}$ for the case $n_2=0.$ Similarly as before we apply $(b)$ in \lemc{L3} with $f(x_1,x_2)=\frac{x_1^{n_1-1}}{P_1(x_1,x_2)}$ to get that $T^2_2=\frac{1}{2}(U_2-V_2+W_2)$. In this case some long but easy computations taking account of \tabc{atlas} give
\begin{align*}
U_2&=(\sigma_{212}R_{21}^2+\sigma_{211}R_{22})f(\sigma_{210},0)+\sigma_{211}^2R_{21}^2\partial_1f(\sigma_{210},0)+2\sigma_{211}\sigma_{221}R_{21}^2\partial_2f(\sigma_{210},0)\\
&=\varepsilon_2^{2\lambda}\sigma_{210}^{n_1}
\left(\frac{L_2(\sigma_{210})}{\sigma_{221}\sigma_{210}^\lambda}\right)^2\left(2Z-\frac{\sigma_{211}S_2}{\sigma_{210}P_1(\sigma_{210},0)}\right),
\end{align*}
where we use that $R_{21}=\varepsilon_2^\lambda\frac{L_2(\sigma_{210})}{\sigma_{221}\sigma_{210}^\lambda}$ from \refc{Aeq27} and \refc{Aeq17} and that $R_{22}=-\varepsilon_2^{2\lambda}S_2\left(\frac{L_2(\sigma_{210})}{\sigma_{221}\sigma_{210}^\lambda}\right)^2$ from \refc{Aeq26} and, for the sake of shortness, we denote
\[
Z\!:=\frac{\sigma_{212}\sigma_{210}+(n_1-1)\sigma_{211}^2}{2\sigma_{210}^2P_1(\sigma_{210},0)}+\frac{\sigma_{211}^2}{2\sigma_{210}}\partial_1P_1^{-1}(\sigma_{210},0)+\frac{\sigma_{211}\sigma_{221}}{\sigma_{210}}\partial_2P_1^{-1}(\sigma_{210},0).
\]
Since $\tau_{210}$, $\tau_{211}$ and $\tau_{212}$ vanish at $\varepsilon_2=0$, one can also verify that
\[
 V_2=\tau_{212}f_2(\tau_{210},0)+\tau_{211}^2\partial_1f_2(\tau_{210},0)+2\tau_{211}\tau_{221}\partial_2f_2(\tau_{210},0)
 =\varepsilon_2^{n_1}\varphi_{16}(\varepsilon_2).
\]
Furthermore, on account of the definition of the function $C_2$ given in \refc{def_fun} and applying \teoc{L8},
\begin{align*}
W_2&=\int_{\tau_{210}}^{\sigma_{210}}\left(
(\varepsilon_2^\lambda x^{-\lambda}L_2(x))^2x^{n_1}\partial_2^2P_1^{-1}(x,0)+2\varepsilon_2^{2\lambda}x^{n_1-2\lambda}L_2(x)\gorro{M}_2(\lambda,x)\partial_2P_1^{-1}(x,0)
\right)
\frac{dx}{x}\\[10pt]
&=\varepsilon_2^{2\lambda}\int_{\tau_{210}}^{\sigma_{210}}C_2(x)x^{n_1-2\lambda}\frac{dx}{x}=\varepsilon_2^{2\lambda}
\bigg(
\sigma_{210}^{n_1-2\lambda}\gorro{C}_2(2\lambda-n_1,\sigma_{210})
-\underbrace{\tau_{210}^{n_1-2\lambda}\gorro{C}_2(2\lambda-n_1,\tau_{210})}_{\varepsilon_2^{n_1-2\lambda}\varphi_{17}(\varepsilon_2)}
\bigg).
\end{align*}
Notice at this point that, from \refc{Aeq20}, \refc{Aeq21} and \refc{Aeq22}, $T_{02}=T_{02}^+=(T_2^2\varepsilon_2^{-2\lambda})(\varepsilon_1R_{11}^\lambda)^2$, which shows in particular that $T_2^2\varepsilon_2^{-2\lambda}$ does not depend on $\varepsilon$ because this is the case for $T_{02}$ and, see \refc{Aeq16}, $\varepsilon_1R_{11}^\lambda=\alpha_{11}$. Consequently $U_2-V_2+W_2$ does not depend on $\varepsilon$ and so $\varepsilon_2^{n_1-2\lambda}(\varphi_{16}(\varepsilon_2)-\varphi_{17}(\varepsilon_2))=c$. Since $\lambda\approx\lambda_0\notin\Q$, this implies that $\varphi_{16}=\varphi_{17}$ and therefore
\begin{align*}
 T_{02}&=\left(\frac{\sigma_{111}^\lambda\sigma_{120}}{L_1^\lambda(\sigma_{210})}\right)^2
 \left(
  \sigma_{210}^{n_1}
\left(\frac{L_2(\sigma_{210})}{\sigma_{221}\sigma_{210}^\lambda}\right)^2\left(Z-\frac{\sigma_{211}S_2}{2\sigma_{210}P_1(\sigma_{210},0)}\right)
+\frac{1}{2}\sigma_{210}^{n_1-2\lambda}\gorro{C}_2(2\lambda-n_1,\sigma_{210})
 \right)\\
 &=\Delta_{00}^2\sigma_{210}^{n_1}\left(
 Z-\frac{\sigma_{211}S_2}{2\sigma_{210}P_1(\sigma_{210},0)}
 +\frac{\sigma_{221}^2}{2L_2^2(\sigma_{210})}\gorro{C}_2(2\lambda-n_1,\sigma_{210})
 \right).
\end{align*}
for all $\np\in V.$ Exactly as before, by applying \lemc{fun_ok} and \teoc{L8} we can assert that $\gorro{C}_2(2\lambda-n_1,\sigma_{210})$ is $\cc^\infty$ in a neighbourhood of any $(\lambda_\star,\mu_\star)\in\hat W$ such that $\left\{\lambda_\star,2\lambda_\star-n_1\right\}\cap\Z_{\geq 0}=\emptyset.$ The corresponding condition for the function~$S_2$, see \refc{def_S}, is $\lambda_\star\notin\Z_{\geq 0}$. Thus
the function on the right hand side in the above equality is $\cc^\infty$ on $\big((0,+\infty)\setminus\tilde D_{02}^n\big)\times W$ with $\tilde D_{02}^n\!:=\N\cup \frac{\N_{\geq n_1}}{2}$. Due to $D_{02}^n=\frac{\N}{2}$ from \obsc{domains}, it turns out that $D_{20}^n\cup\tilde D_{20}^n=D_{20}^n$ and, on account of the considerations in the second paragraph of the proof, this shows that the above equality is true in the domain given in the statement.
This concludes the proof of the result.
\end{prooftext}

\begin{lem}\label{toma}
Let $\Phi(x,y)$, with $x=(x_1,x_2,\ldots,x_n)\in\R^N$ and $y\in\R$, be a continuous function in a neighbourhood of $(0,0)\in\R^N\times\R.$ If $y\Phi(x,y)$ is analytic in a neighbourhood of $(0,0)$ then $\Phi(x,y)$ is analytic in a neighbourhood of $(0,0)$.
\end{lem}

\begin{prova}
By the Weierstrass Division Theorem (see \cite[Theorem 1.8]{Greuel} or \cite[Theorem 6.1.3]{Krantz}) there exist a neighbourhood $U$ of $0\in\R^N$ and an open interval $I$ containing $y=0$ such that $y\Phi(x,y)=yg(x,y)+r(x)$ with $g\in\cc^\omega(U\times I)$ and $r\in\cc^\omega(I)$. The evaluation of this equality at $y=0$ yields $r\equiv 0.$ Consequently $\Phi(x,y)=g(x,y)$ for all $(x,y)\in U\times (I\setminus\{0\})$ and, by the continuity of $\Phi$ in a neighbourhood of $(0,0)$, we easily get $\Phi\equiv g$ on $U\times I$. This proves the result because $g\in\cc^\omega(U\times I)$. 
\end{prova}

\begin{prop}\label{pre-analitico}
In the analytic setting $($see \obsc{analytic_setting}$)$, the following assertions hold:
\begin{enumerate}[$(a)$]
\item The coefficient $\Delta_{ij}$ of the Dulac map is $\cc^\omega$ on 
$((0,+\infty)\setminus D_{ij}^0)\times W$ for $(i,j)\in\{(0,0),(1,0),(0,1),(1,1)\}$.

\item For each $(i,j)\in\{(n_1,0),(0,n_2),(n_1+1,0),(0,n_2+1)\}$, the coefficient $T_{ij}$ of the Dulac time is analytic on 
$((0,+\infty)\setminus D_{ij}^n)\times W$. This is also the case for $(i,j)=(2,0)$ and $(i,j)=(0,2)$ assuming $n_1=0$ and $n_2=0$, respectively. 

\end{enumerate}
\end{prop}

\begin{prova}
By applying \lemc{fun_ok} we know that, for $i=1,2$, the functions $L_i(u;\np)$, $M_i(u;\np)$ and $A_i(u;\np)$ given in \refc{def_fun} are analytic on $I_i\times\hat W$. In addition, 
\begin{itemize}
\item the functions $B_1(u;\np)$ and $C_1(u;\np)$ are analytic on 
         $I_1\times((0,+\infty)\setminus\frac{1}{\N})\times W$, and
\item the functions $B_2(u;\np)$ and $C_2(u;\np)$ are analytic on $I_2\times((0,+\infty)\setminus\N)\times W$.
\end{itemize}
Moreover, since the parametrization $\sigma_i(s;\np)$ of the transverse section $\Sigma_i$ is analytic by assumption for $i=1,2$, from \refc{def_S} we get that $S_1(\lambda,\mu)$ and $S_2(\lambda,\mu)$ are analytic on $((0,+\infty)\setminus\frac{1}{\N})\times W$ and $((0,+\infty)\setminus\N)\times W,$ respectively. 

The fact that each coefficient $\Delta_{ij}(\lambda,\mu)$ in assertion $(b)$ of \teoc{A} is analytic on $((0,+\infty)\setminus D_{ij}^0)\times W$ follows readily from regularity properties stated in the previous paragraph because, see \obsc{domains}, 
\[
D_{00}^0=\emptyset,\; D_{01}^0=\N,\; D_{10}^0=\frac{1}{\N}\text{ and }D_{11}^0=\N\cup\frac{1}{\N}.
\]
%Taking this into account, that each coefficient $\Delta_{ij}(\lambda,\mu)$ is meromorphic on $\hat W=(0,+\infty)\times W$ with poles of order at most $2$ along $\lambda\in D_{ij}^0$ follows by \propc{poles1} and \obsc{analitico2}. 
This proves assertion~$(a)$.

By the first assertion in $(d)$ of \teoc{L8}, the regularity properties established in the first paragraph also imply that each coefficient $T_{ij}(\lambda,\mu)$ listed in $(c)$ of \teoc{A} is analytic on $((0,+\infty)\setminus D_{ij}^n)\times W$, with the exception of the special values 
\begin{itemize}
\item $\lambda=\frac{1}{k}$ with $k\in\big\{1,2,\ldots,\lceil\frac{n_2}{n_1+1}\rceil-1\big\}$ for $T_{n_1+1,0}(\lambda,\mu)$, and 
\item $\lambda=\frac{1}{k}$ with $k\in\big\{1,2,\ldots,\lceil\frac{n_2}{2}\rceil-1\big\}$ for $T_{20}(\lambda,\mu)$,
\end{itemize} 
where the respective formula does not hold. Indeed this follows using that, see \obsc{domains} again,
$D_{00}^n=\emptyset$,
\[
D_{n_1,0}^n=
\bigcup_{i=1}^{n_1}\frac{i}{\N_{\geq n_2}},
\; 
D_{0,n_2}^n=\left\{\begin{array}{cl}
\frac{\N_{\geq n_1}}{n_2} & \text{ if $n_2\geqslant 1$,}\\[5pt]
\emptyset & \text{ if $n_2=0$,}
\end{array}
\right.
\;
D_{n_1+1,0}^n=\bigcup_{i=1}^{n_1+1}\frac{i}{\N_{\ge n_2}}
\text{ and }
D_{0,n_2+1}^n=\frac{\N_{\ge n_1}}{n_2+1}\cup\N,
\]
together with $D_{20}^n=\frac{2}{\N_{\geq n_2}}$ for $n_1=0$ and 
$D_{02}^n=\frac{\N}{2}$ for $n_2=0$.
For instance, due to $A_2(u;\np)\in\cc^\omega(I_2\times\hat W)$, the first assertion in $(d)$ of \teoc{L8} implies that $\hat A_2(\alpha,u;\np)$ is analytic on $(\R\setminus\Z_{\geq 0})\times I_2\times\hat W$ and hence
\[
 T_{0,n_2}(\np)=\Delta_{00}^{n_2}\frac{\sigma_{210}^{n_1}\sigma_{221}^{n_2}}{L_2^{n_2}(\sigma_{210})}
 \gorro A_2(n_2\lambda-n_1,\sigma_{210})
\]
is analytic at $\lambda=\lambda_0$ provided that $n_2\lambda_0-n_1\notin\Z_{\geq 0},$ i.e., $\lambda_0\notin D_{0,n_2}^n$. The analysis of the other coefficients follows similarly and the details are omitted for the sake of brevity. 

So let us focus on the analyticity of $T_{n_1+1,0}$ and $T_{20}$ at the special values listed above. In order to study the first case let us fix $\lambda_0=\frac{1}{k}$ with $k\in\{1,\ldots,\lceil\frac{n_2}{n_1+1}\rceil-1\}$. Note that we can write, see $(c)$ in \teoc{A},
\begin{equation}\label{corBeq1}
 T_{n_1+1,0}=f_0+f_1S_1\hat A_1(n_1/\lambda-n_2,\sigma_{120})+f_2\hat B_1\big((n_1+1)/\lambda-n_2,\sigma_{120}\big)
\end{equation}
where, see \refc{def_fun}, $B_1(u)=g_1(u)\hat M_1(1/\lambda,u)+g_2(u)$ and $S_1=f_3+f_4\hat M_1(1/\lambda,\sigma_{120})$ with $g_i(u;\np)\in\cc^\omega(I_1\times\hat W)$ and $f_i(\np)\in\cc^\omega(\hat W)$. That being said we argue as follows:
\begin{enumerate}

\item $\hat A_1(n_1/\lambda-n_2,\sigma_{120})$ is analytic at $\lambda=\lambda_0$ due to $\frac{n_1}{\lambda_0}-n_2=n_1k-n_2\in\Z_{< 0}$ by the first assertion in~$(d)$ of \teoc{L8}.

\item $(\lambda-\lambda_0)\gorro M_1(1/\lambda,u;\np)$, and consequently $(\lambda-\lambda_0)B_1(u;\np)$ and $(\lambda-\lambda_0)S_1(\np)$, extends analytically at $\lambda=\lambda_0$ by the second assertion in $(d)$ of \teoc{L8} since $1/\lambda_0=k\in\Z_{\geq 0},$ 

\item and this implies (in this case by applying the first assertion) that $(\lambda-\lambda_0)\gorro B_1((n_1+1)/\lambda-n_2,\sigma_{120})$ extends analytically at $\lambda=\lambda_0$ because $\frac{n_1+1}{\lambda_0}-n_2=(n_1+1)k-n_2\in\Z_{< 0}.$

\end{enumerate}
Taking this into account, from \refc{corBeq1} it follows readily that $(\lambda-\lambda_0)T_{n_1+1,0}(\np)$ extends analytically at $\lambda=\lambda_0$. On the other hand, since $\lambda_0\notin D_{n_1+1,0}^n$, note that $T_{n_1+1,0}(\np)$ is smooth at $\lambda=\lambda_0$ by $(b)$ in \teoc{oldA}. Accordingly, in view of \lemc{toma}, we can assert that $T_{n_1+1,0}(\np)$ is analytic at $\lambda=\lambda_0$ as desired.

Let us turn next to the second case. So let us fix $\lambda_0=\frac{1}{k}$ with $k\in\{1,\ldots,\lceil\frac{n_2}{2}\rceil-1\}$
and observe that from $(c)$ in \teoc{A} we get that if $n_1=0$ then we can write 
\begin{equation}\label{corBeq2}
 T_{20}=f_0+f_1\hat C_1(2/\lambda-n_2,\sigma_{120})+f_2S_1\hat B_1(1/\lambda-n_2,\sigma_{120})
\end{equation}
with, see \refc{def_fun}, $C_1(u)=B_1(u)\big(L_1(u)+2\hat M_1(1/\lambda,u)\big)$ and $S_1=f_3+f_4\hat M_1(1/\lambda,\sigma_{120})$ for some $f_i\in\cc^\omega(\hat W)$. We point out that in this case, since $n_1=0,$ $B_1(u)=L_1(u)\partial_1P_2^{-1}(0,u)$ is analytic on $I_1\times\hat W.$ Then we proceed as follows:
\begin{enumerate}

\item $\hat B_1(1/\lambda-n_2,\sigma_{120})$ is analytic at $\lambda=\lambda_0$ due to $1/\lambda_0-n_2=k-n_2\in\Z_{<0}$ by the first assertion in $(d)$ of \teoc{L8}. 

\item $(\lambda-\lambda_0)\gorro M_1(1/\lambda,u;\np)$ extends analytically at $\lambda=\lambda_0$ by the second assertion in $(d)$ of \teoc{L8} because $1/\lambda_0=k\in\Z_{\geq 0}$ and, 

\item consequently, this is so for $(\lambda-\lambda_0)S_1(\np)$ and $(\lambda-\lambda_0)\hat C_1(2/\lambda-n_2,\sigma_{120}),$ the latter by the first assertion in $(d)$ of \teoc{L8} since $2/\lambda_0-n_2=2k-n_2\in\Z_{<0}.$   

\end{enumerate}
On account of this, from \refc{corBeq2} we get that $(\lambda-\lambda_0)T_{20}(\np)$ extends analytically at $\lambda=\lambda_0.$ Exactly as before, it happens that $T_{20}(\np)$ is smooth at $\lambda=\lambda_0$ by $(b)$ in \teoc{oldA} due to $\lambda_0\notin D_{20}^n.$ Therefore, by \lemc{toma} again, we can assert that $T_{20}(\np)$ is analytic at $\lambda=\lambda_0$ as desired. This proves the validity of $(b)$.
\end{prova}

%So far we have proved that each $T_{ij}(\lambda,\mu)$ in assertion $(c)$ of \teoc{A} is analytic on $((0,+\infty)\setminus D_{ij}^n)\times W$. Finally the fact that each coefficient is meromorphic on $\hat W=(0,+\infty)\times W$ with poles of order at most $2$ along $\lambda\in D_{ij}^n$ follows, recall \obsc{analitico2}, by applying Propositions~\ref{poles2}, \ref{poles3}, \ref{poles4} and \ref{poles5} in the corresponding case. This proves assertion $(b)$ and concludes the proof of the result.

\section{Poles and residues of the coefficients}\label{sec:poles}

Let us recall, see \teoc{oldA}, that the coefficient $\Delta_{ij}(\lambda,\mu)$ of the Dulac map is $\mathscr C^\infty$ on $((0,+\infty)\setminus D_{ij}^0)\times W$ for each $(i,j)\in\Lambda_0$ and the coefficient $T_{ij}(\lambda,\mu)$ of the Dulac time is $\mathscr C^\infty$ on $((0,+\infty)\setminus D_{ij}^n)\times W$ for each $(i,j)\in\Lambda_n.$ The next result is addressed to the behaviour of these coefficients at the boundaries of their respective domains of definition. 

\begin{lem}\label{polos}
Consider the coefficients $\Delta_{ij}$ and $T_{ij}$ of the Dulac map and the Dulac time, respectively, given by \teoc{oldA}. The following assertions hold:
\begin{enumerate}[$(a)$]
         
\item If $(i,j)\in\Lambda_0$ and $\lambda_0\in D_{ij}^0$ then there exists $\ell\in\Z_{\geq 0}$ such that the function
         $\np\mapsto (\lambda-\lambda_0)^\ell\Delta_{ij}(\np)$ 
         extends $\cc^\infty$ to $\{\lambda_0\}\times W$.
         
\item If $(i,j)\in\Lambda_n$ and $\lambda_0\in D_{ij}^n$ then there exists $\ell\in\Z_{\geq 0}$ such that the function
         $\np\mapsto (\lambda-\lambda_0)^\ell T_{ij}(\np)$ 
         extends $\cc^\infty$ to $\{\lambda_0\}\times W$.         

\end{enumerate}
Moreover, setting $\lambda_0=p/q$ with $\gcd(p,q)=1$, the estimates $\ell\leqslant \frac{i}{p}+\frac{j}{q}\leqslant i+j$ hold in both cases.
\end{lem}

\begin{prova}
For convenience we prove $(b)$ first. Due to $\lambda_0\in D_{ij}^n$, we have $\lambda_0\in\Q$ and we write $\lambda_0=p/q$ with $\gcd(p,q)=1.$ Setting $r_{n}\!:=\max\{r\in\Z_{\geq 0}:(i,j)+r(p,-q)\in\Lambda_n\}$, we define $(i_{n},j_{n})=(i,j)+r_{n}(p,-q)$. Then $\lambda_0\in D_{i_{n},j_{n}}^n$, $\mathscr A_{i_{n} j_{n}\lambda_0}^n\neq\emptyset$, see \defic{alldefi}, and we take $\ell\!:=\max\mathscr A_{i_{n} j_{n}\lambda_0}^n.$
By~$(b2)$ in \teoc{oldA} we know that $\T_{i_{n},j_{n}}^{\lambda_0}(w;\np)\in\cc^\infty(\hat U)[w]$, where $\hat U$ is an open neighbourhood of $\{\lambda_0\}\times W,$ and 
 \[
\T_{i_{n} j_{n}}^{\lambda_0}(w;\np)=\sum_{r\in\mathscr A_{i_{n} j_{n}\lambda_0}^n}T_{i_{n}-rp,j_{n}+rq}(\np)(1+\alpha w)^r\text{ for $\lambda\neq\lambda_0$,}
\]
where $\alpha=p-\lambda q.$ Let us write $\T_{i_{n} j_{n}}^{\lambda_0}(w;\np)=\sum_{k=0}^\ell A_k(\np)w^k$ with $A_k\in \cc^\infty(\hat U)$. For convenience we define $u\!:=1+\alpha w,$ so that $w=\alpha^{-1}(u-1)$ for $\alpha\neq 0$. Thus $w^k=\alpha^{-k}\sum_{r=0}^k {k\choose r} (-1)^{k-r}u^r$ and, for $\lambda\neq\lambda_0,$ 
\[
 \T_{i_{n} j_{n}}^{\lambda_0}(w;\np)=\sum_{r=0}^\ell\left(\sum_{k=r}^\ell A_k(\np)\alpha^{-k}{k\choose r}(-1)^{k-r}\right)(1+\alpha w)^r.
\]
Accordingly this shows that $T_{i_{n}-rp,j_{n}+rq}(\np)=\sum_{k=r}^\ell A_k(\np)\alpha^{-k}{k\choose r}(-1)^{k-r}$ provided that $r\in \mathscr A_{i_{n} j_{n}\lambda_0}^n$ and  $\lambda\neq\lambda_0$. With regard to the first condition let us observe that $r_{n}\in\mathscr A_{i_{n} j_{n}\lambda_0}^n$ by construction. Hence $T_{i,j}(\np)=\sum_{k=r_{n}}^\ell A_k(\np)\alpha^{-k}{k\choose r_{n}}(-1)^{k-r_{n}}$ and, due to $\alpha=q(\lambda_0-\lambda),$
\[
  (\lambda-\lambda_0)^\ell T_{i,j}(\np)=(-1)^{r_{n}}\sum_{k=r_{n}}^\ell q^{-k} A_k(\np)(\lambda-\lambda_0)^{\ell-k}{k\choose r_{n}}\text{ for $\lambda\neq \lambda_0$}.
\]
Since $A_k\in \cc^\infty(\hat U)$, this shows that $\np\mapsto (\lambda-\lambda_0)^\ell T_{ij}(\np)$ 
extends $\cc^\infty$ to $\{\lambda_0\}\times W$ and proves $(b)$.

The proof of $(a)$ follows verbatim replacing $n=(n_1,n_2)$ by $0=(0,0)$ and is omitted for the sake of shortness. Let us turn now to the proof of the last assertion in the statement. The estimate for the the case in $(a)$, i.e., $(i,j)\in\Lambda_0$ and $\lambda_0\in D_{ij}^0,$ is clear because
\[
 \max\mathscr A_{i_{0} j_{0}\lambda_0}^0\leqslant \frac{i_0}{p}=\frac{i}{p}+r_0\leqslant\frac{i}{p}+\frac{j}{q}\leqslant i+j.
\]
Here the first inequality follows using that $\mathscr A_{i_{0} j_{0}\lambda_0}^0\neq\emptyset$ and $(i_0-rp,j_0+rq)\in\Lambda_0=\Z_{\geq 0}\times\Z_{\geq 0}$ for all $r\in \mathscr A_{i_{0} j_{0}\lambda_0}^0$, see \defic{alldefi}, the equality is due to $(i_0,j_0)\!:=(i,j)+r_0(p,-q)$, the second inequality is a consequence of $j-r_0q=j_0\geqslant 0$ and the third inequality is evident since $p,q\in\N.$ Finally, the estimate for the case in $(b)$, i.e., $(i,j)\in\Lambda_n$ and $\lambda_0\in D_{ij}^n,$ is a consequence of the previous discussion and the fact that, by construction, $\mathscr A_{i_{n} j_{n}\lambda_0}^n\neq\emptyset$ and $\max\mathscr A_{i_{n} j_{n}\lambda_0}^n\leqslant \max\mathscr A_{i_{0} j_{0}\lambda_0}^0$. This completes the proof of the result.
\end{prova}

By \lemc{polos} the coefficients $\Delta_{ij}$ and $T_{ij}$ have poles at $D_{ij}^0\times W$ and $D_{ij}^n\times W$, respectively, of order at most $i+j.$ This is a general result, meaning that it holds for any $(i,j).$  
\teoc{A} provides the explicit expression of some of these coefficients and the rest of the present section is devoted to give sharps bounds for the order of their poles. We will also compute the residues of these coefficients at their poles, which determine the values of the leading terms of the polynomials $\boldsymbol{\Delta}_{ij}^{\lambda_0}(\omega;\np)$ at $\lambda_0\in D_{ij}^0$ and  $\T_{ij}^{\lambda_0}(\omega;\np)$ at $\lambda_0\in D_{ij}^n$
(see \teoc{3punts} and \teoc{9punts}, respectively, in \secc{quarta}). We illustrate the use of the residues for this purpose in \exc{ex1}. Let us also advance that at the end of the section we will finish the proof of \coryc{analitico}, which shows that in the analytic setting these coefficients are meromorphic on $\hat W=(0,+\infty)\times W$.

With regard to the next statement we recall that $D_{01}^0=\N$, $D_{10}^0=\frac{1}{\N}$ and $D_{11}^0=\N\cup\frac{1}{\N}$ (see \obsc{domains}).

\begin{prop}\label{poles1}
The following assertions hold:
\begin{enumerate}[$(a)$]

\item For any $\np_0=(\lambda_0,\mu_0)\in  D_{10}^0\times W$, the function $\np\mapsto (\lambda-\lambda_0)\Delta_{10}(\np)$ extends $\cc^\infty$ at $\np=\np_0$, and if $\lambda_0=\frac{1}{i}$ with $i\in\N$ then 
$\lim\limits_{\np\to\np_0}(\lambda-\lambda_0)\Delta_{10}(\np)=-\frac{\Delta_{00}\sigma_{111}\sigma_{120}^i}{L_1(\sigma_{120})i^3}\frac{M_1^{(i)}(0)}{i!}\big|_{\np=\np_0}.$

\item For any $\np_0=(\lambda_0,\mu_0)\in  D_{01}^0\times W$, the function $\np\mapsto (\lambda-\lambda_0)\Delta_{01}(\np)$ extends $\cc^\infty$ at $\np=\np_0$, and if $\lambda_0=i\in\N$ then $\lim\limits_{\np\to\np_0}(\lambda-\lambda_0)\Delta_{01}(\np)=-\frac{\Delta_{00}^2\sigma_{221}\sigma_{210}^i}{L_2(\sigma_{210})}\frac{M_2^{(i)}(0)}{i!}\big|_{\np=\np_0}.$

\item For any $\np_0=(\lambda_0,\mu_0)\in  (D_{11}^0\setminus\{1\})\times W$, the function $\np\mapsto (\lambda-\lambda_0)\Delta_{11}(\np)$ extends $\cc^\infty$ at $\np=\np_0$ and

\begin{enumerate}[$(c1)$]

\item if $\lambda_0=\frac{1}{i}$ with $i\in\N_{\geq 2}$ then $\lim\limits_{\np\to\np_0}(\lambda-\lambda_0)\Delta_{11}(\np)=\frac{2\Delta_{00}^2\sigma_{111}\sigma_{120}^i}{L_1(\sigma_{120})i^3}\frac{M_1^{(i)}(0)}{i!}S_2\big|_{\np=\np_0},$

\item if $\lambda_0=i\in\N_{\geq 2}$ then $\lim\limits_{\np\to\np_0}(\lambda-\lambda_0)\Delta_{11}(\np)=-\frac{2i\Delta_{00}^2\sigma_{221}\sigma_{210}^i}{L_2(\sigma_{210})}\frac{M_2^{(i)}(0)}{i!}S_1\big|_{\np=\np_0}$.

\end{enumerate}

Finally, for any $\np_0=(\lambda_0,\mu_0)\in  \{1\}\times W$, the function $\np\mapsto (\lambda-\lambda_0)^2\Delta_{11}(\np)$ extends $\cc^\infty$ at $\np=\np_0$ and 
$\lim\limits_{\np\to\np_0}(\lambda-\lambda_0)^2\Delta_{11}(\np)=2\Delta_{00}^2\frac{\sigma_{111}\sigma_{120}M_1'(0)}{L_1(\sigma_{120}))}
\frac{\sigma_{221}\sigma_{210}M_2'(0)}{L_2(\sigma_{210})}\big|_{\np=\np_0}$.
%$\lim\limits_{\np\to\np_0}(\lambda-\lambda_0)^2\Delta_{11}(\np)=\frac{2\Delta_{00}^2\sigma_{111}\sigma_{120}\sigma_{221}\sigma_{210}}{L_1(\sigma_{120})L_2(\sigma_{210})}M_1'(0)M_2'(0)\big|_{\np=\np_0}$.

\end{enumerate}
\end{prop}

\begin{prova}
In order to show $(a)$ we fix $\np_0=(1/i,\mu_0)\in D_{10}^0\times W$ with $i\in\N$ and note that, by $(b)$ in \teoc{A}, $\Delta_{10}=\Delta_{00}\lambda S_1$ where $\Delta_{00}\in\cc^\infty(\hat W)$ and,  see \refc{def_S}, $S_1=f_1-\frac{\sigma_{111}}{L_1(\sigma_{120})}\gorro{M}_1(1/\lambda,\sigma_{120})$ with $f_1\in\cc^\infty(\hat W)$.
On account of this and $(c)$ in \teoc{L8}, the function $(\lambda-1/i)\Delta_{10}(\np)$ extends $\cc^\infty$ at $\np=\np_0$ and 
\begin{equation}\label{42eq1}
 \lim_{\np\to\np_0}(\lambda-1/i)S_1=\frac{-\sigma_{111}}{L_1(\sigma_{120})}\Big|_{\np=\np_0}
 \lim_{\np\to\np_0}\frac{i-1/\lambda}{i/\lambda}\gorro{M}_1(1/\lambda,\sigma_{120})
 =\frac{-\sigma_{111}}{L_1(\sigma_{120})i^2}\frac{M_1^{(i)}(0)}{i!}\sigma_{120}^i\Big|_{\np=\np_0}.
\end{equation}
Therefore $\lim_{\np\to\np_0}(\lambda-1/i)\Delta_{10}(\np)=-\frac{\Delta_{00}\sigma_{111}\sigma_{120}^i}{L_1(\sigma_{120})i^3}\frac{M_1^{(i)}(0)}{i!}\big|_{\np=\np_0}$. 

To prove $(b)$ we fix $\np_0=(i,\mu_0)\in D_{01}^0\times W$ with $i\in\N$ and note that, by $(b)$ in \teoc{A}, $\Delta_{01}=-\Delta_{00}^2S_2$ where $S_2=f_2-{\frac{\sigma_{221}}{L_2(\sigma_{210})}}\gorro{M}_2(\lambda,\sigma_{210})$ with $f_2\in\cc^{\infty}(\hat W).$ Exactly as before, $(c)$ in \teoc{L8} implies that the function $(\lambda-i)\Delta_{01}(\np)$ extends $\cc^\infty$ at $\np=\np_0$ and, moreover, that
\begin{equation}\label{42eq2}
 \lim_{\np\to\np_0}(\lambda-i)S_2=\frac{\sigma_{221}}{L_2(\sigma_{210})}\Big|_{\np=\np_0}\lim_{\np\to\np_0}(i-\lambda)\gorro{M}_2(\lambda,\sigma_{210})
 =\frac{\sigma_{221}}{L_2(\sigma_{210})}\frac{M_2^{(i)}(0)}{i!}\sigma_{210}^i\Big|_{\np=\np_0}
\end{equation}
and, consequently, $\lim_{\np\to\np_0}(\lambda-i)\Delta_{01}(\np)=-\frac{\Delta_{00}^2\sigma_{221}\sigma_{210}^i}{L_2(\sigma_{210})}\frac{M_2^{(i)}(0)}{i!}\big|_{\np=\np_0}$.

Let us turn to the proof of $(c)$. To this end we note that, by $(b)$ in \teoc{A}, $\Delta_{11}=-2\Delta_{00}^2\lambda S_1S_2.$ If $\np_0=(1/i,\np_0)\in D_{11}^0\times W$ with $i\in\N_{\geq 2}$ then $S_2$ is smooth at $\np=\np_0$ by $(a)$ in \teoc{L8} and therefore from~\refc{42eq1} it follows that
\[
 \lim_{\np\to\np_0}(\lambda-1/i)\Delta_{11}(\np)=\frac{2\Delta_{00}^2\sigma_{111}\sigma_{120}^i}{L_1(\sigma_{120})i^3}\frac{M_1^{(i)}(0)}{i!}S_2\Big|_{\np=\np_0}.
\]
Exactly as before, the fact that $(\lambda-1/i)\Delta_{11}(\np)$ extends $\cc^{\infty}$ at $\np=\np_0$ follows by $(c)$ in \teoc{L8}. This shows the assertion in $(c1)$. Similarly if $\np_0=(i,\np_0)\in D_{11}^0\times W$ with $i\in\N_{\geq 2}$ then $S_1$ is smooth at $\np=\np_0$ by $(a)$ in \teoc{L8} and, from~\refc{42eq2}, 
\[
\lim_{\np\to\np_0}(\lambda-i)\Delta_{11}(\np)=-\frac{2i\Delta_{00}^2\sigma_{221}\sigma_{210}^i}{L_2(\sigma_{210})}\frac{M_2^{(i)}(0)}{i!}S_1\Big|_{\np=\np_0}
\]   
which proves $(c2)$. Finally, if $\np_0=(1,\mu_0)$ with $\mu_0\in W$, the combination of \refc{42eq1} and \refc{42eq2} easily implies that
\begin{align*}
 \lim_{\np\to\np_0}(\lambda-1)^2\Delta_{11}(\np)
 &=2\Delta_{00}^2\Big|_{\np=\np_0}\lim_{\np\to\np_0}(\lambda-1)S_1\lim_{\np\to\np_0}(\lambda-1)S_2 \\
 &=\frac{2\Delta_{00}^2\sigma_{111}\sigma_{120}\sigma_{221}\sigma_{210}}{L_1(\sigma_{120})L_2(\sigma_{210})}M_1'(0)M_2'(0)\Big|_{\np=\np_0}
\end{align*}
and, on the other hand, $(c)$ in \teoc{L8} shows that $(\lambda-1)^2\Delta_{11}(\np)$ extends $\cc^{\infty}$ at $\np=\np_0$. This proves the last assertion in $(c)$ and concludes the proof of the result.
\end{prova}

We omit the proof of the next result for the sake of brevity since it is very similar to the previous one. With regard to its statement we recall that $D_{0,n_2}^n=\frac{\N_{\geq n_1}}{n_2}$ and $ D_{n_1,0}^n=
\bigcup_{i=1}^{n_1}\frac{i}{\N_{\geq n_2}}$ (see \obsc{domains}). 

\begin{prop}\label{poles2}
The following assertions hold:
\begin{enumerate}[$(a)$]

\item For any $\np_0=(\lambda_0,\mu_0)\in  D_{0n_2}^n\times W$ with $n_2>0$, the function $\np\mapsto (\lambda-\lambda_0)T_{0n_2}(\np)$ extends $\cc^\infty$ at $\np=\np_0$, and if $\lambda_0=\frac{n_1+i}{n_2}$ with 
$i\in\Z_{\ge 0}$ then $\lim\limits_{\np\to\np_0} (\lambda-\lambda_0)T_{0n_2}(\np)=-\frac{\Delta_{00}^{n_2}}{n_2}\frac{\sigma_{210}^{n_1+i}\sigma_{221}^{n_2}}{L_2^{n_2}(\sigma_{210})}\frac{A_2^{(i)}(0)}{i!}\big|_{\np=\np_0}$.

\item For any $\np_0=(\lambda_0,\mu_0)\in  D_{n_10}^n\times W$ with $\lambda_0\notin\frac{n_1}{\N_{\ge n_2}}$, the function $T_{n_10}(\np)$ extends $\cc^\infty$ at $\np=\np_0.$ In the case that $\lambda_0=\frac{n_1}{n_2+i}$ with $i\in\Z_{\ge 0}$, then the function $\np\mapsto (\lambda-\lambda_0)T_{n_10}(\np)$ extends $\cc^\infty$ at $\np=\np_0$ and $\lim\limits_{\np\to\np_0} (\lambda-\lambda_0)T_{n_10}(\np)=-\frac{n_1}{(n_2+i)^2}\frac{\sigma_{111}^{n_1}\sigma_{120}^{n_2+i}}{L_1^{n_1}(\sigma_{210})}\frac{A_1^{(i)}(0)}{i!}\big|_{\np=\np_0}.$

\end{enumerate}
\end{prop}

Let us recall in regard to the next statement that $D_{0,n_2+1}^n=\frac{\N_{\ge n_1}}{n_2+1}\cup\N,$ see \obsc{domains}.

\begin{prop}\label{poles3}
The following assertions hold:
\begin{enumerate}[$(a)$]

\item For any $\np_0=(\lambda_0,\mu_0)\in  D_{0,n_2+1}^n\times W$ with $\lambda_0\in\N_{\geq\frac{n_1}{n_2}},$ the function $\np\mapsto (\lambda-\lambda_0)^2T_{0,n_2+1}(\np)$ extends~$\cc^\infty$ at $\np=\np_0$, and if $\lambda_0=i\in\N_{\geq\frac{n_1}{n_2}}$ then 
\[
\lim_{\np\to\np_0}(\lambda-\lambda_0)^2T_{0,n_2+1}(\np)=\frac{n_2\Delta_{00}^{n_2+1}\sigma_{210}^{(n_2+1)i}\sigma_{221}^{n_2+1}}{(n_2+1)L_2^{n_2+1}(\sigma_{210})}\frac{M_2^{(i)}(0)}{i!}\frac{A_2^{(n_2i-n_1)}(0)}{(n_2i-n_1)!}\Big|_{\np=\np_0}.
\]

\item For any $\np_0=(\lambda_0,\mu_0)\in  D_{0,n_2+1}^n\times W$ with $\lambda_0\notin\N_{\geq\frac{n_1}{n_2}},$ the function $\np\mapsto (\lambda-\lambda_0)T_{0,n_2+1}(\np)$ extends~$\cc^\infty$ at $\np=\np_0$, and

\begin{enumerate}[$(b1)$]

\item if $\lambda_0=i\in\N_{<\frac{n_1}{n_2}}$ then, setting $i_1\!:=(n_2+1)i-n_1,$ 
\begin{align*}
\lim_{\np\to\np_0}(\lambda-\lambda_0)T_{0,n_2+1}(\np)&=-\frac{\Delta_{00}^{n_2+1}\sigma_{221}^{n_2+1}\sigma_{210}^{n_1}}{L_2^{n_2+1}(\sigma_{210})}
\Bigg(n_2\frac{M_2^{(i)}(0)}{i!}\sigma_{210}^i\hat A_2(in_2-n_1,\sigma_{210})\\
&\qquad\qquad+\frac{n_2\sigma_{210}^{i_1}}{(n_2+1)i_0!}\sum_{j=0}^{i_1}{i_1 \choose j}\frac{M_2^{(j)}(0)A_2^{(i_1-j)}(0)}{j-i}+R
\Bigg)\Bigg|_{\np=\np_0},
\end{align*}
where $R=\frac{\sigma_{210}^{i_1}}{(n_2+1)i_1!}\partial^{i_1}_u\left(L_2^{n_2+1}(u)\partial_2P_1^{-1}(u,0)\right)\big|_{u=0}$ for $i_1\geqslant 0$ and $R=0$ otherwise,

\item if $\lambda_0=\frac{n_1+i}{n_2+1}\notin\N$ with $i\in\Z_{\ge 0}$, then $\lim\limits_{\np\to\np_0}(\lambda-\lambda_0)T_{0,n_2+1}(\np)=-\frac{\Delta_{00}^{n_2+1}\sigma_{221}^{n_2+1}\sigma_{210}^{n_1+i}}{(n_2+1)L_2^{n_2+1}(\sigma_{210})}\frac{B_2^{(i)}(0)}{i!}\big|_{\np=\np_0}.$ 

\end{enumerate}

\end{enumerate}
\end{prop}

\begin{prova}
For the sake of convenience we write $T_{0,n_2+1}$, see $(c)$ in \teoc{A}, as
\begin{equation}\label{44eq1}
T_{0,n_2+1}=f_0\left(f_1+f_2\hat B_2((n_2+1)\lambda-n_1,\sigma_{210})\right)
\end{equation}
with $f_0\!:=\Delta_{00}^{n_2+1}\sigma_{210}^{n_1}\sigma_{221}^{n_2}$, $f_1\!:=\frac{\sigma_{211}}{\sigma_{210}P_1(\sigma_{210},0)}$, $f_2\!:=\frac{\sigma_{221}}{L_2^{n_2+1}(\sigma_{221})}$ and where, recall \refc{def_fun}, 
\[
 B_2(u)=n_2A_2(u)\hat M_2(\lambda,u)+f_3(u) \text{ with $f_3(u)\!:=L_2^{n_2+1}(u)\partial_2 P_1^{-1}(u,0)$.}
\]  
That being said we begin with the proof of $(b2)$. With this aim we note first that $B_2(u;\lambda,\mu)$ is smooth along $\lambda=\lambda_0\notin\Z_{\geq 0}$ because so is $\hat M_2(\lambda,u;\np)$ by $(a)$ in \teoc{L8}. For this reason, since $\frac{n_1+i}{n_2+1}\notin\Z_{\geq 0}$ by assumption, we can apply~\coryc{B22} taking $\alpha=\lambda,$ $\nu=(\lambda,\mu)$, $\alpha_0=\frac{n_1+i}{n_2+1}$, $\nu_0=(\frac{n_1+i}{n_2+1},\mu_0)$, $\kappa_1=n_2+1$ and $\kappa_2=-n_1$ to conclude that 
\[
 \lim_{\np\to\np_0}\left(\frac{n_1+i}{n_2+1}-\lambda\right)\hat B_2\big((n_2+1)\lambda-n_1,\sigma_{210}\big)=\frac{B_2^{(i)}(0)}{(n_2+1)i!}\sigma_{210}^i\Big|_{\np=\np_0}
\]
Hence, on account of \refc{44eq1} and by applying \coryc{B22}, the function $\np\mapsto\big(\lambda-\frac{n_1+i}{n_2+1}\big)T_{0,n_2+1}(\np)$ extends~$\cc^\infty$ at $\np=\np_0$ and tends to $-\frac{\Delta_{00}^{n_2+1}\sigma_{221}^{n_2+1}\sigma_{210}^{n_1+i}}{(n_2+1)L_2^{n_2+1}(\sigma_{210})}\frac{B_2^{(i)}(0)}{i!}\big|_{\np=\np_0}$ as $\np\to\np_0$ and this shows $(b2)$. 

Let us turn now to the proof of assertion $(a)$. So assume that $\lambda_0=i\in\N$ with $n_2i-n_1\geqslant 0$ and observe that, by \coryc{B22}, the function $\np\mapsto (\lambda-i)^2\hat f_3((n_2+1)\lambda-n_1,\sigma_{210})$ extends $\cc^\infty$ at $\np=\np_0$ and tends to 0 as $\np\to\np_0.$ Thus, by applying firstly~$(a)$ in \coryc{B21} and secondly $(a)$ in \lemc{gorrobis} with $\{\alpha=\lambda, \nu=(\lambda,\mu), p=n_1,q=n_2\}$, from \refc{44eq1} we can assert that $\np\mapsto (\lambda-i)^2T_{0,n_2+1}(\np)$ extends $\cc^\infty$ at $\np=\np_0$ and, moreover,
\[
\lim_{\np\to\np_0}(\lambda-i)^2T_{0,n_2+1}(\np)=n_2f_0f_2\big|_{\np=\np_0}
\frac{\sigma_{210}^{n_2i-n_1}}{n_2+1}\frac{M_2^{(i)}(0)}{i!}\frac{A_2^{(n_2i-n_1)}(0)}{(n_2i-n_1)!}\Big|_{\np=\np_0},
\]
which proves $(a)$. In order to show $(b1)$ we consider $\lambda_0=i\in\N$ with $n_2i-n_1<0.$ In this case, if 
$i_1\!:=(n_2+1)i-n_1\geqslant 0$ then $\lim_{\np\to\np_0}(\lambda-i)\hat f_3((n_2+1)\lambda-n_1,\sigma_{210})=\frac{-\sigma_{210}^{i_1}}{n_2+1}\frac{f_3^{(i_1)}(0)}{i_1!}\big|_{\np=\np_0}$ by \coryc{B22}, whereas if $i_1<0$ then $\lim_{\np\to\np_0}(\lambda-i)\hat f_3((n_2+1)\lambda-n_1,\sigma_{210})=0$ by $(a)$ in \teoc{L8}.
Taking this into account the assertion in $(b1)$ follows by applying firstly~$(a)$ in \coryc{B21} and secondly $(b)$ in \lemc{gorrobis} with $\{\alpha=\lambda, \nu=(\lambda,\mu), p=n_1,q=n_2\}$. This concludes the proof of the result.
\end{prova}

Regarding the next statement let us recall, see \obsc{domains}, that $D_{n_1+1,0}^n=\bigcup_{i=1}^{n_1+1}\frac{i}{\N_{\ge n_2}}$.

\begin{prop}\label{poles4}
Let us consider any $\np_0=(\lambda_0,\mu_0)\in  D_{n_1+1,0}^n\times W$. Then the following assertions hold:

\begin{enumerate}[$(a)$]

\item Case $\lambda_0\in\frac{1}{\N}.$ 

\begin{enumerate}

\item[$(a1)$]  If $\lambda_0=\frac{1}{i}$ with $i\in\N_{\ge\frac{n_2}{n_1}}$ then the function $\np\mapsto (\lambda-\lambda_0)^2T_{n_1+1,0}(\np)$ extends~$\cc^\infty$ at $\np=\np_0$ and
 \[
 \lim\limits_{\np\to\np_0} (\lambda-\lambda_0)^2T_{n_1+1,0}(\np)=-\frac{\sigma_{111}^{n_1+1}\sigma_{120}^{(n_1+1)i}}{(n_1+1)i^2L_1^{n_1+1}(\sigma_{120}) }\frac{M_1^{(i)}(0)}{i!}\frac{A_1^{(n_1i-n_2)}(0)}{(n_1i-n_2)!}\Bigg|_{\np=\np_0}.
\]

\item[$(a2)$] If $\lambda_0=\frac{1}{i}$ with $i\in\N\cap[\frac{n_2}{n_1+1},\frac{n_2}{n_1})$ then the function $\np\mapsto (\lambda-\lambda_0)T_{n_1+1,0}(\np)$ extends~$\cc^\infty$ at $\np=\np_0$ and, setting $i_0=(n_1+1)i-n_2,$ 
\begin{align*}
\lim_{\np\to\np_0}(\lambda-\lambda_0)T_{n_1+1,0}(\np)=
-&\frac{\sigma_{111}^{n_1+1}\sigma_{120}^{i(n_1+1)}}{(n_1+1)i^2i_0!L_1^{n_1+1}(\sigma_{120})}
\Bigg(n_1\sum_{j=0}^{i_0}{i_0\choose j}\frac{M_1^{(j)}(0)A_1^{(i_0-j)}(0)}{j-i}\\
&\qquad\qquad+\partial_u^{i_0}\left(L_1^{n_1+1}(u)\partial_1P_2^{-1}(u,0)\right)\big|_{u=0}\Bigg)\Bigg|_{\np=\np_0}.
\end{align*}

\item[$(a3)$] If $\lambda_0=\frac{1}{i}$ with $i\in\N_{<\frac{n_2}{n_1+1}}$ then $T_{n_1+1,0}(\np)$ extends $\cc^\infty$ to $\{\lambda_0\}\times W$.

\end{enumerate}

\item Case $\lambda_0\in\left(\frac{n_1}{\N_{\ge n_2}}\cup\frac{n_1+1}{\N_{\ge n_2}}\right)\setminus\frac{1}{\N}$.

\begin{enumerate}[$(b1)$]

\item If $\lambda_0=\frac{n_1}{n_2+i}\notin\frac{n_1+1}{\N_{\ge n_2}}$ with $i\in\Z_{\ge0}$ then the function $\np\mapsto (\lambda-\lambda_0)T_{n_1+1,0}(\np)$ extends~$\cc^\infty$ at $\np=\np_0$ and 
$
\lim\limits_{\np\to\np_0}(\lambda-\lambda_0)T_{n_1+1,0}(\np)=-\frac{n_1\lambda_0\sigma_{111}^{n_1}\sigma_{120}^{n_2+i}}{(n_2+i)L_1^{n_1}(\sigma_{120})}\frac{A_1^{(i)}(0)}{i!}S_1\Big|_{\np=\np_0}.
$

\item If $\lambda_0=\frac{n_1+1}{n_2+i}\notin\frac{n_1}{\N_{\ge n_2}}$ with $i\in\Z_{\ge0}$ then  $\np\mapsto (\lambda-\lambda_0)T_{n_1+1,0}(\np)$ extends~$\cc^\infty$ at $\np=\np_0$ and
\[
\lim_{\np\to\np_0}(\lambda-\lambda_0)T_{n_1+1,0}(\np)=\frac{n_1\lambda_0\sigma_{111}^{n_1+1}\sigma_{120}^{n_2+i}}{(n_2+i)L_1^{n_1+1}(\sigma_{120})}\frac{(A_1\hat M_1(\frac{1}{\lambda_0},\cdot))^{(i)}(0)}{i!}\Bigg|_{\np=\np_0}.
\]

\item If $\lambda_0=\frac{n_1}{n_2+i_1}=\frac{n_1+1}{n_2+i_2}$ for some $i_1,i_2\in\Z_{\geq 0}$ then the function $\np\mapsto (\lambda-\lambda_0)T_{n_1+1,0}(\np)$ extends~$\cc^\infty$ at $\np=\np_0$ and 
\[
\lim_{\np\to\np_0}(\lambda-\lambda_0)T_{n_1+1,0}(\np)=\frac{n_1\lambda_0\sigma_{111}^{n_1}}{L_1^{n_1}(\sigma_{120})}\left(-\frac{\sigma_{120}^{n_2+i_1}}{n_2+i_1}\frac{A_1^{(i_1)}(0)}{i_1!}S_1+\frac{\sigma_{120}^{n_2+i_2}}{n_2+i_2}\frac{(A_1\hat M_1(\frac{1}{\lambda_0},\cdot))^{(i_2)}(0)}{i_2!}
\right)\Bigg|_{\np=\np_0}.
\]

\end{enumerate}

\item Finally, if $\lambda_0\notin\frac{1}{\N}\cup\frac{n_1}{\N_{\ge n_2}}\cup\frac{n_1+1}{\N_{\ge n_2}}$ then $T_{n_1+1,0}(\np)$ extends $\cc^\infty$ at $\np=\np_0$.

\end{enumerate}
\end{prop}

For the sake of brevity we omit the proof of \propc{poles4}. Let us only mention for reader's convenience that, by $(c)$ in \teoc{A}, 
\[
T_{n_1+1,0}=f_0\left(f_1+f_2\hat B_1((n_1+1)/\lambda-n_2,\sigma_{120})+f_3S_1\hat A_1(n_1/\lambda-n_2,\sigma_{120})\right)
\]
with $f_i\in\cc^\infty(\hat W).$ This expression is similar to the one in \refc{44eq1} for $T_{0,n_2+1}$ that we analysed in the proof of \propc{poles3}, but with the additional summand $f_3S_1\hat A_1$. This extra term increases the number of cases to be studied in terms of $\lambda_0$ but they follow using exactly the same arguments as those explained in the proofs of Propositions~\ref{poles1} and~\ref{poles3}. 

Lastly we state a result concerning the poles of the coefficients $T_{20}$ and $T_{02}$ in the cases $n_1=0$ and $n_2=0$, respectively. For the sake of shortness we do not specify the value of the residues, which can be computed using the same techniques as in the previous results. For the same reason we neither include the proof. With regard to its statement let us recall that $D_{20}^n=\frac{2}{\N_{\ge n_2}}$ and $D_{02}^n=\frac{\N}{2}$, see \obsc{domains}.

\begin{prop}\label{poles5}
The following assertions hold:
\begin{enumerate}[$(a)$]
\item Assume that $n_1=0$ and consider any $\np_0=(\lambda_0,\mu_0)\in D_{20}^n\times W$.
\begin{enumerate}

\item[$(a1)$] If $\lambda_0\in\frac{1}{\N_{\ge n_2}}$ then the function $\np\mapsto(\lambda-\lambda_0)^2T_{20}(\np)$ extends $\cc^\infty$ at $\np_0$.

\item[$(a2)$] If $\lambda_0\notin\frac{1}{\N_{\ge n_2}}$ then the function $\np\mapsto(\lambda-\lambda_0)T_{20}(\np)$ extends $\cc^\infty$ at $\np_0$.

\end{enumerate}
\item Assume that $n_2=0$ and consider any $\np_0=(\lambda_0,\mu_0)\in D_{02}^n\times W$.
\begin{enumerate}
\item[$(b1)$] If $\lambda_0\in\N_{\ge n_1}$ then  the function $\np\mapsto(\lambda-\lambda_0)^2T_{02}(\np)$ extends $\cc^\infty$ at $\np_0$.

\item[$(b2)$] If $\lambda_0\in\N_{<n_1}\cup\left(\frac{\N_{\geq n_1}}{2}\setminus\N\right)$ then the function  $\np\mapsto(\lambda-\lambda_0)T_{02}(\np)$ extends $\cc^\infty$ at~$\np_0$.

\item[$(b3)$] If $\lambda_0\in \frac{\N_{< n_1}}{2}\setminus\N$ then $T_{0n_2}(\np)$ extends $\cc^\infty$ at $\np_0$.

\end{enumerate}
\end{enumerate}
\end{prop}

We are now in position to conclude the proof of \coryc{analitico}. 

\begin{prooftext}{Proof of \coryc{analitico}.}
In the analytic setting (see \obsc{analytic_setting}) we know by \propc{pre-analitico} that the coefficients $\Delta_{ij}$
and $T_{ij}$ listed in \teoc{A} are analytic on $((0,+\infty)\setminus D_{ij}^0)\times W$ and $((0,+\infty)\setminus D_{ij}^0)\times W$, respectively. The fact that each $\Delta_{ij}$ is meromorphic on $\hat W=(0,+\infty)\times W$ with poles of order at most two along $D_{ij}^0\times W$ follows by realising that in the analytic setting the statement of \propc{poles1} is true replacing $\cc^\infty$ by $\cc^\omega,$ i.e., that the extensions are analytic. Indeed, the proof of this analytic version is literally the same but appealing to the analytic assertions in \teoc{L8} instead of the smooth counterparts. More specifically, using $(d)$ in the place of $(a)$ and $(c).$ Similarly, the fact that each $T_{ij}$ is meromorphic on $\hat W=(0,+\infty)\times W$ with poles of order at most two along $D_{ij}^n\times W$ follows by noting that in the analytic setting the statements of Propositions~\ref{poles2}, \ref{poles3}, \ref{poles4} and~\ref{poles5} are true replacing $\cc^\infty$ by $\cc^\omega.$ In this case, besides appealing to $(d)$ in \teoc{L8} in the place of $(a)$ and $(c)$, we apply the analytic versions of \coryc{B22} and \lemc{gorrobis}, i.e., taking $\varpi=\omega$ instead of $\varpi=\infty$. This completes the proof of the result. 
\end{prooftext}

\section{First monomials in the asymptotic expansions}\label{quarta}

\teoc{A} is the main result of the present paper and it is intended to be applied in combination with \teoc{oldA} (which in fact gathers our main results in \cite{MV20}). Because of this, in order to ease the applicability, we next particularise \teoc{oldA} to specify the first monomials appearing in the asymptotic expansion of the Dulac map, see \teoc{3punts}, and the Dulac time, see \teoc{9punts}, for arbitrary hyperbolicity ratio~$\lambda_0$. In both statements,  the order $L$ ranges in a certain interval depending on $\lambda_0.$ The left endpoint of this interval is only given for completeness to guarantee that none of the monomials in the principal part can be included in the remainder.

\begin{theo}\label{3punts}
Let $D(s;\np)$ be the Dulac map of the hyperbolic saddle \refc{X} from $\Sigma_1$ and $\Sigma_2$. 
\begin{enumerate}[$(1)$]

\item If $\lambda_0<1$ then $D(s;\np)=\Delta_{00}(\np) s^\lambda+\Delta_{01}(\np)s^{2\lambda}+\F_L^\infty(\{\lambda_0\}\times W)$ for any $L\in \big[2\lambda_0,\min(3\lambda_0,1+\lambda_0)\big)$.

\item If $\lambda_0=1$ then $D(s;\np)=\Delta_{00}(\np) s^\lambda+\boldsymbol{\Delta}_{10}^{\lambda_0}(\omega;\np)s^{1+\lambda}+\F_L^\infty(\{\lambda_0\}\times W)$ for any $L\in [2,3),$ where 
\[
 \boldsymbol{\Delta}_{10}^{\lambda_0}(\omega;\np)=\Delta_{10}(\np)+\Delta_{01}(\np)(1+\alpha\omega),
\] 
$\alpha=1-\lambda$ and $\omega=\omega(s;\alpha)$.

\item If $\lambda_0>1$ then $D(s;\np)=\Delta_{00}(\np) s^\lambda+\Delta_{10}(\np)s^{\lambda+1}+\F_L^\infty(\{\lambda_0\}\times W)$ for any $L\in \big[\lambda_0+1,\min(2+\lambda_0,2\lambda_0)\big)$.

\end{enumerate}
\end{theo}

\begin{prova}
\begin{enumerate}[(1)]

\item We begin by showing that the assumptions on $\lambda_0$ and $L$ imply $\mathscr B^0_{\lambda_0,L-\lambda_0}=\{(0,0),(0,1)\}$. Let us prove first that $L<\min(3\lambda_0,1+\lambda_0)$ implies $\mathscr B^0_{\lambda_0,L-\lambda_0}\subset\{(0,0),(0,1)\}$. Indeed, we claim that if $(i,j)\in\Lambda_0\setminus\{(0,0),(0,1)\}$ then $(i,j)\notin \mathscr B^0_{\lambda_0,L-\lambda_0}$, i.e., $i+\lambda_0j> L-\lambda_0.$ It is clear that the claim will follow once we prove its validity for $(i,j)=(0,2)$ and $(i,j)=(1,0).$ For the first case observe that $2\lambda_0>L-\lambda_0$ holds because $L<3\lambda_0$ and, for the second one, $1>L-\lambda_0$ holds due to $L<1+\lambda_0$. One can verify similarly that the reverse inclusion $\mathscr B^0_{\lambda_0,L-\lambda_0}\supset\{(0,0),(0,1)\}$ is guaranteed by $2\lambda_0\leqslant L.$ 

Let us show next that $\lambda_0<1$ implies $\lambda_0\notin D^0_{L-\lambda_0}$. To prove this we use firstly that $D_{00}^0\cup D_{01}^0=\N$ by \obsc{domains}, so that $\lambda_0\notin D_{00}^0\cup D_{01}^0.$ Secondly, see \defic{alldefi}, we use that $\lambda_0\in D^0_{L-\lambda_0}$ if and only if there exists $(i,j)\in\mathscr B^0_{\lambda_0,L-\lambda_0}$ such that $\lambda_0\in D_{ij}^0$, which is not possible since $\mathscr B^0_{\lambda_0,L-\lambda_0}=\{(0,0),(0,1)\}$ and $\lambda_0\notin D_{00}^0\cup D_{01}^0.$ Hence $\lambda_0\notin D^0_{L-\lambda_0}$ and the asymptotic expansion follows by $(a1)$ in \teoc{oldA}.

\item Exactly as we did in the previous case, $\lambda_0=1$ and $L\in [2,3)$ yields $\mathscr B_{\lambda_0,L-\lambda_0}^0=\{(0,0),(1,0),(0,1)\}$. This implies, due to $\lambda_0=1\in D_{10}^0=\N$ by \obsc{domains}, that $\lambda_0\in D_{L-\lambda_0}^0$. Then, by $(a2)$ in \teoc{oldA},
\[
 D(s;\np)=\Delta_{00}(\np) s^\lambda+\boldsymbol{\Delta}_{10}^{\lambda_0}(\omega;\np)s^{1+\lambda}+\F_L^\infty(\{\lambda_0\}\times W)
\] 
with $\omega=\omega(s;\alpha)$, $\alpha=1-\lambda$ and $\boldsymbol{\Delta}_{10}^{\lambda_0}(\omega;\np)=\sum_{r=0}^1\Delta_{1-rp,0+rq}(\np)(1+\alpha\omega)^r=\Delta_{10}(\np)+\Delta_{01}(\np)(1+\alpha\omega)$ because, see \defic{alldefi}, $\mathscr A_{01\lambda_0}^0=\{0,1\}$, $\mathscr A_{10\lambda_0}^0=\emptyset$ and $\mathscr A_{00\lambda_0}^0=\{0\}$.

\item Similarly as we argue in $(1)$, in this case the assumptions on $\lambda_0$ and $L$ imply $\mathscr B^0_{\lambda_0,L-\lambda_0}=\{(0,0),(1,0)\}$. Then, since $D_{00}^0\cup D_{10}^0=\frac{1}{\N}$ and $\lambda_0>1$, it turns out that $\lambda_0\notin D_{L-\lambda_0}^0$ and thus the asymptotic expansion in the statement follows by $(a1)$ of \teoc{oldA}.
\end{enumerate}

\noindent This proves the validity of the result. 
\end{prova}

\begin{ex}\label{ex1}
By \teoc{3punts}, if $\lambda_0=1$ then $D(s;\np)=\Delta_{00}(\np) s^\lambda+\boldsymbol{\Delta}_{10}^{\lambda_0}(\omega;\np)s^{1+\lambda}+\F_L^\infty(\{\lambda_0\}\times W)$ for any $L\in [2,3),$ where 
\[
 \boldsymbol{\Delta}_{10}^{\lambda_0}(\omega;\np)=\Delta_{10}(\np)+\Delta_{01}(\np)(1+\alpha\omega),
\] 
$\alpha=1-\lambda$ and $\omega=\omega(s;\alpha)$. The order of monomials in the principal part as $s\to 0^+$ is $s^\lambda\prec_{\lambda_0}s^{1+\lambda}\omega\prec_{\lambda_0}s^{1+\lambda}$, see \cite[Definition 1.7]{MV20} for details. The coefficient of $s^\lambda$ at $\np_0=(1,\mu_0)$ follows directly by evaluating the expression of $\Delta_{00}$ given in assertion $(b)$ of \teoc{A}. The subsequent coefficient is the one of $s^{1+\lambda}\omega$ and, by applying $(b)$ in \propc{poles1} with $i=1,$ its expression at $\np_0=(1,\mu_0)$ is equal to
\[
\lim\limits_{\np\to\np_0}(1-\lambda)\Delta_{01}(\np)=\frac{\Delta_{00}^2\sigma_{221}\sigma_{210}}{L_2(\sigma_{210})}M_2'(0)\big|_{\np=\np_0}.
\]
Moreover some easy computations on account of the definitions given in \refc{def_fun} show that
 \[
  M_2'(0)%=L_2'(0)\partial_2\!\left(\frac{P_2}{P_1}\right)\!(0,0)+L_2(0)\partial_{12}\!\left(\frac{P_2}{P_1}\right)\!(0,0)
  =\partial_1\!\left(\frac{P_2}{P_1}\right)\!(0,0)\partial_2\!\left(\frac{P_2}{P_1}\right)\!(0,0)+\partial_{12}\!\left(\frac{P_2}{P_1}\right)\!(0,0).
 \]
Let us also remark that, more generally, one can compute all the derivatives of $L_i(u),$ $M_i(u),$ $A_i(u),$ $B_i(u)$ and $C_i(u)$ at $u=0$, for $i=1,2,$ in terms of the derivatives of $P_1(x,y)$ and $P_2(x,y)$ at $(x,y)=(0,0)$.
\end{ex}

The second part of \teoc{oldA} provides the asymptotic expansion of the Dulac time associated to a vector field \refc{X} having poles of arbitrary order $n=(n_1,n_2)\in\Z_{\ge 0}^2.$ In \teoc{9punts} we restrict ourselves to the case $n_1=0$ and $n_2\geqslant 1$ for several reasons. Firstly, for the sake of simplicity in the exposition, since dealing with the general situation will increase very much the number of cases to consider. Secondly because the study of the Dulac time of a hyperbolic saddle at infinity of any polynomial vector field of degree~$d$ yields to the case $n_1=0$ and $n_2=d-1.$ Thirdly, and more important for us, because it allows to tackle the conjectural bifurcation diagram of the period function of the quadratic centers that we undertook in~\cite{MMV2}. 

\begin{theo}\label{9punts}
Assuming $n_1=0$ and $n_2\geqslant 1$, let $T (s;\np)$ be the Dulac time of the hyperbolic saddle \refc{X} from $\Sigma_1$ and $\Sigma_2$. 
\begin{enumerate}[$(1)$]
\item If $\lambda_0\in(0,\frac{1}{n_2+1})$ then
$T(s;\np)=T_{00}(\np)+T_{0n_2}(\np)s^{\lambda n_2}+T_{0,n_2+1}(\np)s^{\lambda (n_2+1)}+\F_{L}^\infty(\{\lambda_0\}\times W)$ for any $L\in \big[\lambda_0(n_2+1),\min(1,\lambda_0(n_2+2))\big)$.

\item If $\lambda_0\in(\frac{1}{n_2+1},\frac{2}{n_2+1})\setminus\{\frac{1}{n_2}\}$ then
$$T(s;\np)=T_{00}(\np)+T_{0n_2}(\np)s^{\lambda n_2}+T_{10}(\np)s+T_{0,n_2+1}(\np)s^{\lambda (n_2+1)}+\F_{L}^\infty(\{\lambda_0\}\times W)$$ for any $L\in \big[\max(1,\lambda_0(n_2+1),\min(2,\lambda_0n_2+1,\lambda_0(n_2+2))\big)$.

\item If $\lambda_0\in(\frac{2}{n_2+1},\frac{2}{n_2})$ then $T(s;\np)=T_{00}(\np)+T_{10}(\np)s+T_{0n_2}(\np)s^{\lambda n_2}+T_{20}(\np)s^2+\F_{L}^\infty(\{\lambda_0\}\times W)$ for any $L\in\big[\max(2,\lambda_0n_2),\lambda_0n_2+\min(1,\lambda_0)\big)$.

\item If $\lambda_0>\frac{2}{n_2}$ then $T(s;\np)=T_{00}(\np)+T_{10}(\np)s+T_{20}(\np)s^2+\F_{L}^\infty(\{\lambda_0\}\times W)$ for any $L\in \big[2,\min(3,\lambda_0 n_2)\big)$.

\item If $\lambda_0=\frac{1}{n_2+1}$ then $T(s;\np)=T_{00}(\np)+T_{0n_2}(\np)s^{\lambda n_2}+s\T_{10}^{\lambda_0}(\omega;\np)+\F_L^\infty(\{\lambda_0\}\times W)$ for any $L\in [1,\frac{n_2+2}{n_2+1})$, where 
\[
 \T_{10}^{\lambda_0}(\omega;\np)=T_{10}(\np)+T_{0,n_2+1}(\np)(1+\alpha\omega),
\]  
$\alpha=1-\lambda(n_2+1)$ and $\omega=\omega(s;\alpha)$.

\item If $\lambda_0=\frac{1}{n_2}$ with $n_2>1$ then $T(s;\np)=T_{00}(\np)+s\T_{10}^{\lambda_0}(\omega;\np)+T_{0,n_2+1}(\np)s^{\lambda(n_2+1)}+\F_L^\infty(\{\lambda_0\}\times W)$ for any $L\in[\frac{n_2+1}{n_2},\frac{n_2+2}{n_2})$, where 
\[
 \T_{10}^{\lambda_0}(\omega;\np)=T_{10}(\np)+T_{0 n_2}(\np)(1+\alpha\omega),
\]
$\alpha=1-\lambda n_2$ and $\omega=\omega(s;\alpha)$.

\item If $\lambda_0=\frac{2}{n_2+1}$ with $n_2>1$ then $T(s;\np)=T_{00}(\np)+T_{10}(\np)s+T_{0n_2}(\np)s^{\lambda n_2}+s^2\T_{20}^{\lambda_0}(\omega;\np)+\F_L^\infty(\{\lambda_0\}\times W)$ for any $L\in \big[2,\min\big(\frac{2n_2+4}{n_2+1},\frac{3n_2+1}{n_2+1}\big)\big)$, where
\[
 \T_{20}^{\lambda_0}(\omega;\np)=T_{20}(\np)+T_{0,n_2+1}(\np)\left(1+\alpha\omega\right)^d,
 \] 
$d=\gcd(2,n_2+1)$, $\alpha=\frac{2-\lambda(n_2+1)}{d}$ and $\omega=\omega(s;\alpha).$

\item If $\lambda_0=1$ and $n_2=1$ then
$T(s;\np)=T_{00}(\np)+s\T_{10}^{\lambda_0}(\omega;\np)+s^2\T_{20}^{\lambda_0}(\omega;\np)+\F_{L}^\infty(\{\lambda_0\}\times W)$ for any $L\in [2,3)$, where 
\[
\T_{r0}^{\lambda_0}(\omega;\np)=\sum_{i=0}^rT_{r-i,i}(\np)(1+\alpha \omega)^i\text{, for $r=1,2$,}
\]
$\alpha=1-\lambda$ and $\omega=\omega(s;\alpha)$. 

\item If $\lambda_0=\frac{2}{n_2}$ then $T(s;\np)=T_{00}(\np)+T_{10}(\np)s+s^2\T_{20}^{\lambda_0}(\omega;\np)+\F_L^\infty(\{\lambda_0\}\times W)$ for any $L\in \big[2,\min\big(3,2+\frac{2}{n_2}\big)\big)$, where 
\[
\T_{20}^{\lambda_0}(\omega;\np)=T_{20}(\np)+T_{0n_2}(\np)\left(1+\alpha\omega\right)^d,
\]
$d=\gcd(2,n_2)$, $\alpha=\frac{2-\lambda n_2}{d}$ and $\omega=\omega(s;\alpha).$
\end{enumerate}

\end{theo}

\begin{prova}

The asymptotic expansions in $(1),$ $(2)$, $(3)$ and $(4)$ will follow by applying $(b1)$ in \teoc{oldA} once we determine the grids $\mathscr B_{\lambda_0,L}^n$ and show that, under the respective assumptions on $\lambda_0$ and $L$, we have $\lambda_0\notin D^n_L$. Next we particularise the arguments leading to this in each case: 

\begin{enumerate}[$(1)$]

\item In this case the hypothesis $\lambda_0(n_2+1)\leqslant L<\min(1,\lambda_0(n_2+2))$ yield $\mathscr B_{\lambda_0,L}^n=\{(0,0),(0,n_2),(0,n_2+1)\}$. For instance let us show that $L<\min(1,\lambda_0(n_2+2))$ implies $\mathscr B_{\lambda_0,L}^n\subset\{(0,0),(0,n_2),(0,n_2+1)\}$. To prove this it suffices to check that $(1,0)$ and $(0,n_2+2)$ do not belong to $\mathscr B_{\lambda_0,L}^n$, which is indeed a consequence of $L<1$ and $L<\lambda_0(n_2+2)$, respectively. The reverse inclusion $\supset$ follows similarly taking $\lambda_0(n_2+1)\leqslant L$ into account. Since the assumption $\lambda_0\in(0,\frac{1}{n_2+1})$ and \obsc{domains} imply that $\lambda_0\notin D_{00}^n\cup D_{0n_2}^n\cup D_{0,n_2+1}^n=\emptyset\cup\frac{\N}{n_2}\cup\big(\frac{\N}{n_2+1}\cup\N\big)$, we can assert that $\lambda_0\notin D_L^n.$ 

\item In this case it turns out that $\max(1,\lambda_0(n_2+1)\leqslant L<\min(2,\lambda_0n_2+1,\lambda_0(n_2+2))$ implies that the grid is given by $\mathscr B_{\lambda_0,L}^n=\{(0,0),(0,n_2),(1,0),(0,n_2+1)\}.$ For instance, to show the inclusion $\subset$ is enough to verify that $(2,0)$, $(1,n_2)$ and $(0,n_2+2)$ do not belong to $\mathscr B_{\lambda_0,L}^n$, which is a consequence of $L<2$, $L<1+\lambda_0n_2$ and $L<\lambda_0(n_2+2)$, respectively. That being said, we know by \obsc{domains} that $D_{0n_2}^n=\frac{\N}{n_2}$, $D_{10}^n=\frac{1}{\N_{\geq n_2}}$ and $D_{0,n_2+1}^n=\frac{\N}{n_2+1}\cup\N.$ Thus, on account of the assumption  
$\lambda_0\in(\frac{1}{n_2+1},\frac{2}{n_2+1})\setminus\{\frac{1}{n_2}\}$, we get $\lambda_0\notin D_{00}^n\cup D_{0n_2}^n\cup D_{10}^n\cup D_{0,n_2+1}^n$. Hence, see \defic{alldefi}, $\lambda_0\notin D^n_L$.

\item If $\max(2,n_2\lambda_0)\leqslant L<\min(\lambda_0 n_2+1,\lambda_0(n_2+1))$ then $\mathscr B_{\lambda_0,L}^n=\{(0,0),(1,0),(2,0),(0,n_2)\}$. Indeed, the lower bound gives the inclusion $\supset.$ To prove the inclusion $\subset$ it suffices to check that $(3,0)$, $(1,n_2)$ and $(0,n_2+1)$ do not belong to  
$\mathscr B_{\lambda_0,L}^n$, which is a consequence of $L<3$, $L<1+\lambda_0n_2$ and $L<\lambda_0(n_2+1)$, respectively. These three inequalities follow by the assumption $L<\min(\lambda_0 n_2+1,\lambda_0(n_2+1))$ together with the fact that $\lambda_0n_2<2$ due to $\lambda_0\in(\frac{2}{n_2+1},\frac{2}{n_2})$. This last condition, taking \obsc{domains} also into account, implies $\lambda_0\notin D_{00}^n\cup D_{10}^n\cup D_{0,n_2}^n\cup D_{20}^n=\emptyset\cup\frac{1}{\N_{\geq n_2}}\cup\frac{\N}{n_2}$ and then $\lambda_0\notin D^n_L$.

\item Similarly as in the previous cases, if $2\leqslant L<\min(3,\lambda_0 n_2)$ then $\mathscr B^n_{\lambda_0,L}=\{(0,0),(1,0),(2,0)\}$. Moreover, by \obsc{domains} and the hypothesis $\lambda_0>\frac{2}{n_2}$, we get $\lambda_0\notin D^n_{00}\cup D^n_{10}\cup D^n_{20}=\emptyset\cup\frac{1}{\N_{\geq n_2}}\cup\frac{2}{\N_{\geq n_2}}$. Therefore $\lambda_0\notin D^n_L$.

\end{enumerate}
The remaining assertions follow by applying $(b2)$ in \teoc{oldA}. To this end we need to verify that $\lambda_0\in D^n_L$ and determine the grid $\mathscr B_{\lambda_0,L}^n$ together with the corresponding sets $\mathscr A_{ij\lambda_0}^n$. As before we next particularise this in each case: 

\begin{enumerate}[$(1)$]

\setcounter{enumi}{4}

\item If  $\lambda_0=\frac{1}{n_2+1}$ and
$1\leqslant L<1+\frac{1}{n_2+1}$ then 
$\mathscr  B^n_{\lambda_0,L}=\{(0,0),(0,n_2),(0,n_2+1),(1,0)\}$. 
Indeed, to show the inclusion $\subset$ it suffices to check that $(1,n_2)$, $(0,n_2+2)$ and $(2,0)$ do not belong to $\mathscr  B^n_{\lambda_0,L}$, which is equivalent to $L<1+\lambda_0n_2=1+\frac{n_2}{n_2+1}$, $L<\lambda_0(n_2+2)=1+\frac{1}{n_2+1}$ and $L<2$, respectively. These three conditions are a consequence of the assumption $L<1+\frac{1}{n_2+1}$. With regard to the inclusion $\supset$, the fact that $(0,0)$, $(0,n_2),$ $(0,n_2+1)$ and $(1,0)$ belong to $\mathscr B^n_{\lambda_0,L}$ is written as $L\geqslant 0,$ $L\geqslant \lambda_0n_2=\frac{n_2}{n_2+1}$, $L\geqslant \lambda_0(n_2+1)=1$ and $L\geqslant 1,$ respectively, which are guaranteed by the assumption $L\geqslant 1.$ Since, on the other hand, $\lambda_0=\frac{1}{n_2+1}\in D_{0,n_2+1}^n=\frac{\N}{n_2+1}$ by \obsc{domains}, it turns out that $\lambda_0\in D^n_L$, see \defic{alldefi}.

Finally the result follows, see \defic{alldefi} again, using that $\mathscr A_{10\lambda_0}^n=\{0,1\},$ $\mathscr A_{00\lambda_0}^n=\mathscr A_{0n_2\lambda_0}^n=\{0\}$ and $\mathscr A_{0,n_2+1,\lambda_0}^n=\emptyset,$ together with $p=1$ and $q=n_2+1$, so that  $\alpha=1-\lambda(n_2+1).$

\setcounter{enumi}{5}
\item If $\lambda_0=\frac{1}{n_2}$ with $n_2>1$ and $\frac{n_2+1}{n_2}\leqslant L<\frac{n_2+2}{n_2}$  then, just as we argue in the previous cases, we get that $\mathscr B_{\lambda_0,L}^n=\{(0,0),(1,0),(0,n_2),(0,n_2+1)\}$. Furthermore, since $\lambda_0=\frac{1}{n_2}\in D_{10}^n=\frac{1}{\N_{\geq n_2}}$ by \obsc{domains}, it turns out that $\lambda_0\in D^n_L$. On account of this the result follows using that $\mathscr A_{10\lambda_0}^n=\{0,1\},$ $\mathscr A_{00\lambda_0}^n=\mathscr A_{0,n_2+1,\lambda_0}^n=\{0\}$ and $\mathscr A_{0n_2\lambda_0}^n=\emptyset,$ together with the fact that $\alpha=1-\lambda n_2$, which in turn follows due to $p=1$ and $q=n_2$.

\item If $\lambda_0=\frac{2}{n_2+1}$ with $n_2>1$ and $2\leqslant L<\min\big(\frac{2n_2+4}{n_2+1},\frac{3n_2+1}{n_2+1}\big)$ then 
\[
\mathscr B^n_{\lambda_0,L}=\{(0,0),(1,0),(0,n_2),(2,0),(0,n_2+1)\}.
\] 
As usual, the inequality $L<\min\big(\frac{2n_2+4}{n_2+1},\frac{3n_2+1}{n_2+1}\big)$ gives the inclusion $\subset$, in this case by showing that $(3,0),(1,n_2),(0,n_2+2)\notin\mathscr B^n_{\lambda_0,L}$, whereas the inequality $2\leqslant L$ implies the reverse inclusion $\supset$. Hence, since $\lambda_0=\frac{2}{n_2+1}\in D_{02}^n=\frac{2}{\N_{\geq n_2}}$ by \obsc{domains}, we conclude that $\lambda_0\in D_L^n.$ On the other hand, due to $n_2>1$, 
one can verify that $\mathscr A_{00\lambda_0}^n=\mathscr A_{10\lambda_0}^n=\mathscr A_{0n_2\lambda_0}^n=\{0\}$, $\mathscr A_{0,n_2+1,\lambda_0}^n=\emptyset$ and $\mathscr A_{02\lambda_0}^n=\{0,d\},$ where $d=\gcd(2,n_2+1).$ Since $p=\frac{2}{d}$ and $q=\frac{n_2+1}{d}$, the last equality yields
 \[
  \T_{20}^{\lambda_0}(\omega;\np)=\sum_{r\in\{0,d\}}T_{2-\frac{2}{d}r,\frac{n_2+1}{d}r}(\np)\left(1+\alpha\omega\right)^r
  =T_{20}(\np)+T_{0,n_2+1}(\np)\left(1+\alpha\omega\right)^d,
 \]
where $\omega=\omega(s;\alpha)$ and $\alpha=\frac{2-\lambda(n_2+1)}{d}$. This proves the validity of the statement.  

\item If $\lambda_0=1$, $n_2=1$ and $2\leqslant L<3$ then one can readily show that
 \[
 \mathscr B^n_{\lambda_0,L}=\{(0,0),(1,0),(0,1),(2,0),(1,1),(0,2)\}.
 \]
On account of this, since $\lambda_0=1\in D_{01}^n=\N$ by \obsc{domains} due to $n=(0,1)$, we can assert that $\lambda_0\in D_L^n.$ In this case one can easily verify that $\mathscr A_{00\lambda_0}^n=\{0\}$, $\mathscr A_{01\lambda_0}^n=\mathscr A_{02\lambda_0}^n=\mathscr A_{11\lambda_0}^n=\emptyset$, $\mathscr A_{10\lambda_0}^n=\{0,1\}$ and $\mathscr A_{20\lambda_0}^n=\{0,1,2\}.$ Since $p=q=1$, the two last equalities show, respectively,
\[
\T_{r0}^{\lambda_0}(\omega;\np)=\sum_{i=0}^rT_{r-i,i}(\np)(1+\alpha \omega)^i\text{, for $r=1,2$,}
\]
where $\alpha=1-\lambda$ and $\omega=\omega(s;\alpha).$

\item If $\lambda_0=\frac{2}{n_2}$ and $2\leqslant L<\min(3,2+\frac{2}{n_2})$ then $\mathscr B^n_{\lambda_0,L}=\{(0,0),(1,0),(2,0),(0,n_2)\}$. Consequently, due to $\lambda_0=\frac{2}{n_2}\in D^n_{0n_2}=\frac{\N}{n_2}$ by \obsc{domains}, we have $\lambda_0\in D^n_L$. Moreover $\mathscr A_{00\lambda_0}^n=\mathscr A_{10\lambda_0}^n=\{0\}$, $\mathscr A_{0n_2\lambda_0}^n=\emptyset$ and $\mathscr A_{20\lambda_0}^n=\{0,d\}$ with 
$d=\gcd(2,n_2).$ Since $p=\frac{2}{d}$ and $q=\frac{n_2}{d}$, from the last equality it follows that
 \[
  \T_{20}^{\lambda_0}(\omega;\np)=\sum_{r\in\{0,d\}}T_{2-\frac{2}{d}r,\frac{n_2}{d}r}(\np)\left(1+\alpha\omega\right)^r
  =T_{20}(\np)+T_{0,n_2}(\np)\left(1+\alpha\omega\right)^d,
 \]
where $\omega=\omega(s;\alpha)$ and $\alpha=\frac{2-\lambda n_2}{d}$.

\end{enumerate}
\noindent This concludes the proof of the result.
\end{prova}

\begin{figure}[t]
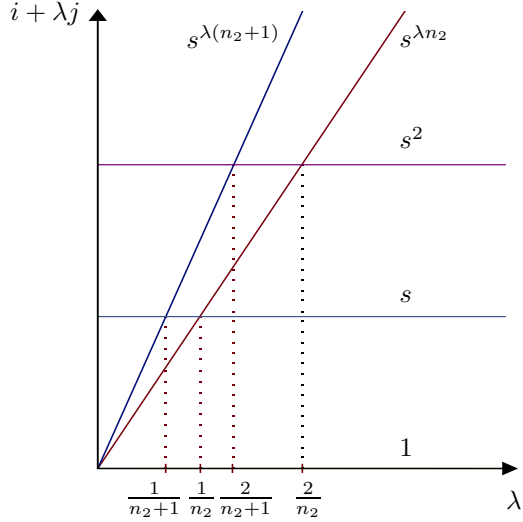

 \centering
 \begin{lpic}[l(0mm),r(0mm),t(0mm),b(0mm)]{diagram(0.4,0.3)}
   \lbl[l]{133,18;$1$}
   \lbl[l]{133,85;$s$}
   \lbl[l]{133,155;$s^2$}
   \lbl[l]{63,200;$s^{\lambda (n_2+1)}$}
   \lbl[l]{133,200;$s^{\lambda n_2}$}
   \lbl[l]{43,-5;$\frac{1}{n_2+1}$}
   \lbl[l]{63,-5;$\frac{1}{n_2}$}
    \lbl[l]{73,-5;$\frac{2}{n_2+1}$}
    \lbl[l]{98,-5;$\frac{2}{n_2}$}
    \lbl[l]{168,-5;$\lambda$}
    \lbl[l]{5,210;$i+\lambda j$}
 \end{lpic}\bigskip     
 \caption{Going upward from each abscissa $\lambda_0\in (0,+\infty)$, order of monomials $s^{i+\lambda j}$ as $s\to 0^+$ and $\lambda\approx\lambda_0$ for $(i,j)\in\{(0,0),(1,0),(2,0),(0,n_2),(0,n_2+1)\}$.}\label{ordre_monomis}
 \end{figure}

Let us finish this section by pointing out that the formula of every coefficient $T_{ij}$ appearing in \teoc{9punts} is given in assertion $(c)$ of \teoc{A}, except for $T_{11}$ in point $(8)$, that corresponds to $\lambda_0=n_2=1$. The formula of this coefficient follows by applying also assertions 
 $(a)$ and $(b)$, which show that $T_{11}=\Omega_{10}T_{01}$ and $\Omega_{10}=\lambda S_1$. Also with regard to this statement, it is worth noting that the order as $s\to 0^+$ of the monomials in points from~$(1)$ to~$(4)$ follow readily from \figc{ordre_monomis}. 
For instance, $1\prec_{\lambda_0}s^{\lambda n_2}\prec_{\lambda_0}s^{\lambda(n_2+1)}\prec_{\lambda_0}s\prec_{\lambda_0}s^{2}$ for $\lambda_0\in (0,\frac{1}{n_2+1})$ and $1\prec_{\lambda_0}s^{\lambda n_2}\prec_{\lambda_0}s\prec_{\lambda_0}s^{\lambda(n_2+1)}\prec_{\lambda_0}s^{2}$ for $\lambda_0\in (\frac{1}{n_2+1},\frac{1}{n_2})$, see \cite[Definition 1.7]{MV20} for details. For $\lambda\approx\lambda_0=\frac{1}{n_2+1}$,
which corresponds to an intersection between two straight-lines in \figc{ordre_monomis}, the compensators come into play and we have $1\prec_{\lambda_0}s^{\lambda n_2}\prec_{\lambda_0}s\omega(s;\alpha)\prec_{\lambda_0}s\prec_{\lambda_0}s^{2}$ with $\alpha=1-\lambda(n_2+1),$ see point $(5)$ in \teoc{9punts}.
This type of information is very relevant in order to apply \cite[Theorem C]{MV20} to bound the number of critical periods or limit cycles that bifurcate from a hyperbolic polycycle. 

\appendix

\section{Derivatives of regular transition map and transition time}\label{ap_regular}
In this section we consider a family of vector fields of the form 
 \begin{equation}\label{ap1}
  Y_{\nu}=\frac{1}{y^{\n}f(x,y;\nu)}\bigl(\partial_x+yh(x,y;\nu)\partial_y\bigr),
 \end{equation}
where 
\begin{itemize}
\item $\n\in\Z$ and $\nu\in U$, where $U$ is some open set of $\R^N,$
\item $f,h\in\mathscr C^K(V\!\times\! U)$ with $V\!:=(a,b)\!\times\!(-c,c)\subset\R^2$, $a<b$ and $c>0,$ 
\item $f(x,0;\nu)\neq 0$ for all $x\in (a,b)$ and $\nu\in U$.
\end{itemize}
We also consider two $\mathscr C^K$ families of transverse sections \map{\xi(\,\cdot\,;\nu)}{(-\varepsilon,\varepsilon)}{\Pi_1} and \map{\zeta(\,\cdot\,;\nu)}{(-\varepsilon,\varepsilon)}{\Pi_2}
to the straight line $\{y=0\}$, i.e., verifying $\xi_2(0)=\zeta_2(0)=0$ together with $\xi_2'(0)\neq 0$ and $\zeta_{2}'(0)\neq 0.$ Our goal is to give the first non-trivial terms of the transition map $P(\,\cdot\,;\nu)$ and the transition time $T(\,\cdot\,;\nu)$ between $\Pi_1$ and $\Pi_2.$ More precisely, denoting by $\varphi(t,p_0;\nu)$ the solution of $Y_{\nu}$ with initial condition $p_0\in V$, we define $P(s;\nu)$ and $T(s;\nu)$ by means of $\varphi(T(s),\xi(s))=\zeta(P(s)).$ The smoothness assumption for the results in this appendix is $K\geqslant 3.$ 

In what follows $\phi(t,p_0;\nu)$ denotes the solution of $Z_{\nu}\!:=\partial_x+yh(x,y;\nu)\partial_y$ with initial condition at $p_0=(x,y).$ It is clear that $\phi(t,p_0;\nu)=\big(x+t,\phi_2(t,p_0;\nu))$. With regard to the second component we prove the next result:

\begin{lem}\label{L1}
Let us define $\displaystyle H(x,y;\nu)=\exp\left(\int_y^xh(u,0;\nu)du\right)$. Then the following hold:
\begin{enumerate}[$(a)$]
\item $\partial_x\phi_2(t,(x,0))=0$ and $\partial_{xx}^2\phi_2(t,(x,0))=0$,
\item $\partial_y\phi_2(t,(x,0))=H(x+t,x)$ and $\partial_{xy}^2\phi_2(t,(x,0))=H(x+t,x)\big(h(x+t,0)-h(x,0)\big)$,
\item $\partial_{yy}^2\phi_2(t,(x,0))=2H(x+t,x)\int_0^tH(x+v,x)\partial_2 h(x+v,0)dv$.
\end{enumerate}
\end{lem}

\begin{prova}
On account of $\partial_t\phi_2(t,(x,y))=\phi_2(t,(x,y))h\big(x+t,\phi_2(t,(x,y))\big)$ and $\phi_2(t,(x,0))=0$ we obtain
\[
\partial_t\partial_x\phi_2(t,(x,0))=h(x+t,0)\partial_x\phi_2(t,(x,0)).
\]
Since $\partial_x\phi_2(0,(x,0))=0$ due to $\phi_2(0,(x,y))=y$, we get $\partial_x\phi_2(t,(x,0))=0$. Accordingly $\partial^2_{xx}\phi_2(t,(x,0))=0$ and this shows $(a)$.
Similarly we obtain $\partial_t\partial_y\phi_2(t,(x,0))=h(x+t,0)\partial_y\phi_2(t,(x,0))$ and $\partial_y\phi_2(0,(x,0))=1$. Consequently
\begin{align}\label{L1eq1}
 &\partial_y\phi_2(t,(x,0))=\exp\left(\int_0^th(x+u,0)du\right)=H(x+t,x)
 \intertext{and}\notag
 &\partial^2_{xy}\phi_2(t,(x,0))=\exp\left(\int_0^th(x+u,0)du\right)\int_0^t\partial_1h(x+u,0)du=H(x+t,x)\big(h(x+t,0)-h(x,0)\big),
\end{align}
which shows the validity of $(b)$. Finally, using that
\begin{align*}
 \partial_t\partial_{yy}^2\phi_2(t,(x,0))&=\left.\partial_{yy}^2\Big(h(x+t,\phi_2(t,(x,y)))\phi_2(t,(x,y))\Big)\right|_{y=0}\\
 &=2\partial_2h(x+t,0)(\partial_y\phi_2(t,(x,0)))^2+h(x+t,0)\partial_{yy}^2\phi_2(t,(x,0)),
 \end{align*}
together with $\partial_{yy}^2\phi_2(0,(x,0))=0$ and \refc{L1eq1}, we get
\[
 \partial_{yy}^2\phi_2(t,(x,0))=2\exp\left(\int_0^th(x+u,0)du\right)\int_0^t\exp\left(\int_0^vh(x+u,0)du\right)\partial_2 h(x+v,0)dv.
\]
Taking \refc{L1eq1} into account once again, the above equality shows $(c)$ and concludes the proof of the result.
\end{prova}

Let us remark that in the previous result (and in what follows when there is no risk of ambiguity) we omit the dependence with respect to the parameter $\nu$ for the sake of shortness. Note on the other hand that the solution $\varphi(t,\xi(s))$ of~$Y_{\nu}$ is inside $\{y=\phi_2(x-\xi_1(s),\xi(s))\}$. Thus, in order to obtain the first coefficients of the Taylor expansion of $T(s)$ and $P(s)$ at $s=0$, we compute first the ones of 
\[
 s\longmapsto\Omega(x,s;\nu)\!:=\phi_2\big(x-\xi_1(s;\nu),\xi(s;\nu);\nu\big).
\]
This is done in the next result, where $H(x,y)=\exp\left(\int_y^xh(u,0)du\right)$, see \lemc{L1}, and we use the compact notation 
$\xi_{ik}=\xi_i^{(k)}(0)$ for $i=1,2.$

\begin{lem}\label{L2}
The function $\Omega(x,s;\nu)$ is $\mathscr C^K$ on $(a,b)\!\times\!(-\varepsilon,\varepsilon)\!\times\! U.$ Moreover it verifies
$\Omega(x,0;\nu)=0$, $\rho_1(x;\nu)\!:=\partial_s\Omega(x,0;\nu)=\xi_{21}H(x,\xi_{10})$ and 
\[
 \rho_2(x;\nu)\!:=\partial_{ss}^2\Omega(x,0;\nu)=H(x,\xi_{10})\left(\xi_{22}-2\xi_{11}\xi_{21}h(\xi_{10},0)+2\xi_{21}^2\int_{\xi_{10}}^xH(u,\xi_{10})\partial_2h(u,0)du\right).
\]
\end{lem}

\begin{prova}
The fact that $\Omega$ is $\mathscr C^K$ on $(a,b)\!\times\!(-\varepsilon,\varepsilon)\!\times\! U$ follows from the smooth dependence of solutions with respect to initial conditions and parameters (see for instance \cite[Theorem 1.1]{Hale}) and that $\Omega(x,0;\nu)=0$ is due to the invariance of the straight line $\{y=0\}.$

Since $\phi(t,(x,y))$ is the solution of $Z_{\nu}$ with initial condition at $(x,y)$, in order to to avoid any ambiguity we 
consider $\Omega(z,s)=\phi_2(z-\xi_1(s),\xi(s))$ and so we keep the notation 
$\partial_t$, $\partial_x$ and $\partial_y$ for the partial derivatives of $\phi_2(t,(x,y)).$ In doing so we obtain
 \begin{align}\label{L2eq1}\notag
  \rho_1(z)=\left.\partial_s\phi_2\big(z-\xi_1(s),\xi(s)\big)\right|_{s=0}=&-\partial_t\phi_2(z-\xi_1(s),\xi(s))\xi_1'(s)\\\notag&+\left.\partial_x\phi_2(z-\xi_1(s),\xi(s))\xi_1'(s)+\partial_y\phi_2(z-\xi_1(s),\xi(s))\xi_2'(s)\right|_{s=0}\\\notag
  =&-\phi_2(z-\xi_1(s),\xi(s))h\big(z,\phi_2(z-\xi_1(s),\xi(s))\big)\xi_1'(s)\\&\left.+\partial_x\phi_2(z-\xi_1(s),\xi(s))\xi_1'(s)+\partial_y\phi_2(z-\xi_1(s),\xi(s))\xi_2'(s)\right|_{s=0}\\\notag
  =&\,\xi_{21}H(z,\xi_{1s}),
 \end{align}
where in the third equality we use that $\phi$ is the flow of $Z_{\nu}=\partial_x+yh(x,y;\nu)\partial_y$ and in the fourth one that $\phi_2(x-\xi_1(0),\xi(0))=0$ due to $\xi_2(0)=0$, together with $\partial_x\phi_2(t,(x,0))=0$ and $\partial_y\phi_2(t,(x,0))=H(x+t,x)$, as established by \lemc{L1}. 

Next we proceed with the computation of $\rho_2(z).$ With this aim in view note that, from \refc{L2eq1},
\begin{align}\notag
\rho_2(z)=\left.\partial_{ss}^2\phi_2(z-\xi_1(s),\xi(s))\right|_{s=0}=&-\xi_{11}h(z,0)\partial_s\phi_2\big(z-\xi_1(s),\xi(s)\big)+\xi_{11}\partial_s\partial_x\phi_2(z-\xi_1(s),\xi(s))\\
&+\left.\partial_s\Big(\partial_y\phi_2(z-\xi_1(s),\xi(s))\xi_2'(s)\Big)\right|_{s=0}.\label{L2eq2}
\end{align}
By applying \lemc{L1}, some computations show that
\begin{align*}
\left.\partial_s\partial_x\phi_2(z-\xi_1(s),\xi(s))\right|_{s=0}&=\xi_{11}\Big(-\partial_t\partial_x\phi_2+\partial_{xx}^2\phi_2\Big)(z-\xi_{10},(\xi_{10},0))+\xi_{21}\partial_{xy}^2\phi_2(z-\xi_{10},(\xi_{10},0))\\
&=\xi_{21}H(z,\xi_{10})\big(h(z,0)-h(\xi_{10},0)\big).
\intertext{and}
\left.\partial_s\partial_y\phi_2(z-\xi_1(s),\xi(s))\right|_{s=0}&=\xi_{11}\Big(-\partial_t\partial_y\phi_2+\partial_{xy}^2\phi_2\Big)(z-\xi_{10},(\xi_{10},0))+\xi_{21}\partial_{yy}^2\phi_2(z-\xi_{10},(\xi_{10},0))\\
&=H(z,\xi_{10})\left(-\xi_{11}h(\xi_{10},0)+2\xi_{21}\int_{\xi_{10}}^zH(u,\xi_{10})\partial_2h(u,0)du\right).
\end{align*}
Since $\left.\partial_s\phi_2\big(z-\xi_1(s),\xi(s)\big)\right|_{s=0}=\partial_y\phi_2\big(z-\xi_1(0),\xi(0)\big)\xi_2'(0)=\xi_{21}H(z,\xi_{10})$ by \lemc{L1} once again, the substitution of the two previous identities in \refc{L2eq2} yields
 \[
 \rho_2(z)=H(z,\xi_{10})\left(\xi_{22}-2\xi_{11}\xi_{21}h(\xi_{10},0)+2\xi_{21}^2\int_{\xi_{10}}^zH(u,\xi_{10})\partial_2h(u,0)du\right),
 \]
as desired. Hence the result is proved. 
\end{prova}

We are now in position to give the two first non-trivial coefficients of the transition map $P(\,\cdot\,;\nu)$ and the transition time $T(\,\cdot\,;\nu)$ between $\Pi_1$ and $\Pi_2.$ In this regard it is to be quoted a previous result by Chicone (see \cite[Theorem 2.2]{Chicone92}), where it is given the expression of $\partial_s P(0;\nu)$ for vector fields in general position, i.e., not assuming that the straight line $\{y=0\}$ is invariant. He also gives the formula of $\partial_s T(0;\nu)$ in the case that $\ell=0.$ More recently, explicit formulas of $\partial_s P(0;\nu)$ and also $\partial_{ss} P(0;\nu)$ for vector fields in general position are given in \cite[Theorem 4.2]{Leen09}. The proofs in \cite{Chicone92,Leen09} are based on Diliberto's theorem on the integration of the homogeneous variational equations of a plane autonomous differential system in terms of geometric quantities along a given trajectory. (Similar results for the transition map can be found in the book of Andronov {\it et al.}~\cite{Andronov}.) In our next lemma, besides these coefficients, we also give the second coefficient of the transition time, which to the best of our knowledge constitutes a new result. The lemma is in fact an upgrade of \cite[Lemma 2.4]{MV20}, where we study the regularity properties of these maps without giving the expression of the coefficients. In the statement for the sake of shortness we use the compact notation 
$\xi_{ik}=\xi_i^{(k)}(0)$ and $\zeta_{ik}=\zeta_i^{(k)}(0)$, $i=1,2$, for the derivatives of the parametrization of the transverse sections. We also remark that the functions~$\rho_1$ and~$\rho_2$ appearing in these coefficients are the ones given in \lemc{L2}.

\begin{lem}\label{L3} 
Let $P(s;\nu)$ and $T(s;\nu)$ be respectively the transition map and transition time of the flow given by~\refc{ap1} between the transverse sections \map{\xi(\,\cdot\,;\nu)}{(-\varepsilon,\varepsilon)}{\Pi_1} and \map{\zeta(\,\cdot\,;\nu)}{(-\varepsilon,\varepsilon)}{\Pi_2} to $\{y=0\}$. Then the following hold:
 \begin{enumerate}[$(a)$]
  \item The function $P(s;\nu)$ is $\mathscr C^K$ on $\big((-\varepsilon,\varepsilon)\times U\big)$. Moreover 
         $P(0;\nu)=0,$ 
         \begin{align*} 
         &p_1(\nu)\!:=\partial_sP(0;\nu)=\frac{\xi_{21}}{\zeta_{21}}\exp\left(\int_{\xi_{10}}^{\zeta_{10}}h(u,0)du\right)
          %H(\zeta_{10},\xi_{10})
          \intertext{and}
          &p_2(\nu)\!:=\partial_{ss}^2P(0;\nu)=\frac{\big(2\zeta_{11}\zeta_{21}h(\zeta_{10},0)-
          \zeta_{22}\big)p_1^2+\rho_2(\zeta_{10})}{\zeta_{21}}.
         \end{align*} 
 \item $T(s;\nu)=s^{\n}\tilde T(s;\nu)$ with $\tilde T\in\mathscr C^{K-1}\big((-\varepsilon,\varepsilon)\times U\big)$  
          verifying $\tilde T(0;\nu)=\displaystyle\int_{\xi_{10}}^{\zeta_{10}}\!\!\rho_1^{\n}(x)f(x,0)dx$ and
          \begin{align*}
           \partial_s\tilde T(0;\nu)=&\,\zeta_{11}\zeta_{21}^{\n}p_1^{{\n}+1}f(\zeta_{10},0)
   -\xi_{11}\xi_{21}^{\n}f(\xi_{10},0) \\
   &+\frac{1}{2}\int_{\xi_{10}}^{\zeta_{10}}\!\!\rho_1^{{\n}-1}(x)\Big({\n}\rho_2(x)f(x,0)+2\rho_1^2(x)\partial_2f(x,0)\Big)dx.
          \end{align*}
  Moreover if ${\n}=0$ then $T\in\mathscr C^{K}\big((-\varepsilon,\varepsilon)\times U\big)$ and
  \begin{align*}
   \partial_{ss}^2T(0;\nu)=&\,\big(\zeta_{12}p_1^2+\zeta_{11}p_2\big)f(\zeta_{10},0)
   +\zeta_{11}^2p_1^2\partial_1f(\zeta_{10},0)+2\zeta_{11}\zeta_{21}p_1^2\partial_2f(\zeta_{10},0)\\[5pt]
   &-\xi_{12}f(\xi_{10},0)-\xi_{11}^2\partial_1f(\xi_{10},0)-2\xi_{11}\xi_{21}\partial_2f(\xi_{10},0)\\
   &+\int_{\xi_{10}}^{\zeta_{10}}\Big(\rho_1^2(x)\partial_{22}^2f(x,0)+\rho_2(x)\partial_2f(x,0)\Big)dx.
  \end{align*}
 \end{enumerate}
\end{lem}
\begin{prova}
The assertion concerning the smoothness of $P(s;\nu)$ follows by the smooth dependence of solutions with respect to initial conditions and parameters and the application of the implicit function theorem (see for instance \cite[Theorem 1.1]{Hale}). Note on the other hand that, by definition, $\varphi(T(s),\xi(s))=\zeta(P(s))$ where $\varphi(t,p_0)$ is solution of $Y_{\nu}$ with initial condition $p_0\in V.$ Since $Z_{\nu}=y^\ell f(x,y;\nu)Y_{\nu}=\partial_x+yh(x,y)\partial_y$, it follows that
 \[
  \zeta_2(P(s))=\phi_2\big(\zeta_1(P(s))-\xi_1(s),\xi(s)\big)=\Omega\big(\zeta_1(P(s)),s\big),
 \]
where $\phi(t,(x,y))=(t+x,\phi_2(t,(x,y))$ is the flow of $Z_{\nu}$ and, by definition, $\Omega(x,s)=\phi_2(x-\xi_1(s),\xi(s)).$ Accordingly
 \[
  \zeta_2'(P(s))P'(s)=\partial_1\Omega\big(\zeta_1(P(s)),s\big)\zeta_1'(P(s))P'(s)+\partial_2\Omega\big(\zeta_1(P(s)),s\big),
 \]
which, evaluated at $s=0$ and applying \lemc{L2}, gives $\zeta_{21}P'(0)=\partial_2\Omega(\zeta_{10},0)=\rho_1(\zeta_{10})=\xi_{21}H(\zeta_{10},\xi_{10}).$ Therefore $p_1=P'(0)=\frac{\xi_{21}}{\zeta_{21}}H(\zeta_{10},\xi_{10}),$ as desired. By computing an additional derivative with respect to $s$ in the above equality and evaluating at $s=0$ afterwards we get
 \[
  \zeta_{22}p_1^2+\zeta_{21}P''(0)=2\partial_{12}^2\Omega(\zeta_{10},0)\zeta_{11}p_1+\partial_{22}^2\Omega(\zeta_{10},0)
  =2\rho_1(\zeta_{10})h(\zeta_{10},0)\zeta_{11}p_1+\rho_2(\zeta_{10}),
 \]
where we apply \lemc{L2} and take $\rho_1'(\zeta_{10})=\xi_{21}\partial_1H(\zeta_{10},\xi_{10})=\xi_{21}H(\zeta_{10},\xi_{10})h(\zeta_{10},0)=\zeta_{21}p_1h(\zeta_{10},0)$ into account. Consequently,
 \[
  P''(0)=p_2=\frac{\big(2\zeta_{11}\zeta_{21}h(\zeta_{10},0)-\zeta_{22}\big)p_1^2+\rho_2(\zeta_{10})}{\zeta_{21}}
 \]
and this proves $(a)$. Let us turn now to the proof of the assertions in $(b)$. With this aim we note first that the  transition time between $\Pi_1$ and $\Pi_2$ has the following integral expression
 \[
  T(s)=\int_{\xi_1(s)}^{\zeta_1(P(s))}\Omega(x,s)^{\n}f(x,\Omega(x,s))dx.
 \]
By \lemc{L2} we know that $\Omega$ is a $\mathscr C^K$ function such that  $\Omega(x,0)=0$ and $\partial_2\Omega(x,0)=\rho_1(x)$. Hence, the application of \lemc{new-factor1}
shows that $\Omega(x,s)=s(\rho_1(x)+R(x,s))$ for some $\mathscr C^{K-1}$ function~$R$ with $R(x,0)=0$.  Accordingly $T(s)=s^{\n}\tilde T(s)$ with
 \[
  \tilde T(s)\!:=\int_{\xi_1(s)}^{\zeta_1(P(s))}(\rho_1(x)+R(x,s))^{\n}f(x,\Omega(x,s))dx.
 \] 
Then, since $\rho_1$ does not vanish, by a well-known result on the regularity properties of integrals depending on parameters (see~\cite[page 411]{Zorich2}) it follows that $\tilde T$ is $\mathscr C^{K-1}$ as well. Let us compute now $\tilde T(0)$ and $\tilde T'(0)$. This is easy for the first one because $\tilde T(0)=\int_{\xi_{10}}^{\zeta_{10}}\rho_1^{\n}(x)f(x,0)dx.$ Concerning the second one we note that
 \[
  \tilde T'(0)=\rho_1^{\n}(\zeta_{10})f(\zeta_{10},0)\zeta_{11}p_1-\rho_1^{\n}(\xi_{10})f(\xi_{10},0)\xi_{11}
  +\int^{\zeta_{10}}_{\xi_{10}}\!\!\rho_1^{{\n}-1}(x)\left(\frac{1}{2}{\n}\rho_2(x)f(x,0)+\partial_2f(x,0)\rho_1^2(x)\right)dx.
 \]
Here we use that, thanks to \lemc{L2}, $\partial_sR(x,0)=\frac{1}{2}\partial_{22}^2\Omega(x,0)=\frac{1}{2}\rho_2(x).$ Now, taking $\rho_1(\xi_{10})=\xi_{21}$ and $\rho_1(\zeta_{10})=\zeta_{21}p_1$ into account, one can verify that the above expression is equal to the one given in the statement. Hence it only remains to prove the assertions concerning the case ${\n}=0.$ The fact that if ${\n}=0$ then~$T$ is $\mathscr C^{K}$ follows from the regularity properties of integrals depending on parameters that we mention above. With regard to the expression of $T''(0)$ we note that if $\ell=0$ then
 \begin{align*}
  T'(s)=&f\big(\zeta_1(P(s)),\Omega(\zeta_1(P(s)),s)\big)\zeta_1'(P(s))P'(s)-f\big(\xi_1(s),\Omega(\xi_1(s),s)\big)\xi_1'(s)\\
  &+\int_{\xi_1(s)}^{\zeta_1(P(s))}\partial_2f(x,\Omega(x,s))\partial_2\Omega(x,s)dx.
 \end{align*} 
Accordingly, since $\partial_1\Omega(x,0)=0$, $\partial_2\Omega(x,0)=\rho_1(x)$ and $\partial_{22}^2\Omega(x,0)=\rho_2(x)$, some easy computations give
 \begin{align*}
 T''(0)=&\,\partial_1f(\zeta_{10},0)\zeta_{11}^2p_1^2+2\partial_2f(\zeta_{10},0)\rho_1(\zeta_{10})\zeta_{11}p_1
             +f(\zeta_{10},0)\big(\zeta_{12}p_1^2+2\zeta_{11}p_2\big)\\[3pt]
           &-\partial_1f(\xi_{10},0)\xi_{11}^2-2\partial_2f(\xi_{10},0)\rho_1(\xi_{10})\xi_{11}-f(\xi_{10},0)\xi_{12} \\[3pt]
           &+\int^{\zeta_{10}}_{\xi_{10}}\Big(\partial_{22}^2f(x,0)\rho_1^2(x)+\rho_2(x)\partial_2f(x,0)
              \Big)dx. 
 \end{align*} 
Finally the substitution of $\rho_1(\xi_{10})=\xi_{21}$ and $\rho_1(\zeta_{10})=\zeta_{21}p_1$ yields to the expression of $T''(0)$ given in the statement. This concludes the proof of the result. 
 \end{prova}

\section{An incomplete Mellin transform}\label{Mellin}\label{apMellin}

In this appendix we introduce a sort of incomplete Mellin transform that is a key tool for giving a closed expression for the coefficients of the first monomials in the asymptotic expansion of the Dulac map and Dulac time. In short, given $\alpha\in\R\setminus\Z_{\geq 0}$ and a smooth function $f(x)$ on an open interval $I$ that contains $x=0$, we consider the singular scalar differential equation 
\[
 xy'-\alpha y=f(x).
\]
It turns out that this differential equation has for each $\alpha$ a unique solution $y=\hat f(\alpha,x)$ which is smooth on~$I$. As we will see, the fact that $0\in I$ turns out to be crucial for the uniqueness.  The idea is to relate this particular solution with the trajectories of the autonomous planar differential system 
\[
  \left\{\!
   \begin{array}{l}
    \dot x=x, \\[2pt] \dot y=\alpha y+f(x),
   \end{array}
  \right.
\] 
that has a hyperbolic critical point at $(0,-f(0)/\alpha)$ being a saddle for $\alpha<0$ and a focus for $\alpha>0.$ In the saddle case, which is the simplest one, $y=\hat f(\alpha,x)$ is no more than the graph of the stable separatrix. This is in fact the idea in the proof of our next result, which is a little more complicated than it should be because in our applications $f$ depends on parameters and we need good regularity properties of the solution with respect to~$\alpha$ and these parameters as well. For that purpose we apply the so-called center-stable manifold theorem (see for instance \cite[Theorem 1]{Kelley}) but instead one may use the parametrization method for invariant manifolds (see \cite{Fontich03,Fontich05}).

\begin{theo}\label{L8}
Let us consider an open interval $I$ of $\R$ containing $x=0$ and an open subset $U$ of $\R^N$.
\begin{enumerate}[$(a)$]

\item Given $f(x;\nu)\in\mathscr C^{\infty}(I\times U)$, there exits a unique $\hat f(\alpha,x;\nu)\in\mathscr C^{\infty}((\R\setminus\Z_{\ge 0})\times I\times U)$ such that 
\begin{equation}\label{Mellin-eq}
 x\partial_x\hat f({\alpha},x;\nu)-\alpha\hat f({\alpha},x;\nu)=f(x;\nu).
\end{equation}

\item If $x\in I\setminus\{0\}$ then $\partial_x(\hat f({\alpha},x;\nu)|x|^{-\alpha})=f(x;\nu)\frac{|x|^{-\alpha}}{x}$ and, taking any $k\in\Z_{\ge0}$ with $k>\alpha$,
\begin{equation}\label{Mellin-int}
\hat f(\alpha,x;\nu)=
\sum_{i=0}^{k-1}\frac{\partial_x^if(0;\nu)}{i!(i-\alpha)}x^i+|x|^{\alpha}\int_0^x\!\left(f(s;\nu)-T_0^{k-1}f(s;\nu)\right)|s|^{-\alpha}\frac{ds}{s},
\end{equation}
where $T_0^kf(x;\nu)=\sum_{i=0}^{k}\frac{1}{i!}\partial_x^if(0;\nu)x^i$ is the $k$-th degree Taylor polynomial of $f(x;\nu)$ at $x=0$.

\item 
For each $(i_0,x_0,\nu_0)\in\Z_{\ge 0}\times I\times W$ the function $(\alpha,x,\nu)\mapsto(i_0-\alpha)\hat f(\alpha,x;\nu)$ extends $\mathscr C^\infty$ at $(i_0,x_0,\nu_0)$ and, moreover, it tends to $\frac{1}{i_0!}\partial_x^{i_0}f(0;\nu_0)x_0^{i_0}$ as $(\alpha,x,\nu)\to (i_0,x_0,\nu_0).$

\item If $f(x;\nu)$ is analytic on $I\times U$ then $\hat f(\alpha,x;\nu)$ is analytic on $(\R\setminus\Z_{\ge 0})\times I\times U$. Finally,
for each $(\alpha_0,x_0,\nu_0)\in\Z_{\ge 0}\times I\times U$ the function
$(\alpha,x,\nu)\mapsto(\alpha_0-\alpha)\hat f(\alpha,x;\nu)$ extends analytically to $(\alpha_0,x_0,\nu_0)$.
\end{enumerate}
\end{theo}

\begin{prova}
The plan to prove $(a)$ is the following. The uniqueness will be proved firstly. We will show, secondly, the existence for $\alpha<0$ and, thirdly, the existence for $\alpha>0$. 

To prove the uniqueness let us suppose that, for some $\alpha\notin\Z_{\geq 0}$, the differential equation $xy'-\alpha y=f(x;\nu)$ has two solutions, $y=\hat f_1(\alpha,x;\nu)$ and $y=\hat f_2(\alpha,x;\nu)$, that are $\cc^\infty$ on $(\R\setminus\Z_{\ge 0})\times I\times U$. Then $\hat f_1-\hat f_2$ is a smooth function that verifies the homogeneous linear differential equation $xy'-\alpha y=0$ which, in the case that $\alpha\notin\Z_{\ge 0}$, has $y=0$ as unique $\mathscr C^\infty$ solution passing through $x=0.$ Consequently $\hat f_1=\hat f_2,$ as desired.

Let us prove now the existence for the case $\alpha<0.$ To this end, related with the scalar differential equation in \refc{Mellin-eq}, note that the planar vector field 
$x\partial_x+(\alpha y+f(x;\nu))\partial_y$ has, for each fixed $\alpha<0$ and $\nu\in U$, a hyperbolic saddle at $(0,-f(0;\nu)/\alpha)$ with a non-vertical stable separatrix. In order to study its regularity with respect to the parameters we consider the augmented system 
\[
  \left\{\!
   \begin{array}{l}
    \dot x=x, \\[2pt] \dot y=\alpha y+f(x;\nu), \\[2pt] \dot\alpha=0, \\[2pt] \dot\nu=0.
   \end{array}
  \right.
\] 
For each fixed $\alpha_0\in (-\infty,0)$ and $\nu_0\in U,$ the application of \cite[Theorem 1]{Kelley} shows that for every $k\in\N$ there exists a local center-stable manifold $W$ at $(0,-f(0;\nu_0)/\alpha_0,\alpha_0,\nu_0)$ that is written as $y=\hat f_{loc}(\alpha,x;\nu)$ where $\hat f_{loc}$ is a~$\mathscr C^k$ function in a neighbourhood $V$ of $(\alpha_0,0,\nu_0).$ In this context, contrary to what happens in general, it turns out that the center-stable manifold is unique, which implies that  $\hat f_{loc}$ is $\mathscr C^\infty$ (see \cite[p. 165]{Takens}). That being said, we assume without lost of generality that $V$ is a cube with center $(\alpha_0,0,\nu_0)$ and edge length $4\varepsilon$. Then for the points in the strip $\mathcal S=\{(\alpha,x,\nu): x\in I\text{ and }(\alpha,0,\nu)\in V\}$ we define
 \begin{equation}\label{Meq1} 
  \hat f(\alpha,x;\nu)\!:=\left\{
  \begin{array}{ll}
   \displaystyle x^{\alpha}\left(\hat f_{loc}(\alpha,\varepsilon;\nu)\varepsilon^{-\alpha}+\int_{\varepsilon}^xf(s;\nu)s^{-\alpha}\frac{ds}{s}\right) & \text{ if $x\in I\cap (0,+\infty),$} \\ [10pt]
   \hat f_{loc}(\alpha,0;\nu)& \text{ if $x=0,$} \\ [10pt]
   \displaystyle (-x)^{\alpha}\left(\hat f_{loc}(\alpha,-\varepsilon;\nu)\varepsilon^{-\alpha}+\int_{-\varepsilon}^xf(s;\nu)(-s)^{-\alpha}\frac{ds}{s}\right) &  \text{ if $x\in I\cap (-\infty,0),$}
  \end{array}
  \right.
 \end{equation}
which is clearly $\mathscr C^\infty$ on $\mathcal S\setminus\{x=0\}$. An easy computation shows that the above function verifies the scalar differential equation \refc{Mellin-eq} for all $(\alpha,x,\nu)\in\mathcal S$ with $x\neq 0.$ Hence, due to $\hat f(\alpha,\pm\varepsilon;\nu)=\hat f_{loc}(\alpha,\pm\varepsilon;\nu)$, by the existence and uniqueness theorem for solutions of differential equations (see~\cite[Theorem 1.1]{Hale} for instance) we have that $\hat f|_{V}=\hat f_{loc}$ and, consequently, $\hat f\in\mathscr C^{\infty}(\mathcal S)$. On account of the uniqueness of $\hat f$ proved firstly, the arbitrariness of $\alpha_0\in (-\infty,0)$ and $\nu_0\in U$ shows that \refc{Meq1} provides a well defined $\mathscr C^{\infty}$ function $\hat f(\alpha,x;\nu)$ on $(-\infty,0)\times I\times U$. This proves the existence for the case $\alpha<0.$
 
Let us show next the existence for the case $\alpha>0.$ In what follows we shall use the more compact notation $\hat \ell_\alpha(x;\nu)=\hat \ell(\alpha,x;\nu)$ omitting also the dependence on $x$ and $\nu$ when there is no risk of ambiguity. Following this notation, some easy computations show that
\begin{enumerate}[$1.$]
\item If $\ell=g+h$ then $\hat \ell_\alpha=\hat g_\alpha+\hat h_\alpha$, provided that $\hat g_\alpha$ and $\hat h_\alpha$ exist.
\item If $\ell(x;\nu)=\sum_{i=0}^kd_i(\nu)x^i$ and $\alpha\notin\{0,1,2,\ldots,k\}$ 
        then $\hat \ell_\alpha(x;\nu)=\sum_{i=0}^k\frac{d_i(\nu)}{i-\alpha}x^i$.
\item If $\ell(x;\nu)=x^mg(x;\nu)$ with $m>\alpha$ then $\hat \ell_\alpha(x;\nu)=x^m\hat g_{\alpha-m}(x;\nu)$.
\end{enumerate}
That being said, let us fix an arbitrary $m\in \N$ and note that, by applying \lemc{new-factor2}, we can write
 \[
  f(x;\nu)=\sum_{i=0}^{m-1}d_i(\nu)x^i+x^mg(x;\nu),
 \]
with $d_i\in\mathscr C^\infty(U)$ and $g\in\mathscr C^\infty(I\times U).$ On account of this, since we have already proved the existence of~$\hat f_{\alpha}$ for $\alpha<0$, the three properties above imply the existence of $\hat f(\alpha,x;\nu)\in\mathscr C^{\infty}\big(((-\infty,m)\setminus\Z_{\geq 0})\times I\times U\big)$ satisfying~\refc{Mellin-eq}. Finally the arbitrariness of $m\in\N$ and the uniqueness of $\hat f$ proved firstly imply that $\hat f(\alpha,x;\nu)$ is a well defined $\mathscr C^{\infty}$ function on $(\R\setminus\Z_{\ge 0})\times I\times U$ verifying~\refc{Mellin-eq}. This concludes the proof of $(a$). 

Let us prove next the assertions in $(b)$. The fact that the equality $\partial_x(\hat f({\alpha},x;\nu)|x|^{-\alpha})=f(x;\nu)\frac{|x|^{-\alpha}}{x}$ holds for all $x\in I\setminus\{0\}$ follows easily from \refc{Mellin-eq}
by considering the cases $x>0$ and $x<0$ separately. In order to prove \refc{Mellin-int} we note first that, thanks to \lemc{new-factor2}, we can write $f(x;\nu)-T_0^{k-1}f(x;\nu)=x^kg(x;\nu)$ with $g\in\mathscr C^\infty(I\times U).$ Taking this into account and performing the coordinate change $s=tx$ we get  
\[
|x|^\alpha\int_0^x(f(s;\nu)-T_0^{k-1}f(s;\nu))|s|^{-\alpha}\frac{ds}{s}=|x|^\alpha\int_0^xs^{k}g(s;\nu)|s|^{-\alpha}\frac{ds}{s}=x^k\int_0^1t^{k-\alpha}g(tx;\nu)\frac{dt}{t}.
\]
We claim that this is a $\mathscr C^\infty$ function of $(\alpha,x,\nu)\in(-\infty,k)\times I\times U$. To prove this we apply assertions $(i)$, $(c)$ and $(g)$ in \lemc{FLK} to conclude that $(t;\alpha,x,\nu)\mapsto t^{k-\alpha-1}g(tx;\nu)$ belongs to $\F_L^{\infty}((-\infty,k-1-L)\times I\times U)$ for any $L\in\R.$ Consequently, if we fix any $\alpha_0\in (-\infty,k)$ and take $L=\frac{k-\alpha_0}{2}-1$ then for any $x_0\in I$, $\nu_0\in U$, $K\in\Z_{\ge 0}$ and $\nu\in\Z_{\ge 0}^{N+2}$ with $|\nu|\leqslant K$ there exist a compact neighborhood $Q$ of $(\alpha_0,x_0,\nu_0)$ and constants $C,t_0>0$ such that the absolute value of 
\[
 \partial^\nu\big(t^{k-\alpha-1}g(tx;\nu)\big)=\frac{\partial^{|\nu|}(t^{k-\alpha-1}g(tx;\nu))}
  {\partial^{\nu_1}\nu_1\cdots\partial^{\nu_N}\nu_N\partial^{\nu_{N+1}}\alpha\partial^{\nu_{N+2}}x}
\]
is bounded by $Ct^L$ for all $(\alpha,x,\nu)\in Q$ and $t\in (0,t_0)$. It is clear on the other hand that there exists $C'>0$ such that $|\partial^\nu(t^{k-\alpha-1}g(tx;\nu))|\leqslant C'$ for all $(\alpha,x,\nu)\in Q$ and $t\in [t_0,1]$. Accordingly $|\partial^\nu(t^{k-\alpha-1}g(tx;\nu))|$ is bounded by an integrable function of $t\in[0,1]$ not depending on $(\alpha,x,\nu)$. Hence, by applying the Dominated Convergence Theorem (see \cite[Theorem 11.30]{Rudin} and also
\cite[pp. 409--410]{Zorich2}) we can assert that the function $(\alpha,x,\nu)\mapsto \int_0^1t^{k-\alpha}g(tx;\nu)\frac{dt}{t}$ is $\mathscr C^\infty$ on a neighbourhood of $(\alpha_0,x_0,\nu_0)$. This proves the claim and shows in particular that the function on the right hand side of the equality in \refc{Mellin-int} is written as
\[
 \psi(\alpha,x;\nu)\!:=\sum_{i=0}^{k-1}\frac{\partial_x^if(0;\nu)}{i!(i-\alpha)}x^i
 +x^k\int_0^1t^{k-\alpha}g(tx;\nu)\frac{dt}{t}\text{ for all $x\in I\setminus\{0\}$.}
\]
Furthermore, on account of the claim, $\psi\in\cc^\infty\big(((-\infty,k)\setminus\Z_{\geq 0})\times I\times U\big).$ 
On the other hand, by applying the integration by parts formula it follows easily that $x\partial_x\psi-\alpha\psi=f$. 
Consequently
\begin{align}\label{B1eq2}
\hat f(\alpha,x;\nu)&=\sum_{i=0}^{k-1}\frac{\partial_x^if(0;\nu)}{i!(i-\alpha)}x^i
 +x^k\int_0^1t^{k-\alpha}g(tx;\nu)\frac{dt}{t}\\\notag
 &=\sum_{i=0}^{k-1}\frac{\partial_x^if(0;\nu)}{i!(i-\alpha)}x^i+|x|^{\alpha}\int_0^x\!\left(f(s;\nu)-T_0^{k-1}f(s;\nu)\right)|s|^{-\alpha}\frac{ds}{s},
\end{align}
where the first equality is true for all $(\alpha,x,\nu)\in ((-\infty,k)\setminus\Z_{\geq 0})\times I\times U$ by the uniqueness of $\hat f$ and the second one holds only for $x\neq 0$ by the variable change $s=tx.$ This completes the proof of $(b)$. 

In order to prove $(c)$ let us fix $({i_0},x_0,\nu_0)\in\Z_{\geq 0}\times I\times U$ and take any $k\in\Z_{\geq 0}$ such that $k>{i_0}.$ Then the equality in \refc{B1eq2} shows that $(\alpha,x,\nu)\mapsto({i_0}-\alpha)\hat f(\alpha,x;\nu)$ extends $\mathscr C^\infty$ at $({i_0},x_0,\nu_0)$ and, moreover, that it tends to $\frac{1}{i_0!}\partial_x^{i_0}f(0;\nu_0)x_0^{i_0}$ as $(\alpha,x,\nu)\to (i_0,x_0,\nu_0).$

Let us turn finally to the proof of $(d)$, so we assume henceforth that $f(x;\nu)$ is analytic on $I\times U$. Fix any $\alpha_0\in\R\setminus\Z_{\geq 0}$ and $\nu_0\in U$. We claim that the singular differential equation $xy'-\alpha y=f(x;\nu)$ has a solution $y=\hat f_{loc}(\alpha,x;\nu)$ with $\hat f_{loc}(\alpha,0;\nu)=-\frac{1}{\alpha}f(0;\nu)$ that is analytic in a neighbourhood of $(\alpha_0,0,\nu_0)$ inside $(\R\setminus\Z_{\geq 0})\times I\times U.$ 

To prove the claim we consider the holomorphic extension $F(x,\nu)$ of $f(x;\nu)$ in a neighbourhood $\Omega$ of $(0,\nu_0)\in\C^{N+1}$ and for each $i\in\Z_{\geq 0}$ we define $G_i(\alpha,x,\nu)\!:=\frac{\partial_x^iF(0,\nu)}{i!(i-\alpha)}x^i$, which is clearly a holomorphic function on $(\C\setminus\Z_{\geq 0})\times\Omega.$ We will see that
\begin{equation}\label{analitico1}
 S(\alpha,x,\nu)\!:=\sum_{i=0}^{\infty}G_i(\alpha,x,\nu)
\end{equation}
is a holomorphic function in a neighbourhood of $(\alpha_0,0,\nu_0)\in (\C\setminus\Z_{\geq 0})\times\Omega$. To this end we observe that:
\begin{enumerate}[$(i)$] 
\item  By Cauchy's Estimates, see for instance \cite{Rudin}, if $|F(x,\nu)|\leqslant M$ for all $(x,\nu)\in\Omega$ 
          with $|x|<R$ and $|\nu-\nu_0|<\varepsilon$ 
          then $|\partial_x^iF(0,\nu)|\leqslant\frac{i!M}{R^i}.$
\item There exist $\delta_1,\delta_2>0$ small enough such that if $|\alpha-\alpha_0|<\delta_1$ then 
         $|i-\alpha|>\delta_2$ for all $i\in\Z_{\geq 0}$. 
\end{enumerate}
Consequently $|G_i(\alpha,x,\nu)|<\frac{M}{\delta_2}\left(\frac{L}{R}\right)^i$ for all $(\alpha,x,\nu)\in\C^{N+2}$ with $|x|< L<R$, $|\nu-\nu_0|<\varepsilon$ and $|\alpha-\alpha_0|<\delta_1.$ This shows that \refc{analitico1} converges uniformly  in a neighbourhood of $(\alpha_0,0,\nu_0)\in (\C\setminus\Z_{\geq 0})\times\Omega$. On account of this, 
and the fact that $G_i(\alpha,x,\nu)$ is holomorphic on $(\C\setminus\Z_{\geq 0})\times\Omega$ for all $i\geqslant 0$, we can assert (see for instance \cite[Proposition~2]{Malgrange}) that $S(\alpha,x,\nu)$ is holomorphic on $(\C\setminus\Z_{\geq 0})\times\Omega$. We have on the other hand that $x\partial_xS-\alpha S=F$ because, by the uniform convergence again,
 \[
  x\partial_x S(\alpha,x,\nu)-\alpha S(\alpha,x,\nu)=
  x\sum_{i=0}^{\infty}\frac{\partial_x^iF(0,\nu)}{i!(i-\alpha)}ix^{i-1}
  -\alpha\sum_{i=0}^{\infty}\frac{\partial_x^iF(0,\nu)}{i!(i-\alpha)}x^{i}=\sum_{i=0}^{\infty}\frac{\partial_x^iF(0,\nu)}{i!}x^{i}
 =F(x;\nu). 
 \]
Therefore the claim follows taking $\hat f_{loc}(\alpha,x;\nu)$ to be the restriction of $S(\alpha,x;\nu)$ to the real domain. 

Suppose that $\hat f_{loc}(\alpha,x;\nu)$ is analytic in some open cube $V$ with center $(\alpha_0,0,\nu_0)$ and edge length $4\varepsilon$. Then from here we follow exactly the same approach as in the proof of~$(a)$, i.e., we define $\hat f(\alpha,x;\nu)$ in $\mathcal S=\{(\alpha,x,\nu): x\in I\text{ and }(\alpha,0,\nu)\in V\}$ by means of~\refc{Meq1} and it turns out that $\hat f(\alpha,x;\nu)$ is analytic on $\mathcal S\setminus\{x=0\}$. Indeed, this follows from the analyticity of $f(x;\nu)$ and that, on account of the previous claim, $(\alpha,\nu)\mapsto\hat f_{loc}(\alpha,\pm\varepsilon;\nu)$ is analytic at $(\alpha_0,\nu_0)$. Then, exactly as for the regularity assertion in $(a)$, by the existence and uniqueness theorem for solutions of differential equations we have that~$\hat f$ is an analytic function on $\mathcal S$. By the arbitrariness of $\nu_0\in U$ and $\alpha_0\in\R\setminus\Z_{\geq 0},$ this shows that $\hat f(\alpha,x;\nu)$ is analytic on $(\R\setminus\Z_{\ge 0})\times I\times U$. 

In order to prove the second assertion in $(d)$ we fix $\alpha_0\in\Z_{\geq 0}$ and $\nu_0\in U.$ Then the proof of the previous claim shows that $(\alpha,x,\nu)\longmapsto (\alpha-\alpha_0)\hat f(\alpha,x,\nu)$ is analytic at $(\alpha_0,x_0,\nu_0)$ for $x_0=0.$ To prove that this is also true for any $x_0\in I$ we argue exactly as before by using the extension defined in \refc{Meq1} and, for the sake of shortness, it is left to the reader. This concludes the proof of the result.
\end{prova}

\begin{obs} There are some previous results related with the function $\hat f(\alpha,x;\nu)$ defined in \teoc{L8} that should be referred here: 
\begin{enumerate}[$(i)$]

\item Bénoit uses in \cite[p. 106]{Benoit} a transformation $M_\alpha:\C[[t]]\to\C[[t]]$ for every fixed $\alpha\in\R_{>0}\setminus\Z$ defined, for each formal series $f\in\C[[t]]$, by means of the differential equation $-t\frac{d}{dt}M_\alpha(f)+\alpha M_\alpha(f)=f$. Hence, by assertion $(a)$ in \teoc{L8}, if $f\in\R[[t]]$ is convergent then  $M_\alpha(f)=-\hat f(\alpha,t)$. 

\item If $\alpha<0$ then we can take $k=0$ in \refc{Mellin-int} and get that 
\[
 \hat f(\alpha,x)=x^{\alpha}\int_0^x f(s)s^{-\alpha}\frac{ds}{s}\text{ for $x>0.$}
\]  
Therefore if $\alpha>0$ then $\lim_{x\to+\infty}x^{\alpha}\hat f(-\alpha,x)$ coincides with the usual Mellin transform (see~\cite{flajolet}) 
\[
 \mathscr Mf(\alpha)=\int_{0}^\infty f(s)s^\alpha\frac{ds}{s}.
\] 
\item Novikov introduces in \cite{Novikov} a truncated (the author calls it one-sided) Mellin transform as
\[
 u\in L_{loc}^1\big((0,1]\big)\longmapsto\mathscr M_1u(\alpha)\!:=\int_0^1s^{\alpha-1}u(s)ds
\] 
and observe in this regard that $\mathscr M_1u(\alpha)=\hat u({-\alpha},1)$ for $\alpha>0$. 
\end{enumerate}
The formula in \refc{Mellin-int} enables to interpret $\hat f(\alpha,x;\nu)$ as a sort of incomplete (and parametric) version of the Mellin transform of $f(x;\nu)$. As we have seen in the proof of \teoc{L8}, \refc{Mellin-int} extends $\cc^\infty$ to $x=0$ by means of the expression~\refc{B1eq2} taking the $\cc^\infty$ function $g(x;\nu)=\frac{f(x;\nu)-T_0^{k-1}f(x;\nu)}{x^k}$, see \lemc{new-factor2}.
\end{obs}

The proof of the following two results is omitted because it is an easy application of \teoc{L8}.

\begin{cory}\label{B21} 
Consider an open interval $I$ of $\R$ containing $x=0$, an open subset $U$ of $\R^N$ and $\alpha\in\R\setminus\Z_{\geq 0}.$ Then the following hold:
\begin{enumerate}[$(a)$]
\item If $f(x;\nu)=g(x;\nu)+h(x;\nu)$ with $g,h\in \cc^{\infty}(I\times U)$ then 
         $\hat f(\alpha,x;\nu)=\hat g(\alpha,x;\nu)+\hat h(\alpha,x;\nu)$.
\item If $f(x;\nu)=c(\nu)g(x;\nu)$ with $g\in \cc^{\infty}(I\times U)$ and $c\in \cc^{\infty}(U)$ then 
         $\hat f(\alpha,x;\nu)=c(\nu)\hat g(\alpha,x;\nu).$
\item If $f(x;\nu)=x^ng(x;\nu)$ with $g\in \cc^{\infty}(I\times U)$ and $n\in\N$ then 
         $\hat f(\alpha,x;\nu)=x^n\hat g(\alpha-n,x;\nu).$ 
\item If $f(x;\nu)\equiv 1$ then $\hat f(\alpha,x;\nu)\equiv-\frac{1}{\alpha}.$                    
\end{enumerate}
\end{cory}

The next two results are equally valid in the smooth category $\cc^\infty$ and the analytic category $\cc^\omega$. For simplicity in the exposition we write $\cc^\varpi$ with the wild card $\varpi\in\{\infty,\omega\}$.

\begin{cory}\label{B22} 
Let us fix $\varpi\in\{\infty,\omega\}$ and 
consider an open interval $I$ of $\R$ containing $x=0$ and an open subset $U$ of $\R^N$. If $f(x;\nu)\in \cc^{\varpi}(I\times U)$ and $\kappa_1,\kappa_2,\alpha_0\in\R$ verify $\kappa_1\neq 0$ and $i_0\!:=\kappa_1\alpha_0+\kappa_2\in\Z_{\geq 0}$ then, for any $(x_0,\nu_0)\in I\times U$, the function $(\alpha,x,\nu)\mapsto (\alpha_0-\alpha)\hat f(\kappa_1\alpha+\kappa_2,x;\nu)$ extends $\cc^{\varpi}$ at $(\alpha_0,x_0,\nu_0)$  and it tends to $\frac{1}{\kappa_1i_0!}\partial_x^{i_0}f(0;\nu_0)x_0^{i_0}$
as $(\alpha,x,\nu)\to (\alpha_0,x_0,\nu_0)$.
\end{cory}

We conclude the present appendix by proving a technical lemma to be applied for studying the poles of the coefficients obtained in \teoc{A}.

\begin{lem}\label{gorrobis} Let us fix $\varpi\in\{\infty,\omega\}$ and 
consider an open interval $I$ of $\R$ containing $x=0$, an open subset~$U$ of $\R^N$ and $\alpha\in\R\setminus\Z_{\geq 0}.$ Let $M(x;\nu)$ and $A(x;\nu)$ be $\cc^{\varpi}$ functions on $I\times U$ and define 
\[
B(x;\alpha,\nu)\!:=A(x;\nu)\hat M(\alpha,x;\nu),
\]
which is a $\cc^{\varpi}$ function on $I\!\times\!(\R\setminus\Z_{\geq 0})\!\times\! U$ by \teoc{L8}. Finally let us take $i_0,p,q\in\Z$, with $i_0\geqslant 0$ and $q\neq -1$, and set $i_1\!:=qi_0-p$ and $i_2\!:=(q+1)i_0-p$. The following assertions hold:
\begin{enumerate}[$(a)$]

\item If $i_1\geqslant 0$ then, for any $(x_0,\nu_0)\in I\times U$, the function $(\alpha,x,\nu)\mapsto (i_0-\alpha)^2\hat B((q+1)\alpha-p,x;\alpha,\nu)$ extends~$\cc^{\varpi}$ at $(i_0,x_0,\nu_0)$ and it tends to 
\[
 \frac{x_0^{i_2}}{q+1}\frac{M^{(i_0)}(0;\nu_0)}{i_0!}\frac{A^{(i_1)}(0;\nu_0)}{i_1!}\text{ as $(\alpha,x,\nu)\to (i_0,x_0,\nu_0)$}.
\]

\item If $i_1<0$ then, for any $(x_0,\nu_0)\in I\times U$, the function $(\alpha,x,\nu)\mapsto (i_0-\alpha)\hat B((q+1)\alpha-p,x;\alpha,\nu)$ extends~$\cc^{\varpi}$ at $(i_0,x_0,\nu_0)$ and it tends to
\begin{align*}
&\frac{x_0^{i_2}}{(q+1)\,i_2!}\sum_{j=0}^{i_2}{i_2\choose j}\frac{M^{(j)}(0;\nu_0)A^{(i_2-j)}(0;\nu_0)}{j-i_0}+x_0^{i_0}\frac{M^{(i_0)}(0;\nu_0)}{i_0!}\hat A(i_1,x_0;\nu_0)\text{ as $(\alpha,x,\nu)\to (i_0,x_0,\nu_0),$}
\end{align*}
where the summation is zero in the case that $i_2<0.$ 
\end{enumerate}
\end{lem}

\begin{prova}
By applying \lemc{new-factor2} we can write $M(x;\nu)=\sum_{j=0}^{i_0}\frac{M^{(j)}(0;\nu)}{j!}x^j+x^{{i_0}+1}g(x;\nu)$ with $g\in \cc^{\varpi}(I\times U)$. Then the application of \coryc{B21} shows that $\hat M(\alpha,x;\nu)=\sum_{j=0}^{i_0}\frac{M^{(j)}(0;\nu)}{j!(j-\alpha)}x^j+x^{{i_0}+1}\hat g(\alpha-{i_0}-1,x;\nu)$. Consequently, on account of $B(x;\alpha,\nu)\!:=A(x;\nu)\hat M(\alpha,x;\nu)$, we get that
\[
 B(x;\alpha,\nu)=\sum_{j=0}^{i_0}\frac{M^{(j)}(0;\nu)}{j!(j-\alpha)}x^jA(x;\nu)+x^{{i_0}+1}N(x;\alpha,\nu),
 \]
where we set $N(x;\alpha,\nu)\!:=A(x;\nu)\hat g(\alpha-{i_0}-1,x;\nu)$ for shortness. Observe that, since $\hat g(\alpha-{i_0}-1,x;\nu)$ is~$\cc^{\varpi}$ along $\alpha={i_0}$ by \teoc{L8}, so is $N(x;\alpha,\nu)$. Hence, by applying \coryc{B21} again with $\alpha'=(q+1)\alpha-p$ and $\nu'=(\alpha,\nu)$,
\[
 \hat B\big((q+1)\alpha-p,x;\alpha,\nu\big)=\sum_{j=0}^{i_0}\frac{M^{(j)}(0;\nu)}{j!(j-\alpha)}x^j\hat A\big((q+1)\alpha-p-j,x;\nu\big)+x^{{i_0}+1}\hat N\big((q+1)\alpha-p-{i_0}-1,x;\alpha,\nu\big). 
\]
Thus multiplying by $({i_0}-\alpha)^k$ on both sides of the above equality we get
\begin{align}
({i_0}-\alpha)^k\hat B\big((q+1)\alpha-p,x;\alpha,\nu\big)=&\sum_{j=0}^{i_0}\frac{M^{(j)}(0;\nu)}{j!}\frac{({i_0}-\alpha)^k}{j-\alpha}\hat A\big((q+1)\alpha-p-j,x;\nu\big)x^j \notag\\  
&\quad+({i_0}-\alpha)^kx^{{i_0}+1}\hat N\big((q+1)\alpha-p-{i_0}-1,x;\alpha,\nu\big).\label{B3eq0}
\end{align}
In order to prove $(a)$ we set $k=2$ above, so that
\begin{align*}
({i_0}-\alpha)^2\hat B\big((q+1)\alpha-p,x;\alpha,\nu\big)&=
\frac{M^{({i_0})}(0;\nu)}{{i_0}!}({i_0}-\alpha)\hat A\big((q+1)\alpha-p-{i_0},x;\nu\big)x^i\\
&\quad+\sum_{j=0}^{{i_0}-1}\frac{M^{(j)}(0;\nu)}{j!}\frac{({i_0}-\alpha)^2}{j-\alpha}\hat A\big((q+1)\alpha-p-j,x;\nu\big)x^j\\  
&\quad+({i_0}-\alpha)^2x^{{i_0}+1}\hat N\big((q+1)\alpha-p-{i_0}-1,x;\alpha,\nu\big).
\end{align*}
By \coryc{B22} this expressions shows that $(\alpha,x,\nu)\mapsto (i_0-\alpha)^2\hat B((q+1)\alpha-p,x;\alpha,\nu)$ extends~$\cc^{\varpi}$ at $(i_0,x_0,\nu_0)$ for any $(x_0,\nu_0)\in I\times U.$ Furthermore, since 
all the summands except the first one tend to zero as $(\alpha,x,\nu)\to (i_0,x_0,\nu_0)$ by \coryc{B22} again, 
\begin{align}\label{B3eq1}
\lim_{(\alpha,x,\nu)\to (i_0,x_0,\nu_0)}({i_0}-\alpha)^2&\hat B((q+1)\alpha-p,x;\nu)\\\notag
&=\frac{M^{({i_0})}(0;\nu_0)}{{i_0}!}x_0^{i_0}\lim_{(\alpha,x,\nu)\to (i_0,x_0,\nu_0)}({i_0}-\alpha)\hat A\big((q+1)\alpha-p-{i_0},x;\nu\big)
\end{align}
provided that the limit on the right hand side exists. In order to compute it we apply \coryc{B22} once again, with $\kappa_1=q+1$ and $\kappa_2=-p-{i_0}$, to conclude that
\[
\lim_{(\alpha,x,\nu)\to (i_0,x_0,\nu_o)}\hat A\big((q+1)\alpha-p-{i_0},x;\nu\big)=\frac{x_0^{i_1}}{q+1}\frac{A^{(i_1)}(0;\nu_0)}{i_1!},
\]
where we also take the assumption $i_1=q{i_0}-p=\kappa_1{i_0}+\kappa_2\in\Z_{\geq 0}$ into account. Consequently, from \refc{B3eq1},
\[
\lim_{(\alpha,x,\nu)\to (i_0,x_0,\nu_0)}({i_0}-\alpha)^2\hat B((q+1)\alpha-p,x;\nu)
=\frac{x_0^{{i_0}+i_1}}{q+1}\frac{M^{({i_0})}(0;\nu_0)}{{i_0}!}\frac{A^{(i_1)}(0;\nu_0)}{i_1!}
\]
and this proves $(a)$. Let us turn next to the assertion in $(b)$. In this case we set $k=1$ in \refc{B3eq0} to obtain 
\begin{align*}
({i_0}-\alpha)\hat B\big((q+1)\alpha-p,x;\alpha,\nu\big)&=
\frac{M^{({i_0})}(0;\nu)}{{i_0}!}\hat A\big((q+1)\alpha-p-{i_0},x;\nu\big)x^{i_0}\\
&\quad+\sum_{j=0}^{{i_0}-1}\frac{M^{(j)}(0;\nu)}{j!}\frac{{i_0}-\alpha}{j-\alpha}\hat A\big((q+1)\alpha-p-j,x;\nu\big)x^j\\  
&\quad+({i_0}-\alpha)x^{{i_0}+1}\hat N\big((q+1)\alpha-p-{i_0}-1,x;\alpha,\nu\big).
\end{align*}
Note that the last summand on the right hand side is $\cc^{\varpi}$ at $(i_0,x_0,\nu_0)$ by applying \teoc{L8} because $(q+1)\alpha-p-{i_0}-1|_{\alpha={i_0}}=i_1-1<0$ due to the hypothesis $i_1\!:=q{i_0}-p<0$. It shows furthermore that it tends to zero as $(\alpha,x,\nu)\to (i_0,x_0,\nu_0)$. Exactly the same reason shows that the first summand is $\cc^{\varpi}$ at $(i_0,x_0,\nu_0)$ and that it tends to $\frac{M^{(i_0)}(0;\nu_0)}{i_0!}\hat A\big(q{i_0}-p,x_0;\nu_0\big)$ as $(\alpha,x,\nu)\to (i_0,x_0,\nu_0)$. Then, by applying \coryc{B22} with $\kappa_1=q+1$ and $\kappa_2=-p-j$,
the remaining summands on the right hand side also extend~$\cc^{\varpi}$ at $(i_0,x_0,\nu_0)$ and
\begin{align*}
\lim_{(\alpha,x,\nu)\to (i_0,x_0,\nu_0)}&(i_0-\alpha)\hat B((q+1)\alpha-p,x;\nu)\\
&=\frac{1}{q+1}\sum_{j=0}^{i_2}\frac{x_0^{i_2}}{j-{i_0}}\frac{M^{(j)}(0;\nu_0)}{j!}\frac{A^{(i_2-j)}(0;\nu_0)}{(i_2-j)!}+x_0^{i_0}\frac{M^{({i_0})}(0;\nu_0)}{{i_0}!}\hat A(q{i_0}-p,x_0;\nu_0).
\end{align*}
Here we also use that $\kappa_1{i_0}+\kappa_2=(q+1){i_0}-p-j\geqslant 0$ if and only if $j\leqslant (q+1){i_0}-p=:\!i_2$. This proves $(b)$ and concludes the proof of the result. 
\end{prova}

\bibliographystyle{plain}

\end{document}